%% file: article.tex
\documentclass{gtart}
\newcommand{\pathtotrunk}{./}
\input{text/article_preamble.tex}
\input{text/top_matter.tex}

\begin{document}

\begin{abstract}
We construct a new subfactor planar algebra, and as a corollary a new subfactor, with the `extended Haagerup' principal graph pair.   This completes the classification of irreducible amenable subfactors with index in the range $(4,3+\sqrt{3})$, which was initiated by Haagerup in 1993.  We prove that the subfactor planar algebra with these principal graphs is unique.  We give a skein theoretic description, and a description as a subalgebra generated by a certain element in the graph planar algebra of its principal graph. In the skein theoretic description there is an explicit algorithm for evaluating closed diagrams.  This evaluation algorithm is unusual because intermediate steps may increase the number of generators in a diagram.
\end{abstract}

\maketitle

% remove table of contents for submitted version
%\setcounter{tocdepth}{1}
%\tableofcontents

\section{Introduction}
\input{text/intro.tex}

\section{Background}\label{sec:background}
\input{text/background.tex}

\section{Uniqueness and existence}
\label{sec:uniqueness}
\input{text/uniqueness.tex}

%\section{Existence}
%\label{sec:existence}
%\input{text/generator.tex}

\section{The jellyfish algorithm}
\label{sec:algorithm}
\input{text/algorithm.tex}

%\section{Relations from moments}
%\label{sec:relations}
%\input{text/relations.tex}
\section{Relations from moments}
\label{sec:relations}
\input{text/newrelations.tex}
%\section{Relations from moments III}
%\label{sec:relations}
%\input{text/relations3.tex}

\section{Properties of the generator}
\label{sec:propertiesofS}
\input{text/propertiesofS}

\appendix
\input{text/fusion}

% ----------------------------------------------------------------
\hfuzz5pt 
\newcommand{\urlprefix}{}
\bibliographystyle{gtart}
%Included for winedt:
%input "bibliography/bibliography.bib"
\bibliography{../../bibliography/bibliography}
% ----------------------------------------------------------------

This paper is available online at \arxiv{0909.4099}, and at
\url{http://tqft.net/EH}.

% A GTART necessity:
% \Addresses
% ----------------------------------------------------------------
\end{document}

%% file: text/article_preamble.tex
\usepackage{ifpdf}

\ifpdf
\usepackage[pdftex,plainpages=false,hypertexnames=false,pdfpagelabels]{hyperref}
\else
\usepackage[dvips,plainpages=false,hypertexnames=false]{hyperref}
\fi

\usepackage[usenames]{color}

\input{\pathtotrunk preamble.tex}

\ifpdf
\usepackage[pdftex]{graphicx}
\else
\usepackage[dvips]{graphicx}
\fi

\newtheorem{example}{Example}

\usepackage{tikz}
\usetikzlibrary{shapes}
\usetikzlibrary{backgrounds}
\usetikzlibrary{decorations,decorations.pathreplacing}
\input{\pathtotrunk diagrams/tikz/TikzStyles}
\input{\pathtotrunk diagrams/tikz/InnerProducts/STrains.tex}%
\usepackage{microtype}

%% file: preamble.tex
%auto-ignore
%this ensures the arxiv doesn't try to start TeXing here.

\usepackage{amsmath,amssymb,amsfonts}
\usepackage{enumerate}

\usepackage{leftidx,soul,ulem}
\usepackage[section]{placeins}

% ----------------------------------------------------------------
\vfuzz2pt % Don't report over-full v-boxes if over-edge is small
\hfuzz2pt % Don't report over-full h-boxes if over-edge is small
% ----------------------------------------------------------------

% diagrams -------------------------------------------------------
% figures ---------------------------------------------------------
%%% borrowed from Dror's cobordisms paper, use this to include eps or pdf graphics.
\ifpdf
\newcommand{\pathtodiagrams}{\pathtotrunk diagrams/pdf/}
\else
\newcommand{\pathtodiagrams}{\pathtotrunk diagrams/eps/}
\fi

\newcommand{\mathfig}[2]{{\hspace{-3pt}\begin{array}{c}%
  \raisebox{-2.5pt}{\includegraphics[width=#1\textwidth]{\pathtodiagrams #2}}%
\end{array}\hspace{-3pt}}}

\newcommand{\arxiv}[1]{\href{http://arxiv.org/abs/#1}{\tt arXiv:\nolinkurl{#1}}}
\newcommand{\doi}[1]{\href{http://dx.doi.org/#1}{{\tt DOI:#1}}}

\newcommand{\googlebooks}[1]{(preview at \href{http://books.google.com/books?id=#1}{google books})}

\def\RCS$#1: #2 ${\expandafter\def\csname RCS#1\endcsname{#2}}
\RCS$Revision: 1916 $ \RCS$Date: 2010-07-12 11:05:48 -0700 (Mon, 12 Jul 2010) $

\hyphenation{Tem-per-ley}
\hyphenation{Math-e-mat-ic-a}

% THEOREMS -------------------------------------------------------
\theoremstyle{plain}
\newtheorem{prop}{Proposition}[section]

\newtheorem{thm}[prop]{Theorem}
\newtheorem*{thm*}{Theorem}

\newtheorem{lem}[prop]{Lemma}
\newtheorem{defnlem}[prop]{Definition-Lemma}

\newtheorem*{cor*}{Corollary}

\newtheorem{defn}[prop]{Definition}         % numbered definition
\newtheorem*{defn*}{Definition}             % unnumbered definition
\newtheorem*{notation*}{Notation} 

\newtheorem{assumption}[prop]{Assumption}

\newtheorem*{maintheorem}{Theorem \ref{thm:existence}}
\newenvironment{rem}{\noindent\textsl{Remark.}}{}  % perhaps looks better than rem above?
\numberwithin{equation}{section}
%\numberwithin{figure}{section}

% Marginal notes in draft mode -----------------------------------
     % draft mode
     % draft mode
     % draft mode
     % draft mode
     % draft mode
\newcounter{comment}
\newcommand{\noop}[1]{}

% \mathrlap -- a horizontal \smash--------------------------------
% For comparison, the existing overlap macros:
% \def\llap#1{\hbox to 0pt{\hss#1}}
% \def\rlap#1{\hbox to 0pt{#1\hss}}
\def\clap#1{\hbox to 0pt{\hss#1\hss}}

% MATH -----------------------------------------------------------
\newcommand{\Natural}{\mathbb N}
\newcommand{\Integer}{\mathbb Z}
\newcommand{\Rational}{\mathbb Q}
\newcommand{\Real}{\mathbb R}
\newcommand{\Complex}{\mathbb C}

\newcommand{\cC}{\mathcal{C}}

\newcommand{\cI}{\mathcal{I}}

\newcommand{\cN}{\mathcal{N}}

\newcommand{\cP}{\mathcal{P}}
\newcommand{\cQ}{\mathcal{Q}}

\newcommand{\PA}{\mathcal{PA}}

\renewcommand{\imath}{\mathfrak{i}}
\renewcommand{\jmath}{\mathfrak{j}}

\newcommand{\Ha}[1]{\mathcal{H}_{#1}}
\newcommand{\Bi}[1]{\mathcal{B}_{#1}}
\newcommand{\AH}{\mathcal{AH}}

\newcommand{\into}{\hookrightarrow}

\newcommand{\iso}{\cong}

\newcommand{\setcl}[2]{\left\{ \left. #1 \;\right| \; #2 \right\}}

\newcommand{\pairing}[2]{\left\langle#1 ,#2 \right\rangle}

\newcommand{\ppres}{\cQ^{k}}
\newcommand{\pun}{\ppres\!/\!\cN}

\newcommand{\punone}{\cQ^1\!/\!\cN}
\newcommand{\punzero}{\cQ^0\!/\!\cN}

\newcommand{\code}[1]{{\tt #1}}
\newcommand{\MMA}{\code{Math\-e\-mat\-ic\-a} {}}

\newcommand{\directSum}{\oplus}

\newcommand{\DirectSum}{\bigoplus}
\newcommand{\tensor}{\otimes}
\newcommand{\Tensor}{\bigotimes}

\newcommand{\JW}[1]{f^{(#1)}}
\newcommand{\tr}[1]{\text{tr}\left(#1\right)}
\newcommand{\dtr}[1]{\text{tr}_0\left(#1\right)}

% ----------------------------------------------------------------

%% file: diagrams/tikz/TikzStyles.tex
\tikzstyle{shaded}=[fill=red!10!blue!20!gray!50!white]
\tikzstyle{shaded line}=[double=red!10!blue!20!gray!50!white, double distance=2mm, draw=black]
\tikzstyle{unshaded}=[fill=white]
\tikzstyle{unshaded line}=[double=white, double distance=2mm, draw=black]
\tikzstyle{Tbox}=[circle, draw, thick, fill=white, opaque,]
\tikzstyle{empty box}=[circle, draw, thick, fill=white, opaque, inner sep=2mm]
\tikzstyle{background rectangle}= [fill=red!10!blue!20!gray!40!white,rounded corners=2mm] 
\tikzstyle{on}=[very thick, red!50!blue!50!black]
\tikzstyle{off}=[gray]

\tikzstyle{traces}=[scale=.2, inner sep=1mm,baseline]
\tikzstyle{quadratic}=[scale=.23, inner sep=.5mm, baseline]
\tikzstyle{annular}=[scale=.6, inner sep=1mm, baseline]
\tikzstyle{make triple edge size}= [scale=.4, inner sep=1mm,baseline] 
\tikzstyle{icosahedron network}=[scale=.35, inner sep=1mm, baseline]
\tikzstyle{ATLsix}=[scale=.25, baseline]
\tikzstyle{TL12}=[scale=.15,baseline=-0.5ex]
\tikzstyle{TLEG}=[scale=.5,baseline]
\tikzstyle{PAdefn}=[scale=.7,baseline]
\tikzstyle{STrain}=[baseline=0,scale=2]

%% file: diagrams/tikz/InnerProducts/STrains.tex
\newcommand{\upsidedown}[1]{\begin{scope}[y=-1cm] #1 \end{scope}}

\newcommand{\drawS}[3]{%
	\filldraw[fill=white,thick] (#1,#2) ellipse (3mm and 3mm);
	\node at (#1,#2) {\Large $S$};
	\path(#1,#2) ++(#3:0.37) node {$\star$};
}

\newcommand{\RainbowOne}{
	\fill[shaded] (-0.8,0) -- (-0.8,0.6) arc (180:0:0.8) -- (0.8,0) -- (0.2,0) -- (0.2,1) -- (-0.2,1) -- (-0.2,0);
	\draw (-0.8,0) -- (-0.8,0.6) arc (180:0:0.8) -- (0.8,0);
	\draw (-0.2,0) -- (-0.2,1);
	\draw (0.2,0) -- (0.2,1);
	\node at (0.45,0.5) {\footnotesize$2n-1$};
	\drawS{0}{1}{-90}
}

\newcommand{\RainbowTwo}{
	\draw (0,0) -- (0,1);
	\node at (0.15,0.5) {\footnotesize$2n$};
	\drawS{0}{1}{180}
	\filldraw[shaded] (-1,0) -- (-1,0.5) arc (180:0:1) -- (1,0)  -- (0.9,0) -- (0.9,0.5) arc (0:180:0.9) -- (-0.9,0);
}

\newcommand{\JWPlusTwo}{%
	\filldraw[fill=white,thick] (-1,-0.2) rectangle (1,0.2);
	\node at (0,0) {\Large$\JW{2n+2}$};
}

\newcommand{\JWPlusFour}{%
	\filldraw[fill=white,thick] (-1.2,-0.2) rectangle (1.2,0.2);
	\node at (0,0) {\Large$\JW{2n+4}$};
}

\newcommand{\JWwide}[1]{%
	\filldraw[fill=white,thick] (-1.2,-0.2) rectangle (1.2,0.2);
	\node at (0,0) {\Large$\JW{#1}$};
}

\newcommand{\STrainOneOne}{%
	\foreach \x in {-1,0} {
		\node[anchor=south] at (\x+0.5,1) {\footnotesize$n-1$};
		\draw (\x,1) -- (\x + 1, 1);
	}
	\foreach \x in {-1,0,1} {
		\drawS{\x}{1}{90}
	}
}

\newcommand{\STrainOneStrings}{%
	\fill[shaded] (-0.5,0) rectangle (0.5,1);
	\foreach \x in {-0.5,0.5} {
		\draw (\x,1) -- (\x,0);
		\node at (\x+0.25,0.5) {\footnotesize$n+1$};
	}
}

\newcommand{\STrainStrings}[2]{%
	\fill[shaded] (-0.5,0) rectangle (0.5,1);
	\draw (-0.5,1) -- (-0.5,0);
	\node[anchor=west] at (-0.5,0.5) {\footnotesize#1};
	\draw (0.5,1) -- (0.5,0);
	\node[anchor=west] at (0.5,0.5) {\footnotesize#2};
}

\newcommand{\STrainOne}{%
	\node[anchor=south] at (0,1) {\footnotesize$n-1$};
	\draw (-0.5,1) -- (0.5, 1);
	\foreach \x in {-0.5,0.5} {
		\drawS{\x}{1}{90}
	}
}

\newcommand{\STrainOver}{%
        \draw (0.5,1) .. controls (-0.2,1.7) and (-1,1.7)..              
              (-1,1) .. controls (-1,0.5) and (-0.7,0.5) ..  (-0.7,0.2);
        \node at (-1.1,1) {\footnotesize$2$};
}

\newcommand{\STrainThreeStrings}[3]{%
	\fill[shaded] (-1,1) rectangle (1,0);
	\foreach \x in {-1,0,1} {
		\draw (\x,0) -- (\x,1);
	}
        \node [anchor=west] at (-1,0.5) {\footnotesize#1};
        \node [anchor=west] at (0,0.5) {\footnotesize#2};
	\node [anchor=west] at (1,0.5) {\footnotesize#3};
}

%% file: text/top_matter.tex
\title{Constructing the extended Haagerup planar algebra}

\author{Stephen~Bigelow}
\address{
}%
\author{Kim~Morrison}
\address{
}%
\author{Emily~Peters}
\address{
}%
\author{Noah~Snyder}
\address{
}%
\address{%
\rm URLs:\stdspace \tt \url{http://www.math.ucsb.edu/~bigelow/} \url{http://tqft.net/} \\
\tt \url{http://euclid.unh.edu/~eep} \rm and \tt
\url{http://math.berkeley.edu/~nsnyder} \\
\rm Email:  \tt bigelow@math.ucsb.edu, scott@tqft.net, \\ eep@euclid.unh.edu  \rm and \tt nsnyder@math.columbia.edu}

\date{
  First edition: September 21, 2009.
}

\primaryclass{
46L37 %Subfactors and their classification
} \secondaryclass{
18D10 %Monoidal categories (= multiplicative categories), symmetric monoidal categories, braided categories [See also 19D23]	
} \keywords{
  Planar Algebras, Subfactors, Skein Theory, Principal Graphs
}

%% file: text/intro.tex
%!TEX root = ../article.tex
A subfactor is an inclusion $N \subset M$ of von Neumann algebras with trivial center.  The theory of subfactors can be thought of as a nonabelian version of Galois theory, and has had many applications in operator algebras, quantum algebra, and knot theory.  For example, the construction of a new finite depth subfactor, as in this paper, also yields two new fusion categories (by taking the even parts) and a new $3$-dimensional TQFT (via the Ocneanu-Turaev-Viro construction \cite{MR1191386, MR1317353,MR0281657}).

A subfactor $N \subset M$ has three key invariants.  From strongest to weakest, they are: the standard invariant (which captures all information about ``basic'' bimodules over $M$ and $N$), the principal and dual principal graphs (which together describe the fusion rules for these basic bimodules), and the index (which is a real number measuring the ``size'' of the basic bimodules).  We will use the axiomatization of the standard invariant as a subfactor planar algebra, which is due to Jones \cite{math.QA/9909027}.  Other axiomatizations include Ocneanu's paragroups \cite{MR996454} and Popa's $\lambda$-lattices \cite{MR1334479}.  (For readers more familiar with tensor categories, these three approaches are analogous to the diagram calculus \cite{MR0281657, MR1091619, MR1113284}, basic $6j$ symbols \cite[Chapter 5]{MR1292673}, and towers of endomorphism algebras \cite{MR951511}, respectively.) The standard invariant is a complete invariant of amenable subfactors of the hyperfinite $II_1$ factor \cite{MR1055708, MR1278111}.

The index of a subfactor $N \subset M$ must lie in the set $$\setcl{4 \cos^2\left(\frac{\pi}{n}\right)}{n \geq 3} \cup [4,\infty],$$ and all numbers in this set can be realized as the index of a subfactor \cite{MR696688}.  Early work on classifying subfactors of ``small index" concentrated on the case of index less than $4$.  The principal graphs of these subfactors are exactly the Dynkin diagrams $A_n$, $D_{2n}$, $E_6$ and $E_8$.  Furthermore there is exactly one subfactor planar algebra with principal graph $A_n$ or $D_{2n}$ and there is exactly one pair of complex conjugate subfactor planar algebras with principal graph $E_6$ or $E_8$.  (See \cite{MR996454} for the outline of this result, and \cite{MR1193933, MR1145672, MR1313457, MR1308617} for more details.)  The story of the corresponding classification for index equal to $4$ is outlined in \cite[p. 231]{MR1278111}. In this case, the principal graph must be an affine Dynkin diagram. For some principal graphs there are multiple non-conjugate subfactors with the same principal graph, which are distinguished by homological data.

The classification of subfactors of ``small index" greater than $4$ was initiated by Haagerup \cite{MR1317352}.  His main result is a list of all possible pairs of principal graphs of irreducible subfactors of index larger than $4$ but smaller than $3+\sqrt{3}$. Here we begin to see subfactors whose principal graph is different from its dual principal graph.  If $\Gamma$ refers to a pair of principal graphs and we need to refer to one individually, we will use the notation $\Gamma^p$ and $\Gamma^d$.  Any subfactor $N\subset M$ has a dual given by the basic construction $M\subset M_1$.  Taking duals reverses the shading on the planar algebra, switches the principal and dual principal graphs, and preserves index.  Haagerup's list is as follows (we list each pair once).
\newcommand{\ind}{\operatorname{index}}
\begin{itemize}
\item 
 $(A_\infty,A_\infty)$, %which for every index greater than or equal to $4$ \cite{MR1198815},
\item
the infinite family
$$\left\{\Ha{n}  = \left(\mathfig{0.32}{graphs/haagerup-n}, \mathfig{0.32}{graphs/dual-haagerup-n}\right)\right\}_{n\in \Natural}, $$
which has $\ind(\Ha{0})=\frac{5+\sqrt{13}}{2}$,
$\ind(\Ha{1})$ the largest root of $x^3 - 8 x^2 + 17 x -5$,
and $\ind(\Ha{n})$ monotonically increasing with $n$,
converging to the real root of $x^3-6x^2+8x-4$,
(thus $\ind(\Ha{0}) \approx 4.30278$,
      $\ind(\Ha{1}) \approx 4.37720$,
 and $\lim_n \ind(\Ha{n}) \approx 4.38298$),
\item
the infinite family
$$\left\{\Bi{n}  =  \left(\mathfig{0.32}{graphs/hexagon-n}, \mathfig{0.32}{graphs/dual-hexagon-n}\right)\right\}_{n \in \Natural},$$
which has $\ind(\Bi{0}) = \frac{7+\sqrt{5}}{2}$,
and $\ind(\Bi{n})$ monotonically increasing in $n$,
converging to the real root of $x^3-8x^2+19x-16$,
(thus $\ind(\Bi{0}) \approx 4.61803$, and $\lim_n \ind(\Bi{n}) \approx 4.65897$),
\item one more pair of graphs, $$\AH = \left(\mathfig{0.32}{graphs/HA}, \mathfig{0.29}{graphs/dual-HA}\right),$$
which has index $\frac{5+\sqrt{17}}{2} \approx 4.56155$.
\end{itemize}

Haagerup's paper announces this result up to index $3+\sqrt{3}\approx 4.73205$, but only proves it up to index $3+\sqrt{2} \approx 4.41421$; this includes all of the graphs $\Ha{n}$, but none of the graphs $\Bi{n}$ or $\AH$.  Haagerup's proof of the full result has not yet appeared.   In work in progress, Jones, Morrison, Penneys, Peters, and Snyder have independently confirmed his result (following Haagerup's outline except at one point using a result from \cite{quadratic}), and have extended his techniques to give a partial result up to index $5$ (see \cite{index5-part1,index5-part2,index5-part3, index5-part4}).  In this paper, we will only rely on the part of Haagerup's classification that has appeared in print.

Haagerup's original result did not specify which of the possible principal graphs are actually realized.
Considerable progress has since been made in this direction.
Asaeda and Haagerup \cite{MR1686551} proved the existence and uniqueness
of a subfactor planar algebra whose principal graphs are $\Ha{0}$ (called the {\em Haagerup subfactor}),
and a subfactor planar algebra for $\AH$ (called the {\em Asaeda-Haagerup subfactor}).
Izumi \cite{MR1832764} gave an alternate construction of the Haagerup subfactor.
Bisch \cite{MR1625762} showed none of the graphs $\Bi{n}$ can be principal graphs because they give inconsistent fusion rules.
Asaeda \cite{MR2307421}  and Asaeda-Yasuda \cite{MR2472028} proved that $\Ha{n}$ is not a principal graph for $n \ge 2$. To do this, they showed that the index is not a cyclotomic integer, and then appealed to
a result of Etingof, Nikshych and Ostrik \cite{MR2183279}, which in turn is proved by reduction to the case of modular categories, where it was proved in the context of rational conformal field theories by Coste--Gannon \cite{MR1266785} using a result of de Boere--Goeree \cite{MR1120140}.

The main result of our paper is
\begin{maintheorem}
There is a subfactor planar algebra with principal graphs $\Ha{1}$.
\end{maintheorem}
In addition, we prove in Theorem \ref{thm:uniqueness} that this planar algebra is the only one with these principal graphs.
This result completes the classification of all subfactor planar algebras up to index $3+\sqrt{3}$:
\begin{cor*}[\cite{MR1317352}, \cite{MR1686551}, \cite{MR1625762}, \cite{MR2472028}, and Theorem \ref{thm:existence}]
The only irreducible subfactor planar algebras with index in the range $(4, 3+\sqrt{3})$ are
\begin{itemize}
\item the non-amenable Temperley-Lieb planar algebra at every index in this range, with principal graphs $(A_\infty, A_\infty)$,
\item the Haagerup planar algebra with principal graphs $\Ha{0}$, and its dual,
\item the Haagerup-Asaeda planar algebra with principal graphs $\AH$, and its dual, and
\item the extended Haagerup planar algebra with principal graphs $\Ha{1}$, and its dual.
\end{itemize}
\end{cor*}
By Popa's classification \cite{MR1055708} the latter three pairs can each be realized uniquely as the standard invariant of a subfactor of the hyperfinite $II_1$ factor.  This gives a complete classification of amenable subfactors of the hyperfinite $II_1$ factor with index between $(4, 3+\sqrt{3})$.  The non-amenable case remains open because it is unknown for which indices Temperley-Lieb can be realized as the standard invariant of the hyperfinite $II_1$ factor, nor in how many ways it can be realized (see \cite{MR1159284, MR1293872} for some work in this direction).   Furthermore, there remain many interesting questions about small index subfactors of arbitrary factors.

It was already expected that the extended Haagerup subfactor should exist, thanks to approximate numerical evidence coming from computations by Ikeda \cite{MR1633929}. We note that although our construction relies on a computation of the traces of a few  large matrices, this computation consists of exact arithmetic in a number field, and is a very different calculation from the one Ikeda did numerically.

The search for small index subfactors
has so far produced the three pairs of ``sporadic'' examples: the
Haagerup, Asaeda-Haagerup and extended Haagerup subfactors.
These are some of the very few known subfactors that do not seem to fit into the frameworks of groups,
quantum groups, or conformal field theory \cite{MR2468378}. (See also a generalization of the Haagerup subfactor due to Izumi \cite[Example 7.2]{MR1832764}).
You might think of them as analogs of the exceptional simple Lie algebras,
or of the sporadic finite simple groups. (Without a good extension theory, it is not yet clear what ``simple'' should mean in this context.)

In this paper,
we study the extended Haagerup planar algebra.
We construct the extended Haagerup planar algebra by
locating it inside the graph planar algebra \cite{MR1865703} of its principal graph.
By a result of Jones--Penneys \cite{gpa} (generalized in \cite{gpa2}) every subfactor planar algebra occurs in this way.
We find the right planar subalgebra by following a recipe outlined by Jones \cite{MR1865703, quadratic} and further developed by Peters \cite{0902.1294}, who applied it to the Haagerup planar algebra.

We also give a presentation of the extended Haagerup planar algebra
using a single planar generator and explicit relations.
We prove that the subalgebra of the graph planar algebra contains an element also satisfying these relations. This is convenient because different properties
become more apparent in different descriptions of the planar algebra.
For example,
the subalgebra of the graph planar algebra is clearly non-trivial,
which would be difficult to prove directly from the generators and relations.
In the other direction,
in \S \ref{sec:relations} we prove that
our relations result in a space of closed diagrams
that is at most one dimensional,
which would be difficult to prove in the graph planar algebra setting.

In \S \ref{sec:background} we recall the definitions of planar algebras and graph planar algebras \cite{math.QA/9909027, MR1865703}.  We also set some notation for the graph planar algebra of $\Ha{1}^p$.  In \S \ref{sec:uniqueness} we prove our two main theorems, Theorems \ref{thm:uniqueness} and \ref{thm:existence}.
Theorem \ref{thm:uniqueness}, the uniqueness theorem, says
that for each $k$ there is at most one subfactor planar algebra with principal graphs $\Ha{k}$.  Furthermore we give a skein theoretic description by generators and relations of the unique candidate planar algebra.  Theorem \ref{thm:existence}, the existence theorem, constructs a subfactor planar algebra with principal graphs $\Ha{1}$ by realizing the skein theoretic planar algebra as a subalgebra of the graph planar algebra. Proofs of several key results needed for the main existence and uniqueness arguments are deferred to \S \ref{sec:algorithm}, \S \ref{sec:relations}, and \S\ref{sec:propertiesofS}.  In particular, \S \ref{sec:algorithm} describes an evaluation algorithm that uses the skein theory to evaluate any closed diagram (Theorem \ref{thm:zeroboxonedim}).  This is crucial to our proofs of both existence and uniqueness and may be of broader interest in quantum topology.  This section can be read independently of the rest of the paper.  
Section  \ref{sec:relations}  consists of calculations of inner products using generators and relations. Section \ref{sec:propertiesofS} gives the description of  the generator of our subfactor planar algebra inside the graph planar algebra  and verifies its properties. Appendix \ref{sec:fusion} gives the tensor product rules for the two fusion categories associated to the extended Haagerup subfactor.

Part of this work was done while Stephen Bigelow and Emily Peters were visiting the University of Melbourne.  Kim Morrison was at Microsoft Station Q and the Miller Institute for Basic Research during this work. Emily Peters  was supported in part by NSF Grant DMS0401734 and a fellowship from Soroptimist International and Noah Snyder was supported in part by RTG grant DMS-0354321 and in part by an NSF Postdoctoral Fellowship. We would like to thank Vaughan Jones for many useful discussions, and Yossi Farjoun for lessons on Newton's method.

%% file: text/background.tex
\subsection{Planar algebras}

Planar algebras were defined in \cite{math.QA/9909027} and \cite{MR1865703}.
More general definitions have since appeared elsewhere,
but we only need the original notion of a \emph{shaded planar algebra},
which we sketch here.  For further details see  \cite[\S 2]{MR1865703}, \cite[\S 0]{math.QA/9909027}, or \cite{MR1957084}.
%(The original definition is phrased in terms of colored operads, but essentially the same).

\begin{defn}
A {\em (shaded) planar tangle} 
has an outer disk,
a finite number of inner disks,
and a finite number of non-intersecting strings.
A string can be either a closed loop
or an edge with endpoints on boundary circles.
We require that
there be an even number of endpoints on each boundary circle,
and a checkerboard shading of
the regions in the complement of the interior disks.
We further require that there be a marked point on the boundary of each disk, and that the inner disks are ordered.

Two planar tangles are considered equal if they are isotopic (not necessarily rel boundary).

Here is an example of a planar tangle.
$$%
%\beginpgfgraphicnamed{\pathtotrunk diagrams/tikz/#1-external}%
\input{\pathtotrunk diagrams/tikz/SampleTangle.tex}%
%\endpgfgraphicnamed
$$

Planar tangles can be composed by placing one planar tangle inside an interior disk of another,
lining up the marked points,
and connecting endpoints of strands.
The numbers of endpoints and the shadings must match up appropriately.  This composition turns the collection of planar tangles into a colored operad.
\end{defn}

\begin{defn}
A {\em (shaded) planar algebra} consists of
\begin{itemize}
\item A family of vector spaces $\{V_{(n,\pm,)}\}_{n\in \Natural}$, called the positive and negative $n$-box spaces.
\item For each planar tangle, a multilinear map $V_{n_1, \pm_1} \otimes \ldots \otimes V_{n_k, \pm_k} \rightarrow V_{n_0, \pm_0}$ where $n_i$ is half the number of endpoints on the $i$th interior boundary circle, $n_0$ is half the number of endpoints on the outer boundary circle, and the signs $\pm$ are positive (respectively negative) when the marked point on the corresponding boundary circle is in an unshaded region (respectively shaded region).
\end{itemize}

For example, the planar tangle above gives a map $$V_{1,+} \otimes V_{2,+} \otimes V_{2,-} \rightarrow V_{3,+}.$$

The linear map associated to a `radial' tangle (with one inner disc, radial strings, and matching marked points) must be the identity.
We require that the action of planar tangles be compatible with composition of planar tangles.  In other words, composition of planar tangles must correspond to the obvious composition of multilinear maps.  In operadic language this says that a planar algebra is an algebra over the operad of planar tangles.
\end{defn}

We will refer to an element of $V_{n,\pm}$ (and specifically $V_{n,+}$, unless otherwise stated) as an ``$n$-box.''

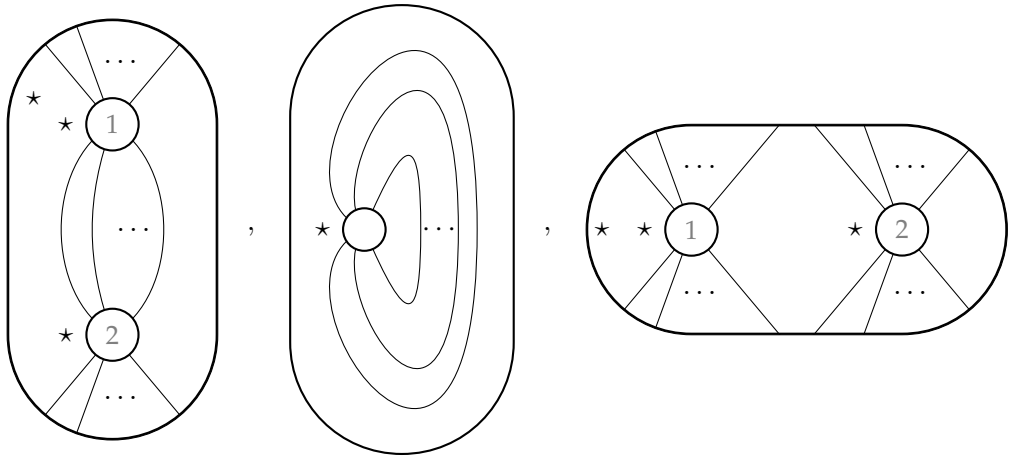
\begin{figure}[ht]
$$%
%\beginpgfgraphicnamed{\pathtotrunk diagrams/tikz/#1-external}%
\input{\pathtotrunk diagrams/tikz/multiplication.tex}%
%\endpgfgraphicnamed
 \quad , \quad %
%\beginpgfgraphicnamed{\pathtotrunk diagrams/tikz/#1-external}%
\input{\pathtotrunk diagrams/tikz/trace.tex}%
%\endpgfgraphicnamed
 \quad , \quad %
%\beginpgfgraphicnamed{\pathtotrunk diagrams/tikz/#1-external}%
\input{\pathtotrunk diagrams/tikz/tensor.tex}%
%\endpgfgraphicnamed
$$
\caption
  {The multiplication, trace, and tensor product tangles.}
\label{fig:timestracetensor}
\end{figure}

We make frequent use of three families of planar tangles
called multiplication, trace, and tensor product, which are shown in Figure \ref{fig:timestracetensor}.
``Multiplication'' gives an associative product
$V_{n,\pm} \otimes V_{n,\pm} \rightarrow V_{n,\pm}$.
``Trace'' gives a map $V_{n,\pm} \rightarrow V_{0,\pm}$.
``Tensor product'' gives an associative product
$V_{m,\pm} \otimes V_{n,\pm} \rightarrow V_{m+n,\pm}$ if $m$ is even,
or $V_{m,\pm} \otimes V_{n,\mp} \rightarrow V_{m+n,\pm}$ if $m$ is odd.

The (shaded or unshaded) empty diagrams can be thought of as elements of $V_{0,\pm}$, since the `empty tangle' induces a map from the empty tensor product $\Complex $ to the space $V_{0,\pm}$.
If the space $V_{0,\pm}$ is one dimensional
then we can identify it with $\Complex$
by sending the empty diagram to one.
In many other cases, we can make do with the following.

\begin{defn}
A {\em partition function} is a pair of linear maps 
$$Z_{\pm}:V_{0,\pm} \rightarrow \Complex$$
that send the empty diagrams to $1$.  

In a planar algebra with a partition function,
let
$$\operatorname{tr} \co V_{n,\pm} \to \Complex$$
denote the composition of the trace tangle with $Z$.
\end{defn}

Sometimes we will need to refer simply to the action of the trace tangle, which we denote $\operatorname{tr}_0 : V_{n,\pm} \to V_{0,\pm}$.

Notice that the above trace tangle is the ``right trace.''  There is also a ``left trace'' where all the strands are connected around the left side.

\begin{defn}
A planar algebra with a partition function can be:
\begin{itemize}
\item \emph{Positive definite:} There is an antilinear adjoint operation $*$ on each $V_{n,\pm}$,
      compatible with the adjoint operation on planar tangles given by reflection.
      The sesquilinear form
      $\langle x,y \rangle = \tr{x y^*}$ is positive definite.
\item \emph{Spherical:} The left trace
      $$\operatorname{tr}_l \co V_{1,\pm} \to V_{0,\mp} \xrightarrow{Z_\mp} \mathbb{C}$$
      and the right trace
      $$\operatorname{tr}_r \co V_{1,\pm} \to V_{0,\pm} \xrightarrow{Z_\pm} \mathbb{C}$$
      are equal.
\end{itemize}
\end{defn}

The spherical property implies that the left and right traces are equal on every $V_{n,\pm}$. Since every planar algebra we consider is spherical, we will usually ignore the distinction between left and right trace.

\begin{defn}
A {\em subfactor planar algebra} is
a positive definite spherical planar algebra
such that
$\dim{V_{0,+}}=\dim{V_{0,-}}=1$ and $\dim{V_{n,\pm}} < \infty$.
\end{defn}

As a consequence of being spherical and having $1$-dimensional $0$-box spaces, subfactor planar algebras always have a well-defined modulus, as described below.

\begin{defn}
We say that the planar algebra has {\em modulus} $d$
if the following relations hold.
$$\begin{tikzpicture}[baseline=.3cm]
		\draw[thick] (.5,.5) circle (.5);
		\filldraw[ shaded] (.5,.5) circle (.3);
	\end{tikzpicture}
	=d \cdot 
	\begin{tikzpicture}[baseline=.3cm]
		\draw[thick] (.5,.5) circle (.5);
	\end{tikzpicture} , \quad
	\begin{tikzpicture}[baseline=.3cm]
		\draw[thick, shaded] (.5,.5) circle (.5);
		\filldraw[fill=white] (.5,.5) circle (.3);
	\end{tikzpicture}
	=d \cdot 
	\begin{tikzpicture}[baseline=.3cm]
		\draw[thick, shaded] (.5,.5) circle (.5);
	\end{tikzpicture}.$$
\end{defn}

The {\em principal graphs} of a subfactor encode the fusion rules for the basic bimodules ${}_N M_M$ and ${}_M M_N$.  The vertices of the principal graph are the isomorphism classes of simple $N$-$N$ and $N$-$M$ bimodules that occur in tensor products of the basic bimodules. The edges give the decompositions of tensor products of simple bimodules with the basic bimodule. The dual principal graph encodes similar information, but for $M$-$N$ and $M$-$M$ bimodules. This definition is due to Connes \cite{MR1303779} and Ocneanu \cite{MR996454}, together with later work of Jones \cite{bimodules} that lets you replace the original Hilbert space bimodules with algebraic bimodules, for example,  ${}_M L^2(M)_N$ with ${}_M M_N$.

In the language of planar algebras, the basic bimodule is a single strand, and the isomorphism classes of simple bimodules are equivalence classes of irreducible projections in the box spaces.  For a more detailed description, see \cite[\S 4.1]{math/0808.0764}.

\begin{defn} 
A subfactor planar algebra $\cP$ is \emph{irreducible} if $\dim{\cP_{1,+}} =1$.

A subfactor planar algebra has \emph{finite depth} if it has finitely many isomorphism classes of irreducible projections, that is, finitely many vertices in the principal graphs.  \cite{MR996454, MR1055708}
\end{defn}

The following is well-known, and combines the results of \cite{MR1055708} and  \cite{math.QA/9909027}.

\begin{thm}
\label{thm:subfactors}
Finite depth finite index subfactors of the hyperfinite $II_1$ factor are in one-to-one correspondence with finite depth subfactor planar algebras.
Irreducible subfactors (those for which $M$ is irreducible as an $N$-$M$ bimodule) correspond to irreducible subfactor planar algebras.
\end{thm}

\begin{proof}
Suppose we are given a subfactor $N \subset M$.  The corresponding planar algebra is constructed as follows.
Let $\cC$ be the $2$-category of all bimodules that appear
in the decomposition of some tensor product of alternating copies of
$M$ as a $N$-$M$ bimodule and $M$ as a $M$-$M$ bimodule.
These are $A$-$B$ bimodules for $A,B \in \{M,N\}$,
and form the $1$-morphisms of $\cC$.
Composition of $1$-morphisms is given by tensor product.
The $2$-morphisms of $\cC$ are the intertwiners.

We can define duals in $\cC$ by taking the contragradient bimodule, which interchanges $M$ as an $N$-$M$ bimodule with $M$ as an $M$-$N$ bimodule.
Now define the associated planar algebra by
$$\cP_{n,\pm} =
 \operatorname{End}_{\cC}\left(\widehat{\Tensor}_n M^{\pm}\right).$$
Here $M^\pm$ means $M$ or $M^*$,
and $\widehat{\Tensor}_n M$
means $M \tensor M^* \tensor M \tensor \cdots \tensor M^{\pm}$.
The action of tangles is via the usual interpretation of string diagrams
as $2$-morphisms in a $2$-category \cite{MR1113284},
with critical points interpreted as evaluation and coevaluation maps.

The difficult direction is to recover a subfactor from a planar algebra.  The proof of this result was given in \cite{MR1055708}.  However in that paper, Popa uses towers of commutants instead of tensor products of bimodules,
and $\lambda$-lattices instead of planar algebras.  See \cite{math.QA/9909027} to translate
from $\lambda$-lattices into planar algebras.  See \cite{MR1424954} and \cite{bimodules} to translate from towers of commutants into tensor products of bimodules.
\end{proof}

\begin{rem}
The above theorem says that a certain kind of subfactor is completely characterized by its  representation theory (that is, the bimodules). This can be thought of as a subfactor version of
the Dopplicher-Roberts theorem \cite{MR1010160}, or more generally, of Tannaka-Krein type theorems \cite{MR1173027}.
\end{rem}

In general, given any extremal finite index subfactor of a $II_1$ factor, the standard invariant is a subfactor planar algebra.  Several other reconstruction results have been proved.  Popa extended his results on finite depth subfactors to amenable subfactors of the hyperfinite $II_1$ factor in \cite{MR1278111}.  The general situation for non-amenable subfactors of the hyperfinite $II_1$ factor is more complicated: some subfactor planar algebras cannot be realized at all (an unpublished result of Popa's, see \cite{MR1159284}), while others can be realized by a continuous family of different subfactors \cite{MR2314611}.    Furthermore, once you move beyond the hyperfinite $II_1$ factor there are many new questions.  On the one hand any subfactor planar algebra comes from a (canonically constructed, but not necessarily unique) subfactor of the free group factor $L(F_{\infty})$ \cite{MR2051399, 0807.4146, 0807.3704, 0712.2904}, while on the other hand there exist factors for which only the trivial planar algebra can be realized as the standard invariant of a subfactor \cite{MR2504433}.

\subsection{Temperley-Lieb}

Everyone's favorite example of a planar algebra is Temperley-Lieb.  It was defined (as an algebra) in \cite{MR0498284}, applied to subfactor theory in \cite{MR696688}, and formulated diagrammatically in \cite{MR899057}. 
% \cite{MR1280463}.
The vector space $TL_{n,\pm}$ is spanned by non-crossing pairings of $2n$ points (where the $\pm$ depends on whether the marked point is in a shaded region or unshaded region).
These pairings are drawn as intervals in a disc, starting from a marked point on the boundary.
For example,
$$TL_{3,+}=\text{span} \left\{ %
%\beginpgfgraphicnamed{\pathtotrunk diagrams/tikz/#1-external}%
\input{\pathtotrunk diagrams/tikz/basis.tex}%
%\endpgfgraphicnamed
 \right\}.$$

Planar tangles act on Temperley-Lieb elements ``diagrammatically:''
the inputs are inserted into the inner disk, strings are smoothed out,
and each loop is discarded in exchange for a factor of $d \in \Complex$.
For example,
$$%
%\beginpgfgraphicnamed{\pathtotrunk diagrams/tikz/#1-external}%
\input{\pathtotrunk diagrams/tikz/TLEG.tex}%
%\endpgfgraphicnamed
.$$

If $d \in \Real$ we can introduce an antilinear involution $*$ by reflecting diagrams.
If $d \geq 2$ then Temperley-Lieb is a subfactor planar algebra.
If $d= 2 \cos \frac{\pi}{n}$ for $n=3,4,5, \ldots$
then we can take a certain quotient to obtain a subfactor planar algebra.
The irreducible projections in the Temperley-Lieb planar algebra
are the Jones-Wenzl idempotents \cite{MR866500, MR873400}.

\begin{defn}\label{defn:JW}
The {\em Jones-Wenzl idempotent} $f^{(n)} \in TL_{n,\pm}$ is characterized by
\begin{align*}
f^{(n)} & \neq 0 \\
f^{(n)} f^{(n)} & = f^{(n)} \\
e_i f^{(n)} = f^{(n)} e_i & = 0 \qquad \text{ for } i=1,2, \ldots , n-1 
\end{align*}
where $e_1,\dots,e_{n-1}$ are the Jones projections
$$e_1 = \frac{1}{d} \; %
%\beginpgfgraphicnamed{\pathtotrunk diagrams/tikz/#1-external}%
\input{\pathtotrunk diagrams/tikz/TLn/e1.tex}%
%\endpgfgraphicnamed
, \quad
  e_2 = \frac{1}{d} \; %
%\beginpgfgraphicnamed{\pathtotrunk diagrams/tikz/#1-external}%
\input{\pathtotrunk diagrams/tikz/TLn/e2.tex}%
%\endpgfgraphicnamed
, \quad \dots \; , \quad
  e_{n-1} = \frac{1}{d} \; %
%\beginpgfgraphicnamed{\pathtotrunk diagrams/tikz/#1-external}%
\input{\pathtotrunk diagrams/tikz/TLn/en-1.tex}%
%\endpgfgraphicnamed
.$$
\end{defn}

The following gives a recursive definition of the Jones-Wenzl idempotents.  We should mention that the following pictures are drawn using rectangles instead of disks, and the marked points are assumed to be on the left side of the rectangles (including the implicit bounding rectangles).  The quantum integers $[n] = \frac{q^n-q^{-n}}{q-q^{-1}}$ which appear below are specialized at a value of $q$ such that $[2]=d$.  

\begin{lem}[\cite{MR1446615}]
\label{lem:basicrecursion}
\begin{align*}
\begin{tikzpicture}[baseline=0,scale=0.3]
\foreach \x in {-3,...,3} { \draw (\x,-3.5) -- (\x,3.5); }
\filldraw[thick,fill=white] (-4,-1.5) rectangle (4,1.5);
\node at (0,0) {$\JW{k}$};
\end{tikzpicture}
= 
\begin{tikzpicture}[baseline=0,scale=0.3]
\foreach \x in {-3,...,2} { \draw (\x,-3.5) -- (\x,3.5); }
\draw (4,-3.5) -- (4,3.5);
\filldraw[thick,fill=white] (-4,-1.5) rectangle (3,1.5);
\node at (-0.5,0) {$\JW{k-1}$};
\end{tikzpicture}
+
\frac{1}{[k]}
\sum_{a = 1}^{k-1} (-1)^{a+k} [a] \;
\begin{tikzpicture}[baseline=0,scale=0.3]
\draw[decorate,decoration={brace,raise=2pt}] (-4,3.5) -- (-2,3.5);
\node at (-3,4.6) {$a$};
\foreach \x in {-4,-3} { \draw (\x+1,0) -- (\x,3.5); }
\draw (-1,3.5) arc (360:180:0.5);
\foreach \x in {0,1} { \draw (\x-1,0) -- (\x,3.5); }
\draw (1,0)--(1.5,1.75) .. controls (2,3.5) and (3,3.5) .. (3,1.5)--(3,-3.5);
\foreach \x in {-3,...,1} { \draw (\x,-3.5) -- (\x,0); }
\filldraw[thick,fill=white] (-4,-1.5) rectangle (2,1.5);
\node at (-1,0) {$\JW{k-1}$};
\end{tikzpicture}.
\end{align*}
\end{lem}

\begin{proof}
It is straightforward to check that the right hand side of this equation satisfies the characterizing relations of Definition \ref{defn:JW}.
\end{proof}

We will also use the following more symmetrical version of the lemma.

\begin{lem}
\label{lem:yetanother}
\begin{align*}
\begin{tikzpicture}[baseline=0,scale=0.3]
\foreach \x in {-3,...,3} { \draw (\x,-3.5) -- (\x,3.5); }
\filldraw[thick,fill=white] (-4,-1.5) rectangle (4,1.5);
\node at (0,0) {$\JW{k}$};
\end{tikzpicture}
= 
\begin{tikzpicture}[baseline=0,scale=0.3]
\foreach \x in {-3,...,2} { \draw (\x,-3.5) -- (\x,3.5); }
\draw (4,-3.5) -- (4,3.5);
\filldraw[thick,fill=white] (-4,-1.5) rectangle (3,1.5);
\node at (-0.5,0) {$\JW{k-1}$};
\end{tikzpicture}
+
\frac{1}{[k][k-1]}
\sum_{a,b = 1}^{k-1} (-1)^{a+b+1} [a][b] \;
\begin{tikzpicture}[baseline=0,scale=0.3]
\draw[decorate,decoration={brace,raise=2pt}] (-4,3.5) -- (-2,3.5);
\node at (-3,4.6) {$a$};
\draw[decorate,decoration={brace,mirror,raise=2pt}] (-4,-3.5) -- (-1,-3.5);
\node at (-2.5,-4.85) {$b$};
\foreach \x in {-4,-3} { \draw (\x,-3.5) -- (\x+1,0) -- (\x,3.5); }
\draw (-2,-3.5) -- (-1,0);
\draw (-1,-3.5) arc (180:0:0.5);
\draw (-1,3.5) arc (360:180:0.5);
\draw (0,3.5) -- (-1,0);
\foreach \x in {1,2} { \draw (\x,-3.5) -- (\x-1,0) -- (\x,3.5); }
\filldraw[thick,fill=white] (-4,-1.5) rectangle (2,1.5);
\node at (-1,0) {$\JW{k-2}$};
\end{tikzpicture}.
\end{align*}
\end{lem}

\begin{proof}
First apply Lemma~\ref{lem:basicrecursion},
and then apply the vertically reflected version of that lemma
to all but the first term in the resulting expression.
\end{proof}

We sometimes consider
the complete expansion of a Jones-Wenzl idempotent
into a linear combination of Temperley-Lieb diagrams.
Suppose $\beta$ is a $k$-strand Tem\-per\-ley-Lieb diagram.
Let $\mathrm{Coeff}_{f^{(k)}}(\beta)$ denote
the coefficient of $\beta$ in the expansion of $f^{(k)}$.
Thus
$$f^{(k)} = \sum \mathrm{Coeff}_{f^{(k)}}(\beta) \beta,$$
where the sum is over all $k$-strand Temperley-Lieb diagrams $\beta$.

We will frequently state values of $\mathrm{Coeff}_{f^{(k)}}(\beta)$
without giving the details of how they are computed.  A convenient formula for these values is given by Frenkel and Khovanov in \cite{MR1446615}. See also \cite{morrison} for a detailed exposition, including a helpful example at the end of \S 4 of that paper.  A related formula was announced by Ocneanu \cite{MR1907188}, and special cases of this were proved by Reznikoff \cite{MR2375712}.
%--------------------------------------

\subsection{The graph planar algebra}\label{sec:gpa}
In this section we define the {\em graph planar algebra} $GPA(G)$ of a bipartite graph $G$ with a chosen base point, and recall some of its basic properties \cite{MR1865703}. Except in degenerate cases, this fails to be a subfactor  planar algebra because $\dim{GPA(G)_{0,+}}$ and $\dim{GPA(G)_{0,-}}$ are greater than $1$.  However, after specifying a certain partition function, all the other axioms for a subfactor planar algebra hold.

The box space $GPA(G)_{n,\pm}$ is the space of functionals on the set of loops on $G$ that have length $2n$ and are based at an even vertex in the case of $GPA(G)_{n,+}$, or an odd vertex in the case of $GPA(G)_{n,-}$.

Suppose $T$ is a planar tangle with $k$ inner disks, 
and $f_1, \ldots, f_k$ are functionals
in the appropriate spaces $GPA(G)_{n_i,\pm_i}$.  
Then we will define $T(f_1,\dots,f_k)$ as a certain ``weighted state sum.''

A {\em state} on $T$ is an assignment of vertices of $G$ to regions of $T$ and edges of $G$ to strings of $T$, such that unshaded regions are assigned even vertices, shaded regions are assigned odd vertices, and the edge assigned to the string between two regions goes between the vertices assigned to those regions.  In particular, a state for any graph is uniquely specified by giving only the assignment of edges, and a state for a simply laced graph is specified by giving only the assignment of vertices.  Since all the graphs we consider are simply laced, we typically specify states by giving the assignment of vertices.  The inner boundaries $\partial_i(\sigma)$ and outer boundary $\partial_0(\sigma)$ of a state 
are the loops obtained by reading the edges assigned to strings clockwise around the corresponding disk.  

We define $T(f_1, \ldots , f_k)$ by describing its value on a loop $\ell$.
This is given by the following weighted state sum.
$$T(f_1, \ldots , f_k)(\ell)
=\sum_\sigma c(T,\sigma) \cdot \prod_{i=1, \ldots, k} f_i (\partial_i(\sigma)).
$$
Here,
the sum is over all states $\sigma$ on $G$
such that $\partial_0(\sigma) = \ell$,
and the weight $c(T,\sigma)$ is defined below.

To specify the weight $c(T, \sigma)$, it helps to draw $T$ in a certain standard form. Each disk is drawn as a rectangle, with the same number of strands meeting the top and bottom edges, no strands meeting the side edges, and the starred region on the left side. The strands are drawn smoothly, with a finite number of local maxima and minima. Then
$$c(T,\sigma) = \prod_{t \in E(T)} \sqrt\frac{d_{\sigma(t_{\text{convex}})}}{d_{\sigma(t_{\text{concave}})}},$$
where $E(T)$ is the set of local maxima and minima of the strands of $T$, $d_{v}$ is the Perron-Frobenius dimension of the vertex $v$, and
$t_{\text{convex}}$ and $t_{\text{concave}}$ are respectively the regions on the convex and concave sides of $t$.  The Perron-Frobenius dimension of a vertex is the corresponding entry in the Perron-Frobenius eigenvector of the adjacency matrix.  This is the largest-eigenvalue eigenvector, normalized so the Perron-Frobenius dimension of the base point is $1$, and its entries are strictly positive.

It is now easy to check that this planar algebra has modulus $d$, the Perron-Frobenius dimension of $G$.

\begin{example}\label{eg:rotation}
Fix $G$, a simply laced graph.
Consider 
$$\rho_8 = %
%\beginpgfgraphicnamed{\pathtotrunk diagrams/tikz/#1-external}%
\input{\pathtotrunk diagrams/tikz/rho8.tex}%
%\endpgfgraphicnamed
,$$
the ``two-click'' rotation on $8$-boxes, already drawn in standard form, and a loop 
$\gamma = \gamma_1 \gamma_2 \ldots \gamma_{16} \gamma_1$.
Then
\begin{align*}
\rho(f)(\gamma)
& = %
%\beginpgfgraphicnamed{\pathtotrunk diagrams/tikz/#1-external}%
\input{\pathtotrunk diagrams/tikz/rho8flabelled.tex}%
%\endpgfgraphicnamed
 \\
& = \sqrt \frac{d_{\gamma_3}d_{\gamma_{11}}}{d_{\gamma_1}d_{\gamma_9}} 
f(\gamma_{3}  \ldots \gamma_{16} \gamma_{1} \gamma_{2} \gamma_{3}).
\end{align*}
\end{example}

It is a general fact about the  Perron-Frobenius dimensions of bipartite graphs that $\sum_{\text{even vertices } v} d_v^2 = \sum_{\text{odd vertices } v} d_v^2$. Call this number $\mathcal{I}$, the \emph{global index}.
We use the partition function
\begin{align*}
Z:& GPA_{0,\pm}  \rightarrow  \mathbb{C} \\
&f   \mapsto \sum_v f(v) \frac{d_v^2}{\mathcal{I}}
\end{align*}
and the involution $*$ given by  reversing loops:
 $$f^*(\gamma_1 \gamma_2 \gamma_3 \ldots \gamma_n \gamma_1) :=
     f( \gamma_1 \gamma_n \ldots \gamma_3 \gamma_2 \gamma_1).$$  

 \begin{prop}\label{GPAsphericalposdef}
 For any bipartite graph $G$ with base point the planar algebra with partition function and involution $(GPA(G), Z, *)$ is a spherical positive definite planar algebra
whose modulus is the Perron-Frobenius eigenvalue for $G$.
 \end{prop}
 
 \begin{proof}
This is due to \cite{MR1865703}, but we recall the easy details here.
 The inner product is positive definite, because the basis of Kronecker-delta functionals on loops
$ \{ \delta_{\gamma} \}_{\gamma \in \Gamma_{2k}} $  is an orthogonal basis and $\langle \delta_\gamma, \delta_\gamma \rangle = \frac{d_{\gamma_1} d_{\gamma_{k+1}}}{\cI} > 0 $.
 
 Sphericality is a straightforward computation:
 \begin{align*}
 	\operatorname{tr}_l (X) =
	Z \left( \begin{tikzpicture}[baseline=.4cm]
		\fill[shaded] (-.7,-.7) rectangle (1.2,1.7);
		\filldraw[fill=white] (.5,1) arc (0:180:.5cm) -- (-.5,0) arc (-180:0:.5cm) arc (-90:90:.5cm);
		\draw[thick] (.5,.5) circle (.5);
		\node at (.5,.5) {$X$};
	\end{tikzpicture} \right)
	= & Z \left( \sum_{\substack{\text{edges $e$ from even}\\\text{to odd vertices}}}  X(e)  \frac{d_{s(e)}}{d_{t(e)}} \cdot \delta_{t(e)} \right) \\
	= &\sum_{\substack{\text{edges $e$ from even}\\\text{to odd vertices}}}  X(e)  d_{s(e)} \frac{d_{t(e)}}{\mathcal{I}}  \\
\intertext{and }
 	\operatorname{tr}_r (X)  =
	Z \left( \begin{tikzpicture}[baseline=.4cm]
		\node at (.5,.5) {$X$};
		\filldraw[shaded] (.5,1) arc (180:0:.5cm) -- (1.5,0) arc (-0:-180:.5cm) arc (-90:90:.5cm);
		\draw[thick] (.5,.5) circle (.5);
	\end{tikzpicture} \right)
	= & Z \left( \sum_{\substack{\text{edges $e$ from even}\\\text{to odd vertices}}}  X(e)  \frac{d_{t(e)}}{d_{s(e)}} \cdot \delta_{s(e)} \right) \\
	= &\sum_{\substack{\text{edges $e$ from even}\\\text{to odd vertices}}}  X(e)  d_{t(e)} \frac{d_{s(e)}}{\mathcal{I}},
 \end{align*}
 where $s(e)$ and $t(e)$ are respectively the even and odd vertices of the edge $e$.
 \end{proof}
 
The main reason for interest in graph planar algebras is the following result from \cite{gpa,gpa2}.
 
\begin{thm}
\label{thm:subalgebra}
Given a finite depth subfactor planar algebra $\cP$ with principal graph $\Gamma$ there is an injective map of planar algebras $$\cP \into GPA(\Gamma).$$
\end{thm}

This theorem assures us that if we believe in the existence of the extended Haagerup subfactor, and have enough perseverance, we will inevitably find it as a subalgebra of the graph planar algebra. Indeed, this paper is the result of such perseverance. On the other hand, nothing is this paper logically depends on the above theorem.

 \subsection{Notation for $\Ha{k}$}\label{sec:names}
Let $d_k$ be the Perron-Frobenius dimension of the graphs $\Ha{k}$. For the Haagerup subfactor, we have $d_0 = \sqrt{\frac{5+\sqrt{13}}{2}} \approx 2.07431$. For the extended Haagerup subfactor,
$d_1$ is the largest root of the polynomial $x^6 - 8 x^4 + 17 x^2 -5$, $$d_1 = \sqrt{\frac{8}{3}+\frac{1}{3} \sqrt[3]{\frac{13}{2} \left(-5-3 i \sqrt{3}\right)}+\frac{1}{3} \sqrt[3]{\frac{13}{2} \left(-5+3 i \sqrt{3}\right)}},$$ approximately $2.09218$.  
 
Throughout, if $d$ is the modulus of a planar algebra, we let $q$ be a solution to $q+q^{-1} = d$, and use the quantum integers $$[n] = \frac{q^{n}-q^{-n}}{q-q^{-1}}.$$
 
By $\Ha{k}^p$ we mean the first graph in the pair of principal graphs $\Ha{k}$.
 When we talk about loops or paths on $\Ha{1}^p$ it is useful to have names
  for the vertices and arms.
$$\begin{tikzpicture}[baseline=-1mm, scale=.7]
	\filldraw (0,0) node [above] {$v_0$} circle (1mm);
	\filldraw (1,0) node [above] {$w_0$}  circle (1mm);
	\filldraw (2,0) node [above] {$x_0$}  circle (1mm);
	\filldraw (3,0) node [above] {$y_0$}  circle (1mm);
	\filldraw (4,0) node [above] {$z_0$}  circle (1mm);
	\filldraw (5,0) node [above] {$a_0$}  circle (1mm);
	\filldraw (6,0) node [above] {$b_0$}  circle (1mm);
	\filldraw (7,0) node [above] {$c$}  circle (1mm);

	\filldraw (7.7,.7) node [above] {$b_2$}  circle (1mm);
	\filldraw (8.7,.7) node [above] {$a_2$}  circle (1mm);
	\filldraw (9.7,.7) node [above] {$z_2$}  circle (1mm);
	\filldraw (7.7,-.7) node [above] {$b_1$}  circle (1mm);
	\filldraw (8.7,-.7) node [above] {$a_1$}  circle (1mm);
	\filldraw (9.7,-.7) node [above] {$z_1$}  circle (1mm);
	
	\draw (0,0) -- (7,0);
	\draw (7,0)--(7.7,.7)--(9.7,.7);
	\draw (7,0)--(7.7,-.7)--(9.7,-.7);
	
	\node at (-2,0) {arm $0 \; \rightarrow$};
	\node at (12,.7) {$\leftarrow$ arm $2$};
	\node at (12,-.7) {$\leftarrow $ arm $1$};
\end{tikzpicture}$$

\begin{lem}
\label{lem:quantum-identity}
If $[2]=d_k$, then $[3][4k+4] = [4k+8]$.
\end{lem}
\begin{proof}
The dimensions of the three vertices on an arm of $\Ha{k}^p$, counting from the branch, are
\begin{align*}
\dim{b_1} & = \frac{[4k+5]}{2}, \\
\dim{a_1} & = \frac{[4k+6]-[4k+4]}{2}, \\
\intertext{and}
\dim{z_1} & = \frac{[4k+7] - [4k+5] - [4k+3]}{2}.
\end{align*}
The condition $[2] \dim{z_1} = \dim{a_1}$ easily gives the desired formula.
\end{proof}

%% file: diagrams/tikz/SampleTangle.tex
\begin{tikzpicture}[scale=.65]
	\clip (0,0) circle (3cm);
	
	\begin{scope}[shift=(10:1cm)]	
		\draw[shaded] (0,0)--(0:6cm)--(90:6cm)--(0,0);	
		\draw[shaded] (0,0) .. controls ++(180:2cm) and ++(-90:2cm) .. (0,0);
	\end{scope}
	
	\draw[shaded] (-150:1cm) -- (120:4cm) -- (180:4cm) -- (-150:1cm);
	\draw[shaded] (-150:1cm) -- (-120:4cm) -- (-60:4cm) -- (-150:1cm);
	
	\begin{scope}[shift=(10:1cm)]	
		\node at (0,0) [Tbox, inner sep=1mm] {\small{\textcolor{gray}{2}}};
		\node at (90:1.5cm) [Tbox, inner sep=1mm] {\small{\textcolor{gray}{1}}};
		\node at (-45:.7cm) {$\star$};
		\node at (120:1.6cm) {$\star$};
	\end{scope}
	\node at (-150:1cm) [Tbox, inner sep=1 mm] {\small{\textcolor{gray}{3}}};
	\node at (-120:1.6cm) {$\star$};
	\node at (-30:2.7cm) {$\star$};
	
	\draw[very thick] (0,0) circle (3cm);
\end{tikzpicture}

%% file: diagrams/tikz/multiplication.tex
 \begin{tikzpicture}[PAdefn]
	\clip [draw] (2,2) arc (0:180:2cm) -- (-2,-2) arc (-180:0:2cm) -- (2,2);
	
	%first draw the lines
	\draw (0,2) .. controls ++(-150:1.5cm) and ++(150:1.5cm) .. (0,-2) .. controls ++(110:1.5cm) and ++(-110:1.5cm) .. (0,2);

	\draw (0,2) .. controls ++(-30:1.5cm) and ++(30:1.5cm) .. (0,-2);
	
	\draw (0,2) -- +(110:3cm) -- +(130:3cm) -- (0,2);
	\draw (0,2) -- ++(50:3cm);
		
	\draw (0,-2) -- ++(-50:3cm);
	\draw (0,-2) -- +(-110:3cm) -- +(-130:3cm) -- (0,-2);
	
	%decorate with dots, labels and stars
	\node at (.4,0) {$\ldots$};
	\node at (.2,3.2) {$\dots$};
	\node at (.2,-3.2) {$\dots$};

	\node at (0,2) [Tbox,inner sep=1.4mm] (A) {\small{\textcolor{gray}{1}}};
	\node at (0,-2) [Tbox,inner sep=1.4mm] (B) {\small{\textcolor{gray}{2}}};
	\node at (A.180) [left] {$\star$};
	\node at (B.180) [left] {$\star$};
	\node at (-1.5,2.5)  {$\star$};	
	
	%redraw boundary
	\draw[ultra thick] (2,2) arc (0:180:2cm) -- (-2,-2) arc (-180:0:2cm) -- (2,2);
\end{tikzpicture}

%% file: diagrams/tikz/trace.tex
\begin{tikzpicture}[scale=.25,baseline]
	\clip (8,6) arc (0:180:6cm) -- (-4,-6) arc (-180:0:6cm) -- (8,6);

%	\filldraw[shaded] (0,0) .. controls ++(-157:3cm) and ++(157:3cm) .. (0,0) .. controls ++(112:2cm) and ++(-112:2cm) .. (0,0);
%	\filldraw[shaded] (0,0) .. controls ++(-22:3cm) and ++(22:3cm) .. (0,0) .. controls ++(67:2cm) and ++(-67:2cm) .. (0,0);
%	\draw (0,0) .. controls ++(22:3cm) and ++(90:4cm) .. (3,0) .. controls ++(-90:4cm) and ++(-22:3cm) .. 
	\draw (0,0) .. controls ++(-67:6cm) and ++(-90:5cm) .. (3,0) .. controls ++(90:5cm) and ++(67:6cm) .. (0,0);
	\draw (0,0) .. controls ++(130:4cm) and ++(90:15cm) .. (5,0) .. controls ++(-90:15cm) and ++(-130:4cm) .. (0,0) .. controls ++(-157:7cm) and ++(-90:20cm) .. (6,0) .. controls ++(90:20cm) and ++(157:7cm) .. (0,0);

	\node at (0,0) [Tbox,inner sep=2mm] (T1) {};
	\node at (T1.180) [left] {$\star$};
	\node at (4.1,0) {$\cdots$};
	
	\draw[ultra thick] (8,6) arc (0:180:6cm) -- (-4,-6) arc (-180:0:6cm) -- (8,6);
\end{tikzpicture}

%% file: diagrams/tikz/tensor.tex
 \begin{tikzpicture}[PAdefn]
	\clip [draw] (2,2) arc (90:-90:2cm) -- (-2,-2) arc (-90:-270:2cm) -- (2,2);
	
	%first draw the lines	
	\draw (2,0) -- +(110:3cm) -- +(130:3cm) -- (2,0);
	\draw (2,0) -- ++(50:3cm);
	\draw (2,0) -- ++(-50:3cm);
	\draw (2,0) -- +(-110:3cm) -- +(-130:3cm) -- (2,0);
		
	\draw (-2,0) -- +(110:3cm) -- +(130:3cm) -- (-2,0);
	\draw (-2,0) -- ++(50:3cm);
	\draw (-2,0) -- ++(-50:3cm);
	\draw (-2,0) -- +(-110:3cm) -- +(-130:3cm) -- (-2,0);
	
	%decorate with dots, labels and stars
	\node at (-1.8,1.2) {$\dots$};
	\node at (-1.8,-1.2) {$\dots$};
	\node at (2.2,1.2) {$\dots$};
	\node at (2.2,-1.2) {$\dots$};

	\node at (2,0) [Tbox,inner sep=1.4mm] (A) {\small{\textcolor{gray}{2}}};
	\node at (-2,0) [Tbox,inner sep=1.4mm] (B) {\small{\textcolor{gray}{1}}};
	\node at (A.180) [left] {$\star$};
	\node at (B.180) [left] {$\star$};
	\node at (-3.7,0)  {$\star$};	
	
	%redraw boundary
	\draw[ultra thick] (2,2) arc (90:-90:2cm) -- (-2,-2) arc (-90:-270:2cm) -- (2,2);
\end{tikzpicture}

%% file: diagrams/tikz/basis.tex
\begin{tikzpicture}[TLEG]
	\filldraw[shaded]  (30:1cm) arc (30:90:1cm) arc (30:-150:5mm) arc (150:210:1cm) -- cycle; 
	\filldraw[shaded]  (0,-1) arc (-90:-30:1cm) arc (30:210:5mm);
	\draw[thick] (0,0) circle (1cm);
	\node at (120:1.3cm) {$\star$};
\end{tikzpicture},
\begin{tikzpicture}[TLEG]
	\filldraw[shaded]  (0,1) arc (90:30:1cm) arc (90:270:5mm) arc (-30:-90:1cm) -- cycle; 
	\filldraw[shaded]  (-1,0) arc (180:210:1cm) arc (-90:90:5mm) arc (150:180:1cm);
	\draw[thick] (0,0) circle (1cm);
	\node at (120:1.3cm) {$\star$};
\end{tikzpicture},
\begin{tikzpicture}[TLEG, rotate=180]
	\filldraw[shaded]  (150:1cm) arc (150:90:1cm) arc (-210:-30:5mm) arc (30:-30:1cm)--cycle;
	\filldraw[shaded] (-90:1cm) arc (-90:-150:1cm) arc (150:-30:5mm);
	\draw[thick] (0,0) circle (1cm);
	\node at (-60:1.3cm) {$\star$};
\end{tikzpicture},
\begin{tikzpicture}[TLEG]
	\filldraw[shaded]  (0,-1) arc (-90:-30:1cm) arc (30:210:5mm);
	\filldraw[shaded]  (-1,0) arc (180:210:1cm) arc (-90:90:5mm) arc (150:180:1cm);
	\filldraw[shaded] (90:1cm) arc (90:30:1cm) arc (-30:-210:5mm);
	\draw[thick] (0,0) circle (1cm);
	\node at (120:1.3cm) {$\star$};
\end{tikzpicture},
\begin{tikzpicture}[TLEG]
	\filldraw[shaded]  (90:1cm) arc (30:-150:5mm) arc (150:210:1cm) arc (150:-30:5mm) arc (-90:-30:1cm) arc (-90:-270:5mm) arc (30:90:1cm);
	\draw[thick] (0,0) circle (1cm);
	\node at (120:1.3cm) {$\star$};
\end{tikzpicture}

%% file: diagrams/tikz/TLEG.tex
\begin{tikzpicture}[TLEG]
	\filldraw[shaded] (90:2.3cm) -- (-90:2.3cm) arc (-90:90:2.3cm);
	\filldraw[unshaded] (0:1cm) circle (.5cm);
	\filldraw[shaded] (180:1cm) circle (.5cm);
	\node[Tbox,inner sep=3.5mm] at (0,0) {};
	\draw[thick] (0,0) circle (2.3cm);
	\node at (120:1.3cm) {$\star$};
	\node at (120:2.6cm) {$\star$};
\end{tikzpicture}
\left(
\begin{tikzpicture}[TLEG]
	\filldraw[shaded]  (0,1) arc (90:30:1cm) arc (90:270:5mm) arc (-30:-90:1cm) -- cycle; 
	\filldraw[shaded]  (-1,0) arc (180:210:1cm) arc (-90:90:5mm) arc (150:180:1cm);
	\draw[thick] (0,0) circle (1cm);
	\node at (120:1.3cm) {$\star$};
\end{tikzpicture}
\right)
=
\begin{tikzpicture}[TLEG]
	\filldraw[shaded] (90:2.3cm) -- (-90:2.3cm) arc (-90:90:2.3cm);
	\filldraw[unshaded] (0:1cm) circle (.5cm);
	\filldraw[shaded] (180:1cm) circle (.5cm);
	\draw[thick] (0,0) circle (2.3cm);
	\node at (120:2.6cm) {$\star$};
\end{tikzpicture}
=d^2
\begin{tikzpicture}[TLEG]
	\draw[thick]  (-2.3,0)arc (180:-180: 2.3cm);
	\filldraw[shaded]  (0,2.3)--(0,-2.3) arc (-90:90:2.3cm);
	\node at (120:2.6cm) {$\star$};
\end{tikzpicture}

%% file: diagrams/tikz/TLn/e1.tex
\begin{tikzpicture}[TL12]
        \draw (1cm,1.5cm) arc (-180:0:.5cm);
        \draw (1cm,-1.5cm) arc (180:0:.5cm);
        \draw (3cm, -1.5cm)--(3cm, 1.5cm);
        \node at (5cm,0) {\tiny{$\dots$}};
        \draw (7cm,-1.5cm)--(7cm,1.5cm);
\end{tikzpicture}

%% file: diagrams/tikz/TLn/e2.tex
\begin{tikzpicture}[TL12]
        \draw (1cm,-1.5cm)--(1cm,1.5cm);
        \draw (2cm,1.5cm) arc (-180:0:.5cm);
        \draw (2cm,-1.5cm) arc (180:0:.5cm);
        \draw (4cm, -1.5cm)--(4cm, 1.5cm);
        \node at (6cm,0) {\tiny{$\dots$}};
        \draw (8cm,-1.5cm)--(8cm,1.5cm);
\end{tikzpicture}

%% file: diagrams/tikz/TLn/en-1.tex
\begin{tikzpicture}[TL12]
        \draw (1cm, -1.5cm)--(1cm, 1.5cm);
        \node at (3cm,0) {\tiny{$\dots$}};
        \draw (5cm, -1.5cm)--(5cm, 1.5cm);
        \draw (6cm,1.5cm) arc (-180:0:.5cm);
        \draw (6cm,-1.5cm) arc (180:0:.5cm);
\end{tikzpicture}

%% file: diagrams/tikz/rho8.tex
\begin{tikzpicture}[scale=.4, baseline]
	\draw[shaded] (0,-2) arc (0:-180:.5cm) -- (-1,2) .. controls ++(90:2cm) and ++(-90:2cm) .. (1,6) -- (0,6) .. controls ++(-90:2cm) and ++(90:2cm) .. (-2,2) -- (-2,-2) arc (-180:0:1.5cm);
	\foreach \x in {0,2,4} \draw[shaded] (\x cm,2 cm) .. controls ++(90:2cm) and ++(-90:2cm) .. (\x cm + 2 cm,6 cm) -- (\x cm + 3 cm, 6 cm) .. controls ++(-90:2cm) and ++(90:2cm) .. (\x cm + 1 cm,2 cm);

\begin{scope}[rotate=180, xshift=-7cm]
		\draw[shaded] (0,-2) arc (0:-180:.5cm) -- (-1,2) .. controls ++(90:2cm) and ++(-90:2cm) .. (1,6) -- (0,6) .. controls ++(-90:2cm) and ++(90:2cm) .. (-2,2) -- (-2,-2) arc (-180:0:1.5cm);
		\foreach \x in {0,2,4} \draw[shaded] (\x cm,2 cm) .. controls ++(90:2cm) and ++(-90:2cm) .. (\x cm + 2 cm,6 cm) -- (\x cm + 3 cm, 6 cm) .. controls ++(-90:2cm) and ++(90:2cm) .. (\x cm + 1 cm,2 cm);
\end{scope}

	\draw[thick] (-.5,-2) rectangle (7.5,2);
	\draw[thick] (-4,-6) rectangle (11,6);
\end{tikzpicture}

%% file: diagrams/tikz/rho8flabelled.tex
\begin{tikzpicture}[scale=.4, baseline]
	\draw[shaded] (0,-2) arc (0:-180:.5cm) -- (-1,2) .. controls ++(90:2cm) and ++(-90:2cm) .. (1,6) -- (0,6) .. controls ++(-90:2cm) and ++(90:2cm) .. (-2,2) -- (-2,-2) arc (-180:0:1.5cm);
	\foreach \x in {0,2,4} \draw[shaded] (\x cm,2 cm) .. controls ++(90:2cm) and ++(-90:2cm) .. (\x cm + 2 cm,6 cm) -- (\x cm + 3 cm, 6 cm) .. controls ++(-90:2cm) and ++(90:2cm) .. (\x cm + 1 cm,2 cm);

\begin{scope}[rotate=180, xshift=-7cm]
		\draw[shaded] (0,-2) arc (0:-180:.5cm) -- (-1,2) .. controls ++(90:2cm) and ++(-90:2cm) .. (1,6) -- (0,6) .. controls ++(-90:2cm) and ++(90:2cm) .. (-2,2) -- (-2,-2) arc (-180:0:1.5cm);
		\foreach \x in {0,2,4} \draw[shaded] (\x cm,2 cm) .. controls ++(90:2cm) and ++(-90:2cm) .. (\x cm + 2 cm,6 cm) -- (\x cm + 3 cm, 6 cm) .. controls ++(-90:2cm) and ++(90:2cm) .. (\x cm + 1 cm,2 cm);
\end{scope}

	\draw[thick] (-.5,-2) rectangle (7.5,2);
	\draw[thick] (-4,-6) rectangle (11,6);

	\node at (3.5,0) {$f$};
	
	%labels
	\foreach \x in {1,...,8} \node at (\x cm - 1.5 cm, 6cm) [above] {\tiny{$\gamma_{\x}$}};
	\foreach \x in {9,...,16} \node at (16.5 cm - \x cm, -6 cm) [below] {\tiny {$\gamma_{\x}$}};
		
\end{tikzpicture}

%% file: text/uniqueness.tex
\subsection{Uniqueness}

The goal of this section is to prove that 
there is at most one subfactor planar algebra with principal graphs $\Ha{k}$.
To prove this, we will give a skein theoretic description of a planar algebra $\pun$ (which is not necessarily a subfactor planar algebra).
We then prove in Theorem~\ref{thm:uniqueness} that any subfactor planar algebra with principal graphs $\Ha{k}$ must be isomorphic to $\pun$.

\begin{defn}\label{defn:uncappable}
 We say that a $n$-box $S$ is {\em uncappable}  if $\epsilon_i(S)=0$ for all $i=1,\ldots, 2n$ where
 $$\epsilon_1 = %
%\beginpgfgraphicnamed{\pathtotrunk diagrams/tikz/#1-external}%
\input{\pathtotrunk diagrams/tikz/epsilon1.tex}%
%\endpgfgraphicnamed
 \, , \qquad 
 \epsilon_2 = %
%\beginpgfgraphicnamed{\pathtotrunk diagrams/tikz/#1-external}%
\input{\pathtotrunk diagrams/tikz/epsilon2.tex}%
%\endpgfgraphicnamed
 \, , \quad \ldots , \qquad 
 \epsilon_{2n}= %
%\beginpgfgraphicnamed{\pathtotrunk diagrams/tikz/#1-external}%
\input{\pathtotrunk diagrams/tikz/epsilon2k.tex}%
%\endpgfgraphicnamed
 \, . $$
 
 We say $S$ is a {\em rotational eigenvector with eigenvalue $\omega$} if   $\rho(S)= \omega S$ where
  $$\rho=%
%\beginpgfgraphicnamed{\pathtotrunk diagrams/tikz/#1-external}%
\input{\pathtotrunk diagrams/tikz/rho.tex}%
%\endpgfgraphicnamed
 \, .$$
 Note that $\omega$ must be a $n^{th}$ root of unity.
\end{defn}

As described in \cite{MR1929335}, every subfactor planar algebra is generated by uncappable rotational eigenvectors.

\begin{defn}\label{defn:relations}
If $S$ is an $n$-box,
we use the following names and numbers for relations on $S$:
\begin{enumerate}
\item $\rho(S)=-S$,
\item $S$ is uncappable,
\item $S^2 = \JW{n}$,
\item {\bf one-strand braiding substitute:}
	\begin{align*}
	\scalebox{0.8}{%
%\beginpgfgraphicnamed{\pathtotrunk diagrams/tikz/#1-external}%
\input{\pathtotrunk diagrams/tikz/InnerProducts/A.tex}%
%\endpgfgraphicnamed
}
	        & = i \frac{\sqrt{[n][n+2]}}{[n+1]}
	\scalebox{0.8}{%
%\beginpgfgraphicnamed{\pathtotrunk diagrams/tikz/#1-external}%
\input{\pathtotrunk diagrams/tikz/InnerProducts/B.tex}%
%\endpgfgraphicnamed
},
	\end{align*}
\item {\bf two-strand braiding substitute:}
	\begin{align*}
	\scalebox{0.8}{%
%\beginpgfgraphicnamed{\pathtotrunk diagrams/tikz/#1-external}%
\input{\pathtotrunk diagrams/tikz/InnerProducts/C.tex}%
%\endpgfgraphicnamed
} 
	         = \frac{[2][2n+4]}{[n+1][n+2]}
        \scalebox{0.8}{%
%\beginpgfgraphicnamed{\pathtotrunk diagrams/tikz/#1-external}%
\input{\pathtotrunk diagrams/tikz/InnerProducts/D.tex}%
%\endpgfgraphicnamed
}.
	\end{align*}
\end{enumerate}
\end{defn}

%Definition \ref{defn:ppres}, Proposition \ref{prop:braidingrelations}
%and Theorem \ref{thm:uniqueness}, below, will justify these seemingly
%ad-hoc definitions. 
We call relations (4) and (5) ``braiding substitutes'' because we
think of them as allowing us to move a generator ``through'' strands, rather like an identity
\begin{equation} \label{eqn:overbraid}
\begin{tikzpicture}[baseline]
    \clip (-1.1,-1.4) rectangle (1.1,1.4);

    \node (S) at (0,0) [circle, draw] {$X$};
    \node[anchor=east, inner sep=0] at (S.180) {\footnotesize$*$};

    \draw[rounded corners=2mm] (-1,-2) -- (-1,1) -- (1,1) -- (1,-2);

    \draw (S.315) -- ++(270:15mm);
    \draw (S.283) -- ++(270:15mm);
    \draw (S.257) -- ++(270:15mm);
    \draw (S.225) -- ++(270:15mm);
\end{tikzpicture}
=
\begin{tikzpicture}[baseline]
    \clip (-1.1,-1.4) rectangle (1.1,1.4);

    \node (S) at (0,0) [circle, draw] {$X$};
    \node[anchor=east, inner sep=0] at (S.180) {\footnotesize$*$};

   \draw (S.315) -- ++(270:15mm);
    \draw (S.283) -- ++(270:15mm);
    \draw (S.257) -- ++(270:15mm);
    \draw (S.225) -- ++(270:15mm);

       \draw[line width = 4pt, white, rounded corners=2mm] (-1,-1.5) --  (-1,-1) -- (1,-1) -- (1,-1.5);    
       \draw[rounded corners=2mm] (-1,-1.5) --  (-1,-1) -- (1,-1) -- (1,-1.5);

\end{tikzpicture}
\end{equation}
 in a braided tensor category. The planar algebras we consider in this paper are not braided, and do not satisfy the Equation \eqref{eqn:overbraid}. Nevertheless, we found it useful to look for relations that could play a similar role.  In particular, the evaluation algorithm described in \S \ref{sec:algorithm} was inspired by the evaluation algorithms in \cite{math/0808.0764} and \cite{MR2577673} for planar algebras of types $D_{2n}$, $E_6$, and $E_8$, where Equation \eqref{eqn:overbraid} holds. % for overcrossings but not for undercrossings.

\begin{defn}\label{defn:ppres}
Let $\ppres$ be the spherical planar algebra of modulus $[2]=d_k$,
generated by a $(4k+4)$-box $S$,
subject to relations (1)-(5) above.
\end{defn}

\begin{defn}
A {\em negligible element} of a spherical planar algebra $\cP$
is an element $x \in \cP_{n,\pm}$
such that the diagrammatic trace $\dtr{xy}$ is zero for all $y \in \cP_{n,\pm}$.
\end{defn}

%\begin{defn}
%A \emph{negligible element} of a planar algebra $\cP$ is an element $x \in \cP_{n,\pm}$ such that every planar composition landing in $\cP_{0,\pm'}$ of $x$ with other elements of the planar algebra is %zero.
%\end{defn}

%Note that when the planar algebra is spherical, this is the same as asking that the diagrammatic trace $\dtr{x y}$ is zero for all $y \in \cP_{n,\pm}$.
The set $\cN$ of all negligible elements of a planar algebra $\cP$
forms a planar ideal of $\cP$. In the presence of an antilinear involution $*$, we can replace $\dtr{xy}$ in the definition with $\dtr{x y^*}$ without changing the ideal. 
If the planar algebra is positive definite, then $\cN = 0$.
The following is well known.

\begin{prop}
\label{prop:simple}
Suppose $\cP$ is a spherical planar algebra with non-zero modulus and $\cN$ is the ideal of negligible morphisms.
If the spaces $\cP_{0,\pm}$ are one-dimensional
then every non-trivial planar ideal is contained in $\cN$.
\end{prop}

\begin{proof}
Suppose that a planar ideal $\cI$ contains a non-negligible element $x$ and without loss of generality assume $x \in \cP_{n,+}$. Then there is some element $y \in \cP_{n,+}$ so $\dtr{x y} \neq 0 \in \cP_{0,+}$. The element $\dtr{x y}$ is itself in the planar ideal, so since $\cP_{0,+}$ is one-dimensional, it must be entirely contained in $\cI$, and so the unshaded empty diagram is in the ideal. Drawing a circle around this empty diagram, and using the fact that the modulus is non-zero, shows that the shaded empty diagram is also in the ideal. Now, every box space $\cP_{m,\pm}$ is a module over $\cP_{0,\pm}$ under tensor product, with the empty diagram acting by the identity. Thus $\cP_{m,\pm} \subset \cI$ for all $m \in \Natural$.
\end{proof}

The sesquilinear pairing descends to $\pun$ and is then nondegenerate. 

Let $\rho^{1/2}$ denote the ``one-click'' rotation from $\cP_{n,+}$ to $\cP_{n,-}$ given by
$$%
%\beginpgfgraphicnamed{\pathtotrunk diagrams/tikz/#1-external}%
\input{\pathtotrunk diagrams/tikz/rhohalf.tex}%
%\endpgfgraphicnamed
.$$

\begin{defn}\label{defn:moments}
Let the {\em Haagerup moments} be as follows:
\begin{itemize}
\item $\tr{S^2}  =[n+1]$,
\item $\tr{S^3} = 0$,
\item $\tr{S^4} = [n+1]$,
\item $\tr{\rho^{1/2}(S)^3}  = i\frac{[2n+2]}{\sqrt{[n][n+2]}} $.
\end{itemize}
\end{defn}

\begin{prop}
\label{prop:braidingrelations}
Suppose $\cP$ is a positive definite spherical planar algebra
with modulus $d_k$,
and $S \in \cP_{n,+}$, where $n = 4k+4$.
If $S$ is uncappable
and satisfies $\rho(S) = -S$
and the Haagerup moments given in Definition~\ref{defn:moments}
then $S$ satisfies the five relations given in Definition~\ref{defn:relations}.
\end{prop}

We defer the proof until \S \ref{sec:relations}.
%subsection \S \ref{sec:provingrelations}.

\begin{thm} \label{thm:zeroboxonedim}
If $\cP$ is a planar algebra that is singly generated by an $n$-box
$S$
satisfying the five relations of Definition \ref{defn:relations}
then any closed diagram in $\cP_{0,+}$
is a multiple of the empty diagram (so $\dim (\cP_{o,+})=1$).
\end{thm}

The proof of Theorem \ref{thm:zeroboxonedim} is given in \S \ref{sec:algorithm}.

\begin{thm}
\label{thm:uniqueness}
If there exists a subfactor planar algebra $\cP$ with principal graphs $\Ha{k}$
then $\cP$ is isomorphic to $\pun$.
\end{thm}

\begin{proof}
The principal and dual principal graphs of $\cP$
each have their first trivalent vertex at depth $4k+3$.
In the language of \cite{quadratic},
$\cP$ has $n$-excess one,
where $n = 4k+4$.
We follow \cite[Section 5.1]{quadratic}.
There exists $S \in \cP_{n,+}$ such that $\pairing{S}{TL_{n,+}}=0,$ so
$$\cP_{n,+} =  TL_{n,+} \directSum \Complex S.$$
Let $r$ be the ratio of dimensions of
the two vertices at depth $4k+4$ on the principal graph,
chosen so that $r \geq 1$.
By \cite{quadratic},
we can choose $S$ to be self-adjoint, uncappable, and a rotational eigenvector,
such that
$$S^2  = (1-r) S + r \JW{n}.$$
(Our $S$ is $-\tilde{R}$ in \cite{quadratic}.)

The symmetry of $\Ha{k}^p$ implies that $r=1$, so $S^2 = \JW{n}$.
We can use this to compute powers of $S$ and their traces.
These agree with the first three Haagerup moments,
as given in Definition~\ref{defn:moments}.

In a similar fashion,
one defines $\check{r} \geq 1$ and
$\check{S}$ from the dual principal graph and calculates the 
moments of $\check{S}$.  Since the complement of $\check{S}$ is one-dimensional, $\check{S}$
and $\rho^{1/2}(S)$ must be multiples of each other; that multiple can be
calculated to be $\check{S} = \sqrt{\frac{\check{r}}{\omega r}} \rho^{1/2}(S)$.  (This can be done by comparing $\left< \check{S}, \check{S} \right>$ and $\left< \rho^{1/2}(S), \rho^{1/2}(S) \right>$, as we learned from an earlier draft of \cite{quadratic}).  Then it follows that \begin{equation}
\label{eq:twistedmoment}
\tr{\rho^{1/2}(S)^3}   =
  \omega^{3/2} \sqrt{\frac{r}{\check{r}}} r (\check{r} - 1) [n+1],
\end{equation}
for some square root $\omega^{1/2}$ of the rotational eigenvalue of $S$.

Although we will not use it here, we record the identity
\begin{equation} \rho^{1/2}(S)^2 = - \omega^{1/2} r^{1/2} (\check{r}^{1/2} - \check{r}^{-1/2}) \rho^{1/2}(S) + \omega^{-1} r \JW{n}, \end{equation}
which is equivalent to $$\check{S}^2 = (1-\check{r})\check{S} + \check{r} \JW{n}.$$

By \cite[Theorem 5.1.11]{quadratic}, whenever $\cP$ has $n$-excess one then $\check{r} = \frac{[n+2]}{[n]}$ and
$$r+\frac{1}{r} = 2+ \frac{2 + \omega +\omega^{-1}}{[n][n+2]}.$$
Since $r=1$, this implies that $\omega = -1$.
Note that we could also compute $\check{r}$
directly from the dual principal graph.
Jones's proof that $\omega=-1$ uses the converse of a result along the lines of Lemma \ref{lem:ABrel},
since there must be some linear relation of the form given in that Lemma.

In the case $r=1$ we also have the freedom to replace $S$ with $-S$,
which we use to (arbitrarily) choose the square root $\omega^{1/2} = -i$.
Substituting all these quantities into Equation \eqref{eq:twistedmoment}, and using the identity $[n+1]([n+2]-[n]) = [2n+2]$,
we now see that $S$ has all of the Haagerup moments,
as given in Definition~\ref{defn:moments}.

By Proposition \ref{prop:braidingrelations},
$S$ satisfies the relations given in Definition~\ref{defn:relations}.
Thus there is a planar algebra morphism
$\ppres \to \cP$ given by sending $S$ to $S$.
Since $\cP$ is positive definite,
this descends to the quotient to give a map 
$$\Phi \co \pun \to \cP.$$

By Proposition \ref{prop:simple},
$\pun$ has no nontrivial proper ideals.
Since $\Phi$ is non-zero,
it must be injective.
The image of $\Phi$ is the planar algebra in $\cP$ generated by $S$.
This is a subfactor planar algebra with the same modulus as $\cP$.   Its
principal graphs are not $(A_\infty, A_\infty)$ because the dimension of the $n$-box space is too large.  Haagerup's
classification shows that the principal graphs of the image of $\Phi$ must be $\Ha{k}$.  However, since the principal graphs determine the dimensions of all box spaces, the image of $\Phi$ must be all of $\cP$.
Thus $\Phi$ is an isomorphism of planar algebras.
\end{proof}

%---------------------------
\subsection{Existence}

The subfactor planar algebra with principal graphs $\Ha{0}$
is called the Haagerup planar algebra,
and is isomorphic to $\punzero$.
The corresponding subfactor was constructed
in \cite{MR1686551} and \cite{MR1832764}.
The subfactor planar algebra was directly constructed in \cite{0902.1294}.
There is no subfactor planar algebra with principal graphs $\Ha{k}$
for $k >1$.
In this case,
by \cite{MR2472028},
$\pun$ cannot be a finite depth planar algebra,
let alone a subfactor planar algebra.
The following theorem deals with the one remaining case.

\begin{thm}
\label{thm:existence}
There is a subfactor planar algebra with principal graphs $\Ha{1}$.
\end{thm}

We prove this by finding $\Ha{1}$ as a sub-planar algebra of the graph planar algebra of one of the extended Haagerup graphs.
The following lemma simplifies the proof of irreducibility 
for subalgebras of graph planar algebras.

\begin{lem}
\label{lem:irreducible}
If $\cP \subset GPA(G)$,
$\dim{\cP_{0,+}}=1$,
and $G$ has an even univalent vertex,
then $\cP$ is an irreducible subfactor planar algebra.
\end{lem}
 
\begin{proof}
To show that $\cP$ is an irreducible subfactor planar algebra, 
we need to show that 
$\cP$ is spherical and positive definite, and that
$\dim \cP_{0,\pm}=1$ and  $\dim \cP_{1,+}=1$.
By Proposition \ref{GPAsphericalposdef}, 
the graph planar algebra is spherical and positive definite.
The subalgebra $\cP$ inherits both of these properties.
We are given $\dim \cP_{0,+}=1$.
Also,
$\cP_{0,-}$ injects into $\cP_{1,+}$,
(by tensoring with a strand on the left).
It remains only to show that $\dim \cP_{1,+} = 1$.

Let $v$ be an even univalent vertex of $G$.
Let $w$ be the unique vertex connected to $v$.
Suppose $X \in \cP_{1,+}$ is some functional on paths of length two in $G$.

Now $\dtr{X}$ is a closed diagram with unshaded exterior. (Note here we use $\operatorname{tr}_0$, the diagrammatic trace, without applying a partition function, even though $\dim{P_{0,+}} =1$.)
This is a functional defined on even vertices, via a state sum.
Since $v$ is univalent,
the state sum for $\dtr{X}(v)$ has only one term,
giving
$$\dtr{X}(v) = X(vw) \frac{d_w}{d_v}.$$
Similarly,
$$\dtr{X^*X}(v) = X(vw) X^*(vw) \frac{d_w}{d_v}$$
Thus
if $\dtr{X}(v)$ is zero
then $\dtr{X^*X}(v)$ is zero also.
Note that $\dtr{X}$ and $\dtr{X^*X}$
are both scalar multiples of the empty diagram, and so $\dtr{X^* X}(v)=0$ implies that $\dtr{X^* X}=0$.

Therefore, if $\dtr{X}$ is zero then $\dtr{X^*X}$ is zero.
Then by positive definiteness,
if $\dtr{X^* X}$ is zero then $X$ is zero.

We conclude that the diagrammatic trace function is injective on $\cP_{1,+}$ and thus $\cP_{1,+}$ is one-dimensional.
\end{proof}

Recall the Haagerup moments from Definition \ref{defn:moments}.
In the current setting,
$n=8$, $[2]=d_1$, and the Haagerup moments are as follows.
\begin{itemize}
\item $\tr{S^2}  =[9] \approx 24.66097$,
\item $\tr{S^3}  =0$,
\item $\tr{S^4}  =[9]$,
\item $\tr{\rho^{1/2}(S)^3} = i\frac{[18]}{\sqrt{[8][10]}} \approx 15.29004i$.
%\tr{\rho^{1/2}(S)^4} & = \left(\frac{[8]}{[10]}-1+\frac{[10]}{[8]}\right)[9]
%                       \approx 34.1409... 
\end{itemize}

\begin{prop}\label{prop:SgivesSPA}
Suppose that $S \in GPA(\Ha{1}^p)_{8,+}$ is self-adjoint, uncappable, a rotational eigenvector with eigenvalue $-1$, and has the above Haagerup moments.
Let $\PA(S)$ be the subalgebra of $GPA(\Ha{1}^p)_{8,+}$ generated by $S$.
Then $\PA(S)$ is an irreducible subfactor planar algebra with principal graphs $\Ha{1}$.
\end{prop}

\begin{proof}
By Proposition \ref{prop:braidingrelations}, $S \in GPA(\Ha{1}^p)$ satisfies all of the relations used to define $\cQ^1$.  
Thus by Theorem \ref{thm:uniqueness}, $\PA(S)$ is isomorphic to $\punone$.
By Theorem \ref{thm:zeroboxonedim}, $(\punone)_{0,+}$ is $1$-dimensional.  By Lemma \ref{lem:irreducible} it follows that $\punone$ is an irreducible subfactor planar algebra.
By Haagerup's classification \cite{MR1317352} it follows that the principal graphs of $\punone$ must be the unique possible graph pair with the correct graph norm, namely $\Ha{1}$.
\end{proof}

To prove Theorem \ref{thm:existence},
it remains to find $S \in GPA(\Ha{1}^p)_{8,+}$
satisfying the requirements of the above proposition.
We defer this to Section \ref{sec:propertiesofS},
where we give an explicit description of $S$
and some long computations of the moments,
assisted by computer algebra software.

%% file: diagrams/tikz/epsilon1.tex
\begin{tikzpicture}[annular]
	\clip (0,0) circle (2cm);

	\filldraw[shaded] (0,0) .. controls ++(170:2cm) and ++(100:2cm) .. (0,0);
	\filldraw[shaded] (-158:4cm)--(0,0)--(-112:4cm);

	\draw[shaded] (68:4cm)--(0,0)--(-68:4cm)--(0:10cm);
	
	\draw[ultra thick] (0,0) circle (2cm);
	
	\node at (0,0)  [empty box] (T) {};
	\node at (T.180) [left] {$\star$};
	\node at (180:2cm) [right] {$\star$};
	\node at (0:1cm) {$\cdot$};
	\node at (20:1cm) {$\cdot$};
	\node at (-20:1cm) {$\cdot$};
	
\end{tikzpicture}

%% file: diagrams/tikz/epsilon2.tex
\begin{tikzpicture}[annular]
	\clip (0,0) circle (2cm);

	\filldraw[shaded] (158:4cm) -- (0,0) .. controls ++(130:2cm) and ++(50:2cm) .. (0,0)--(-68:4cm) arc (-68:158:4cm);
	\filldraw[shaded] (-158:4cm)--(0,0)--(-112:4cm);

	\draw[ultra thick] (0,0) circle (2cm);
	
	\node at (0,0)  [empty box] (T) {};
	\node at (T.180) [left] {$\star$};
	\node at (180:2cm) [right] {$\star$};
	\node at (0:1cm) {$\cdot$};
	\node at (20:1cm) {$\cdot$};
	\node at (-20:1cm) {$\cdot$};
	
\end{tikzpicture}

%% file: diagrams/tikz/epsilon2k.tex
\begin{tikzpicture}[annular]
	\clip (0,0) circle (2cm);

	\filldraw[shaded] (112:4cm) -- (0,0) .. controls ++(140:2cm) and ++(-140:2cm) .. (0,0)--(-112:4cm) arc (-112:-202:4cm);
	
	\draw[shaded] (68:4cm)--(0,0)--(-68:4cm)--(0:10cm);
	
	\draw[ultra thick] (0,0) circle (2cm);
	
	\node at (0,0)  [empty box] (T) {};
	\node at (T.180) [left] {$\star$};
	\node at (90:2cm) [below] {$\star$};
	\node at (0:1cm) {$\cdot$};
	\node at (20:1cm) {$\cdot$};
	\node at (-20:1cm) {$\cdot$};
	
\end{tikzpicture}

%% file: diagrams/tikz/rho.tex
\begin{tikzpicture}[annular]
	\clip (0,0) circle (2cm);

	\filldraw[shaded] (158:4cm)--(0,0)--(112:4cm);
	\filldraw[shaded] (-158:4cm)--(0,0)--(-112:4cm);
	
	\draw[shaded] (68:4cm)--(0,0)--(-68:4cm)--(0:10cm);
	
	\draw[ultra thick] (0,0) circle (2cm);
	
	\node at (0,0)  [empty box] (T) {};
	\node at (T.180) [left] {$\star$};
	\node at (-90:2cm) [above] {$\star$};
	\node at (0:1cm) {$\cdot$};
	\node at (20:1cm) {$\cdot$};
	\node at (-20:1cm) {$\cdot$};
\end{tikzpicture}

%% file: diagrams/tikz/InnerProducts/A.tex
\begin{tikzpicture}[STrain]
	\RainbowOne;
        \draw (0,0)--(0,-0.5);
        \node[anchor=west] at (0,-0.35) {\footnotesize$2n+2$};
	\JWPlusTwo;
\end{tikzpicture}

%% file: diagrams/tikz/InnerProducts/B.tex
\begin{tikzpicture}[STrain]
	\STrainStrings{$n+1$}{$n+1$} \STrainOne
        \draw (0,0)--(0,-0.5);
        \node[anchor=west] at (0,-0.35) {\footnotesize$2n+2$};
	\JWPlusTwo
\end{tikzpicture}

%% file: diagrams/tikz/InnerProducts/C.tex
\begin{tikzpicture}[STrain]
	\RainbowTwo
        \draw (0,0)--(0,-0.5);
        \node[anchor=west] at (0,-0.35) {\footnotesize$2n+4$};
	\JWPlusFour
\end{tikzpicture}

%% file: diagrams/tikz/InnerProducts/D.tex
\begin{tikzpicture}[STrain]
	\STrainThreeStrings{$n+1$}{$2$}{$n+1$} \STrainOneOne
        \draw (0,0)--(0,-0.5);
        \node[anchor=west] at (0,-0.35) {\footnotesize$2n+4$};
	\JWPlusFour
\end{tikzpicture}

%% file: diagrams/tikz/rhohalf.tex
\begin{tikzpicture}[annular]
	\clip (0,0) circle (2cm);

	\filldraw[shaded] (158:4cm)--(0,0)--(112:4cm);
	\filldraw[shaded] (-158:4cm)--(0,0)--(-112:4cm);
	
	\draw[shaded] (68:4cm)--(0,0)--(-68:4cm)--(0:10cm);
	
	\draw[ultra thick] (0,0) circle (2cm);
	
	\node at (0,0)  [empty box] (T) {};
	\node at (T.180) [left] {$\star$};
	\node at (-135:2.1cm) [above right] {$\star$};
	\node at (0:1cm) {$\cdot$};
	\node at (20:1cm) {$\cdot$};
	\node at (-20:1cm) {$\cdot$};
\end{tikzpicture}

%% file: text/algorithm.tex
The aim of this section is to prove that the relations of Definition \ref{defn:relations} enable us to reduce any closed diagram built from copies of $S$ to a scalar multiple of the empty diagram.

Questions of whether you can evaluate an arbitrary closed diagram are ubiquitous in quantum topology.  The simplest such algorithms (e.g., the Kauffman bracket algorithm for knots) involve decreasing the number of generators (in this case, crossings) at each step.  Slightly more complicated algorithms (e.g., HOMFLY evaluations) include steps that leave the number of generators constant while decreasing some other measure of complexity (such as the unknotting number).  Another common technique is to apply Euler characteristic arguments to find a small ``face" (with generators thought of as vertices) that can then be removed.  Again, the simplest such algorithms decrease the number of faces at every step (e.g., Kuperberg's rank $2$ spiders \cite{MR1403861}), while more difficult algorithms require steps that maintain the number of faces before removing a face (e.g., Peters' approach to $\Ha{0}$ in \cite{0902.1294}).  In all of these algorithms, the number of generators is monotonically non-increasing as the algorithm proceeds.  The algorithm we describe below is unusual in that it initially increases the number of generators in order to put them in a desirable configuration.  We hope that this technique will be of wider interest in quantum topology (see \cite{han-2221} for a subsequent application of this technique).  Therefore we have written this section to be independent of the rest of the paper, apart from references to Definitions \ref{defn:uncappable} and \ref{defn:relations}.

The algorithm we will describe gives a proof of Theorem \ref{thm:zeroboxonedim}, which we repeat from above:
\newtheorem*{rthm}{Theorem \ref{thm:zeroboxonedim}}
\begin{rthm}
If $\cP$ is a planar algebra that is singly generated by an $n$-box
$S$
satisfying the five relations of Definition \ref{defn:relations}
then any closed diagram in $\cP_{0,+}$
is a multiple of the empty diagram (so $\dim (\cP_{o,+})=1$).
\end{rthm}

We do not actually need the
full strength of the relations of Definition \ref{defn:relations}.
The theorem is true for any planar algebra
that is singly generated by an $n$-box $S$
such that:
\begin{itemize}
\item $S$ is a rotational eigenvector: $\rho(S)=\omega S$ for some $\omega$,
\item $S$ is uncappable (see Definition \ref{defn:uncappable}),
\item $S^2 = a S + b f^{(n)}$ for some scalars $a$ and $b$ (multiplication is defined in Figure \ref{fig:timestracetensor}), and
\item $S$ satisfies one- and two-strand braiding substitutes of the form:
	\begin{align*}
	\scalebox{0.8}{%
%\beginpgfgraphicnamed{\pathtotrunk diagrams/tikz/#1-external}%
\input{\pathtotrunk diagrams/tikz/InnerProducts/A.tex}%
%\endpgfgraphicnamed
}
	        & =x \;
	\scalebox{0.8}{%
%\beginpgfgraphicnamed{\pathtotrunk diagrams/tikz/#1-external}%
\input{\pathtotrunk diagrams/tikz/InnerProducts/B.tex}%
%\endpgfgraphicnamed
} \displaybreak[1] \\
	\scalebox{0.8}{%
%\beginpgfgraphicnamed{\pathtotrunk diagrams/tikz/#1-external}%
\input{\pathtotrunk diagrams/tikz/InnerProducts/C.tex}%
%\endpgfgraphicnamed
} 
	        & = y
        \scalebox{0.8}{%
%\beginpgfgraphicnamed{\pathtotrunk diagrams/tikz/#1-external}%
\input{\pathtotrunk diagrams/tikz/InnerProducts/D.tex}%
%\endpgfgraphicnamed
} 
	\end{align*}
	for some scalars $x$ and $y$ in $\Complex$.
\end{itemize}

\begin{figure}[!htb]
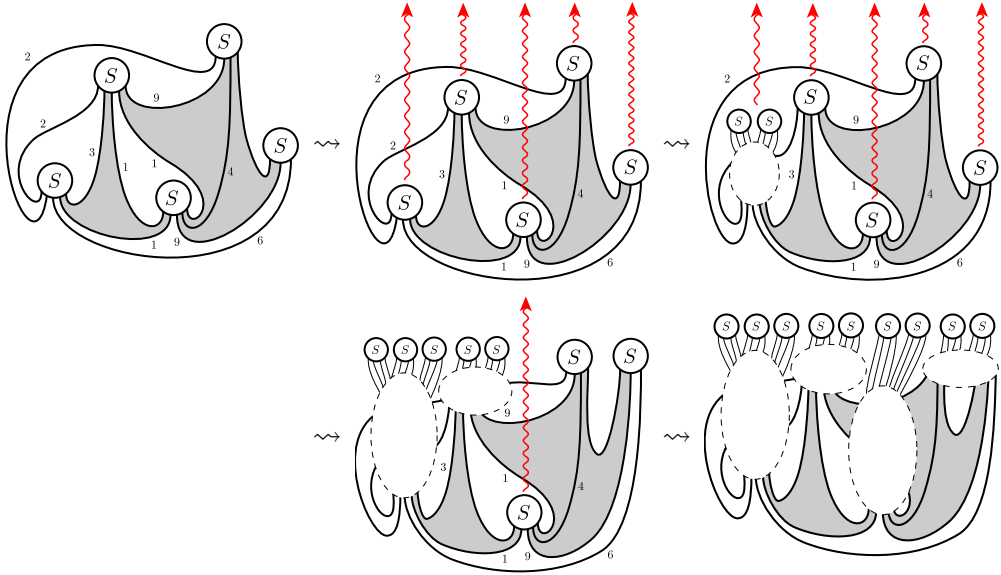

\begin{align*}
\mathfig{0.28}{jellyfish/network} & \rightsquigarrow
\mathfig{0.28}{jellyfish/network-paths} \rightsquigarrow
\mathfig{0.28}{jellyfish/network-paths-2} \\ & \rightsquigarrow
\mathfig{0.28}{jellyfish/network-paths-3} \rightsquigarrow
\mathfig{0.28}{jellyfish/network-paths-4}
\end{align*}
\caption{The initial steps of the jellyfish algorithm. The dotted ovals represent linear combinations of Temperley-Lieb diagrams. This is only a schematic illustration - to be precise, the result should be a linear combination of diagrams with various (sometimes large) numbers of copies of $S$.}
\label{fig:jellyfish}
\end{figure}

Before going through the details, we briefly sketch the idea.
First use the one- and two-strand braiding substitutes
to pull all copies of $S$ to the outside of the diagram.
This will usually increase the number of copies of $S$.
We can then guarantee that there is a pair of copies of $S$ connected by at least $n$ strands.
This is a copy of $S^2$, which we can then express using fewer copies of $S$.
All copies of $S$ remain on the outside,
and so we can again find a copy of $S^2$.
Repeating this
eventually gives an element of the Temperley-Lieb planar algebra, which is evaluated as usual. See Figure \ref{fig:jellyfish} for an example. We like to think of the copies of $S$ as ``jellyfish floating to the surface,'' and hence the name for the algorithm.

\begin{defn}
Suppose $D$ is a diagram in $\cP$.
Let $S_0$ be a fixed copy of the generator inside $D$.
Suppose $\gamma$ is an embedded arc in $D$
from a point on the boundary of $S_0$
to a point on the top edge of $D$.
Suppose $\gamma$ is in general position,
meaning that it intersects the strands of $D$ transversely,
and does not touch any generator except at its initial point on $S_0$.
Let $m$ be the number of points of intersection
between $\gamma$ and the strands of $D$.
If $m$ is minimal over all such arcs $\gamma$
then we say $\gamma$ is a {\em geodesic}
and $m$ is the {\em distance} from $S_0$ to the top of $D$.
\end{defn}

\begin{lem} \label{lem:oneup}
Suppose $X$ is a diagram consisting of
one copy of $S$ with all strands pointing down,
and $d$ parallel strands forming a ``rainbow'' over $S$,
where $d \ge 1$.
Then $X$ is a linear combination of diagrams
that contain at most three copies of $S$,
each having distance less than $d$ from the top of the diagram.
\end{lem}

\begin{proof}
First consider the case $d=1$.
Up to some number of applications of the rotation relation $\rho(S)= \omega S$,
$X$ is as shown in Figure~\ref{fig:T-hat}.

\begin{figure}[!htb]
\center
\begin{tikzpicture}[STrain] \RainbowOne \end{tikzpicture}
\caption{$X$ in the case $d = 1$.}
\label{fig:T-hat}
\end{figure}
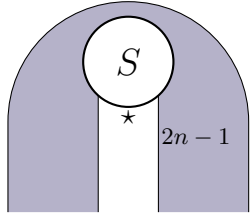

Recall that we have the relation
$$ \scalebox{0.8}{%
%\beginpgfgraphicnamed{\pathtotrunk diagrams/tikz/#1-external}%
\input{\pathtotrunk diagrams/tikz/InnerProducts/A.tex}%
%\endpgfgraphicnamed
}
                 = x \; \scalebox{0.8}{%
%\beginpgfgraphicnamed{\pathtotrunk diagrams/tikz/#1-external}%
\input{\pathtotrunk diagrams/tikz/InnerProducts/B.tex}%
%\endpgfgraphicnamed
}.
$$

Consider what happens to the left side of the above relation
when we write $f^{(2n+2)}$
as a linear combination of Temperley-Lieb diagrams $\beta$.
The term in which $\beta$ is the identity occurs with coefficient one,
and gives the diagram $X$.
Suppose $\beta$ is not the identity.
Then $\beta$ contains a cup that connects two adjacent strands from $X$.
If both ends of the cup are attached to $S$ then the resulting diagram is zero.
If not,
then the cup must be at the far left or the far right of $\beta$.
Such a cup converts $X$ to a rotation of $S$,
so gives distance zero from $S$ to the top of the diagram.

Now consider what happens to the right side of the above relation
when we write $f^{(2n+2)}$
as a linear combination of Temperley-Lieb diagrams $\beta$.
Every term in this expansion
is a diagram with two copies of $S$,
each of having distance zero from the top.

By rearranging terms in the one-strand braiding substitute,
we can write $X$ as a linear combination of diagrams
that contain one or two copies of $S$,
each having distance zero from the top of the diagram.
This completes the case $d = 1$.

The case $d = 2$ is similar,
but we use the two-strand braiding substitute.

Finally,
suppose $d > 2$.
If $d$ is odd
then $\gamma$ begins in a shaded region of $X$.
Then $X$ contains a copy of the diagram shown in Figure~\ref{fig:T-hat},
up to the rotation relation $\rho(S) = \omega S$.
We can therefore use the one-strand braiding substitute, as we did in the case $d=1$.
Similarly,
if $d$ is even then we use the two-strand braiding substitute.
\end{proof}

\begin{defn}
We say a diagram $D$ in $\cP$ is in {\em jellyfish form} if all occurrences of $S$ lie in a row at the top of $D$,
and all strands of $D$
lie entirely below the height of the tops of the copies of $S$.
\end{defn}

\begin{lem}\label{lem:jellyfish}
Every diagram in $\cP$
is a linear combination of diagrams in jellyfish form.
\end{lem}

\begin{proof}
Suppose $D$ is a diagram in $\cP$ (not necessarily closed),
drawn in such a way that all endpoints lie on the bottom edge of $D$.
If every copy of the generator in $D$
is distance zero from the top edge of $D$
then $D$ is already in jellyfish form,
up to isotopy.
If not,
we will use Lemma \ref{lem:oneup} to pull each copy of $S$ to the top $D$.
It is convenient for our proof, 
but not necessary for the algorithm,
to move copies of our generator $S$
along geodesics.

Suppose $S_0$ is a copy of $S$ that has distance $d$ from the top of $D$,
where $d \ge 1$.
Let $\gamma$ be a geodesic from $S_0$ to the top edge of $D$.
Let $X$ be a small neighborhood of $S_0 \cup \gamma$.
By applying an isotopy,
we consider $X$ to be a diagram in a rectangle,
consisting of a copy of $S_0$ with all strands pointing down,
and a ``rainbow'' of $d$ strands over it.

By Lemma \ref{lem:oneup},
$X$ is a linear combination of diagrams
that contain at most three copies of $S$,
each having distance less than $d$ from the top of the diagram.
Let $X'$ be one of the terms
in this expression for $X$.
Let $D'$ be the result of replacing $X$ by $X'$ in $D$.

Suppose $S_1$ is a copy of the generator in $D'$.
If $S_1$ lies in $X'$
then the distance from $S_1$ to the top of $D'$ is at most $d-1$.
Now suppose
$S_1$ does not lie in $X'$.
By basic properties of geodesics,
there is a geodesic in $D$
from $S_1$ to the top of $D$
that does not intersect $\gamma$.
This geodesic is still a path in general position in $D'$,
and still intersects strands in the same number of points.
Thus the distance from $S_1$ to the top of $D$
does not increase when we replace $X$ by $X'$.

In summary,
if we replace $X$ by $X'$,
then $S_0$ will be replaced by one, two or three copies of $S$
that are closer to the top of $D$,
and no other copy of $S$ will become farther from the top of $D$.
Although the number of copies of $S$ may increase,
it is not hard to see that this process must terminate. For example,
we have decreased the sum over each generator $S_0$
of $4$ to the power of the distance from $S_0$ to the top.
\end{proof}

We now prove Theorem \ref{thm:zeroboxonedim}, that $\dim(\cP_{0,+}) = 1$.

\begin{proof}[Proof of Theorem \ref{thm:zeroboxonedim}.]
Suppose $D$ is a closed diagram with unshaded exterior.
We must show that $D$ is a scalar multiple of the empty diagram.
By the previous lemma,
we can assume $D$ is in jellyfish form.
We can also assume there are no closed loops or cups attached to generators,
so that every strand must connect two different copies of the generator.

We argue that there
is a copy of the generator whose strands only go to
one or both of its immediate neighbors.
This is a simple combinatorial fact about this kind of planar graph.
Think of the copies of the generator as vertices,
and consider all strands that do not connect adjacent vertices.
Amongst these,
find one that has the smallest (positive) number of vertices
between its endpoints.
Any vertex between the endpoints of this strand
can connect only to its two neighbors.

Let $S_0$ be a copy of the generator
such that the strands of $S_0$ only go to
one or both of its immediate neighbors.
Then $S_0$ is connected to some neighbor, $S_1$,
by at least $n$ parallel strands.
See Figure \ref{fig:Spair} for an example.
Recall that $S^2 = aS + bf^{(n)}$.
Thus we can replace $S_0$ and $S_1$ with $aS + bf^{(n)}$,
giving a linear combination of diagrams
that are still in jellyfish form,
but contain fewer copies of the generator.
By induction,
$D$ is a scalar multiple of the empty diagram.

\begin{figure}[!htb]
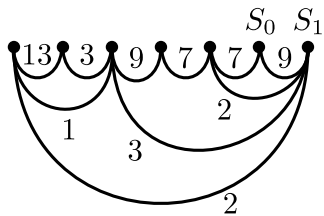

$$\mathfig{0.3}{jellyfish/curved-diagonals}$$
\caption{Jellyfish form, illustrating (with $n=8$) the proof of Theorem \ref{thm:zeroboxonedim}.}
\label{fig:Spair}
\end{figure}

\end{proof}

%% file: text/newrelations.tex
\newcommand{\tablestrut}{\rule[-1.5ex]{0pt}{4.4ex}}
% This is approximately the height of \left( \frac{[n+1]}{[n]} \right)^2

%---------------------------
In this section we prove Proposition \ref{prop:braidingrelations}, which says that certain conditions on an element $S$ imply the five relations of Definition \ref{defn:relations}. We use this Proposition once in the proof of Theorem \ref{thm:uniqueness}, and once in the proof of Theorem \ref{thm:existence}. These are the uniqueness and existence results. In the proof of uniqueness, we must show that a certain subfactor planar algebra $\cP$ is isomorphic to $\pun$. We use Proposition \ref{prop:braidingrelations} to show that an element $S$ of $\cP$ satisfies the defining relations of $\pun$. In the proof of existence, we must show that a certain subalgebra $\cP$ of a graph planar algebra has a one-dimensional space of $n$-boxes. We use Proposition \ref{prop:braidingrelations} to show that the generator $S$ of $\cP$ satisfies relations, which we then use in the algorithm of \S \ref{sec:algorithm}.

Most of this section consists of computations of inner products between diagrams.  Since the values of these inner products may be useful for studying other planar algebras, we strive to use weaker assumptions whenever possible.

\begin{assumption}
\label{ass:a}
$\cP$ is a spherical planar algebra with modulus $[2] = q + q^{-1}$,
where $q$ is not a root of unity (so we can safely divide by quantum integers).
Furthermore, $S \in \cP_{n,+}$ is uncappable
and has rotational eigenvalue $\omega$.
\end{assumption}

Recall that the Haagerup moments are as follows.

\begin{itemize}
\item $\tr{S^2}  =[n+1] $,
\item $\tr{S^3} = 0$,
\item $\tr{S^4}  =[n+1] $,
\item $\tr{\rho^{1/2}(S)^3}  = i\frac{[2n+2]}{\sqrt{[n][n+2]}} $.
\end{itemize}

\begin{assumption}
\label{ass:b}
$\cP$ is positive definite and has modulus $[2] = d_k$, where $n=4k+4$.
Furthermore, $S \in \cP_{n,+}$ has rotational eigenvalue $\omega = -1$,
and has the Haagerup moments.
\end{assumption}

\newtheorem*{propbraidingrelations}{Restatement of Proposition \ref{prop:braidingrelations}}

\begin{propbraidingrelations}
Suppose $\cP$ is a planar algebra, $S \in \cP_{n,+}$,
and Assumptions \ref{ass:a} and \ref{ass:b} hold.
Then $S$ satisfies the five relations given in Definition~\ref{defn:relations}.
\end{propbraidingrelations}

The proof involves some long and difficult computations,
but the basic idea is very simple.
We will define diagrams $A$, $B$, $C$ and $D$.
We must prove certain linear relations hold between $A$ and $B$,
and between $C$ and $D$.
Since $\cP$ is positive definite,
we can do this by computing certain inner products.
In general,
there is a linear relation between $X$ and $Y$
if and only if
$$ \langle X,X \rangle \langle Y,Y \rangle =
  | \langle X,Y \rangle |^2.$$
In this case,
$$ \langle Y,Y \rangle X - \langle X,Y \rangle Y = 0,$$
as can be seen by taking the inner product of this expression with itself.

To compute the necessary inner products
we must evaluate certain closed diagrams.
Most of these closed diagrams involve a Jones-Wenzl idempotent.
In principal,
we could expand this idempotent
into a linear combination of Temperley-Lieb diagrams,
and evaluate each resulting tangle in terms of the moments,
or using relations that have already been proved.
In practice,
we must take care to avoid dealing with
an unreasonably large number of terms.

%-------------------------------
\subsection{Definitions and conventions}

\begin{notation*}
We use the notation that
a thick strand in a Temperley-Lieb diagram
always represents $n-1$ parallel strands.  For example,
$$%
%\beginpgfgraphicnamed{\pathtotrunk diagrams/tikz/#1-external}%
\input{\pathtotrunk diagrams/tikz/TLn/iIiiIi.tex}%
%\endpgfgraphicnamed
 $$
is the identity of $TL_{2n+2}$.
\end{notation*}

\begin{defn}
For $m \ge 0$, let
$$W_m = q^m+q^{-m}-\omega-\omega^{-1},$$
as in \cite[Definition 4.2.6]{quadratic}.
\end{defn}

The diagrams $A$, $B$, $C$ and $D$ 
of Figures \ref{fig:AB} and \ref{fig:CD} 
are the terms in the ``braiding'' relations we wish to prove.  

\begin{figure}[!htb]
$$A=%
%\beginpgfgraphicnamed{\pathtotrunk diagrams/tikz/#1-external}%
\input{\pathtotrunk diagrams/tikz/InnerProducts/A.tex}%
%\endpgfgraphicnamed
, \quad B=%
%\beginpgfgraphicnamed{\pathtotrunk diagrams/tikz/#1-external}%
\input{\pathtotrunk diagrams/tikz/InnerProducts/B.tex}%
%\endpgfgraphicnamed
.$$
\caption{$A$ and $B$}
\label{fig:AB}
\end{figure}

\begin{figure}[!htb]
$$C=%
%\beginpgfgraphicnamed{\pathtotrunk diagrams/tikz/#1-external}%
\input{\pathtotrunk diagrams/tikz/InnerProducts/C.tex}%
%\endpgfgraphicnamed
, \quad D=%
%\beginpgfgraphicnamed{\pathtotrunk diagrams/tikz/#1-external}%
\input{\pathtotrunk diagrams/tikz/InnerProducts/D.tex}%
%\endpgfgraphicnamed
.$$
\caption{$C$ and $D$}
\label{fig:CD}
\end{figure}

\begin{figure}[!htb]
$$\Gamma=%
%\beginpgfgraphicnamed{\pathtotrunk diagrams/tikz/#1-external}%
\input{\pathtotrunk diagrams/tikz/InnerProducts/Tet.tex}%
%\endpgfgraphicnamed
, \quad B'=%
%\beginpgfgraphicnamed{\pathtotrunk diagrams/tikz/#1-external}%
\input{\pathtotrunk diagrams/tikz/InnerProducts/Bprime.tex}%
%\endpgfgraphicnamed
.$$
\caption{$\Gamma$ and $B'$}
\label{fig:GammaBprime}
\end{figure}
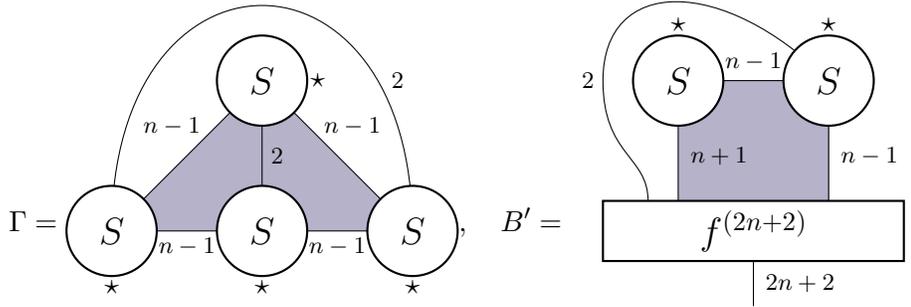

Along the way, we will also use the diagrams $\Gamma$ and $B'$,
as shown in Figure~\ref{fig:GammaBprime}.
%Note that $\Gamma$ is an example of a tetrahedral structure constant
%from the draft versions of \cite{quadratic}.

%---------------------------
\subsection{Computing inner products}

We now calculate the necessary inner products.

\begin{lem}
\label{lem:AA}
If Assumption~\ref{ass:a} holds then
$$ \langle A,A \rangle = \frac{1}{[2n+2]} W_{2n+2} \tr{S^2}.$$
The same holds with the reverse shading.
\end{lem}

\begin{proof}
We must evaluate the closed diagram
$$\pairing{A}{A}=Z\left(\scalebox{0.8}{%
%\beginpgfgraphicnamed{\pathtotrunk diagrams/tikz/#1-external}%
\input{\pathtotrunk diagrams/tikz/InnerProducts/AA.tex}%
%\endpgfgraphicnamed
}\right).$$
The idea is to apply Lemma~\ref{lem:yetanother} to the copy of $f^{(2n+2)}$,
and then evaluate each of the resulting diagrams.

Consider the first term.
Here,
$f^{(2n+2)}$ is replaced by a copy of $f^{(2n+1)}$
together with a single vertical strand on the right.
Since the planar algebra is spherical,
we can drag this strand over to the left.
This results in a partial trace of $f^{(2n+1)}$,
which is equal to
$$ \frac{[2n+2]}{[2n+1]} f^{(2n)}.$$
By noting that $S \cdot f^{(n)} = S$, we  obtain
$$ \frac{[2n+2]}{[2n+1]} \tr{S^2}$$
as the value of the first term.

Now consider the terms in the sum over $a$ and $b$.
Here,
$f^{(2n+2)}$ is replaced by a copy of $f^{(2n)}$
together with a ``cup'' and a ``cap''
in positions given by $a$ and $b$.
In most cases,
the resulting diagram is zero because $S$ is uncappable.
We only need to consider the four cases where $a,b \in \{1,2n+1\}$.
Each of these gives $\tr{S^2}$,
up to some rotation of one or both copies of $S$.
We obtain
$$ \frac{1}{[2n+1][2n+2]}
   ( - 1 - [2n+1]^2 - [2n+1] \omega - [2n+1] \omega^{-1} )
   \tr{S^2}.$$

The result now follows by adding the above two expressions
and writing the quantum integers in terms of $q$.
\end{proof}

\begin{lem}
\label{lem:AB}
If Assumption~\ref{ass:a} holds
and either $\omega=-1$ or $S^2 = f^{(n)}$
then
$$\langle A,B \rangle = \tr{\rho^{1/2}(S)^3}.$$
\end{lem}

\begin{proof}
We must evaluate the closed diagram
$$\pairing{A}{B}=Z\left(\scalebox{0.8}{%
%\beginpgfgraphicnamed{\pathtotrunk diagrams/tikz/#1-external}%
\input{\pathtotrunk diagrams/tikz/InnerProducts/AB.tex}%
%\endpgfgraphicnamed
}\right).$$
Consider the complete expansion of $f^{(2n+2)}$ into
a linear combination of Temperley-Lieb diagrams $\beta \in TL_{2n+2}$.
For most such $\beta$,
the resulting diagram is zero because $S$ is uncappable.
There are only three values of $\beta$ we need to consider.
For each of these,
we compute the corresponding coefficient,
and easily evaluate the corresponding diagram.
The results are shown in Table~\ref{table:AB}.

\begin{table}[ht]
\center
\begin{tabular}{ c | c | c }
$\beta$ & $\mathrm{Coeff}_{f^{(2n+2)}}(\beta)$ & value of diagram \\
\hline
%
%\beginpgfgraphicnamed{\pathtotrunk diagrams/tikz/#1-external}%
\input{\pathtotrunk diagrams/tikz/TLn/iIiiIi.tex}%
%\endpgfgraphicnamed
 
     & $1$ \tablestrut
     & $\tr{\rho^{1/2}(S)^3}$ \\
%
%\beginpgfgraphicnamed{\pathtotrunk diagrams/tikz/#1-external}%
\input{\pathtotrunk diagrams/tikz/TLn/uiInIi.tex}%
%\endpgfgraphicnamed

     & $(-1)^{n+1}\frac{[n+1]}{[2n+2]}$ \tablestrut
     & $\tr{S^3}$ \\
%
%\beginpgfgraphicnamed{\pathtotrunk diagrams/tikz/#1-external}%
\input{\pathtotrunk diagrams/tikz/TLn/iInIiu.tex}%
%\endpgfgraphicnamed

     & $(-1)^{n+1}\frac{[n+1]}{[2n+2]}$ \tablestrut
     & $\omega^{-1} \tr{S^3}$
\end{tabular}
\vspace{6pt}
\caption{The terms of $f^{(2n+2)}$ that contribute to $\langle A,B \rangle$. }
\label{table:AB}
\end{table}

Now take the sum over all $\beta$ in the table
of the coefficient times the value of the diagram.
Note that if $S^2 = f^{(n)}$ then $\tr{S^3} = 0$.
Thus the two non-identity values of $\beta$
either cancel or give zero.
The term where $\beta$ is the identity
gives the desired result.
\end{proof}

\begin{lem}
\label{lem:BB}
If Assumption~\ref{ass:a} holds and $S^2 = f^{(n)}$ then
$$\langle B,B \rangle = \frac{[n+1][2n+2]}{[n][n+2]}.$$
\end{lem}

\begin{proof}
Consider the two diagrams
$$\pairing{B}{B} =Z\left(\scalebox{0.8}{%
%\beginpgfgraphicnamed{\pathtotrunk diagrams/tikz/#1-external}%
\input{\pathtotrunk diagrams/tikz/InnerProducts/BB.tex}%
%\endpgfgraphicnamed
}\right), \quad
  \pairing{f^{(n+1)}}{B} = Z\left(\scalebox{0.8}{%
%\beginpgfgraphicnamed{\pathtotrunk diagrams/tikz/#1-external}%
\input{\pathtotrunk diagrams/tikz/InnerProducts/JWB.tex}%
%\endpgfgraphicnamed
}\right).$$
The second is clearly zero.
We will compare what happens to each of these diagrams
when we expand the copy of $f^{(2n+2)}$
into a linear combination of Temperley-Lieb diagrams.

Let $\beta$ be a Temperley-Lieb diagram
in the expansion of $f^{(2n+2)}$.
Suppose $\beta$ contains a cup
that connects endpoints number $i$ and $i+1$ at the top
for some $i \neq n+1$.
The corresponding diagram for $\langle B,B \rangle$ is zero
because the cup connects two strands from the same copy of $S$.
Similarly,
the corresponding diagram for $\langle f^{(n+1)},B \rangle$ is zero
because the cup connects two strands from the same side of $f^{(n+1)}$.
Thus both diagrams corresponding to $\beta$ are zero.

Now suppose $\beta$ contains a cup in the middle,
connecting endpoints number $n+1$ and $n+2$ at the top.
In the corresponding diagram for $\langle B,B \rangle$,
this cup produces a copy of $S^2$,
which we can replace with $f^{(n)}$.
For $\langle f^{(n+1)},B \rangle$,
this cup produces a partial trace of $f^{(n+1)}$,
which we can replace with $\frac{[n+2]}{[n+1]} f^{(n)}$.
Thus the two diagrams corresponding to $\beta$
differ only by a factor of $\frac{[n+2]}{[n+1]}$.

Finally, suppose $\beta$ is the identity diagram.
The corresponding diagram for $\langle B,B \rangle$
consists of four copies of $S$ arranged in a rectangle.
The left and right sides of this rectangle
consist of $n+1$ parallel strands.
We can replace one of these sides
with a partial trace of $f^{(n)}$.
Thus the diagram is equal to $\frac{[n+1]}{[n]} \tr{S^2}$.
The corresponding diagram for $\langle f^{(n+1)},B \rangle$
is equal to $\tr{S^2}$.

We can now evaluate
$$\langle B,B \rangle
 = \langle B,B \rangle
 - \frac{[n+1]}{[n+2]} \langle f^{(n+1)},B \rangle.$$
For both terms on the right hand side,
we express $f^{(2n+2)}$ as a linear combination
of Temperley-Lieb diagrams $\beta$.
For every $\beta$ except the identity,
these terms cancel.
The identity term gives the following.
$$\langle B,B \rangle
 = \frac{[n+1]}{[n]} \tr{S^2}
 - \frac{[n+1]}{[n+2]} \tr{S^2}.$$
The result now follows from a quantum integer identity, and that fact that $\tr{S^2} = [n+1]$.
\end{proof}

\begin{lem}
\label{lem:CC}
If Assumption~\ref{ass:a} holds
then
$$\langle C,C \rangle
= \frac {[2][2n+2]}{[2n+3][2n+4]} W_{2n+4} \langle A,A \rangle.$$
\end{lem}

\begin{proof}
We must evaluate the closed diagram
$$\pairing{C}{C}=Z\left(\scalebox{0.8}{%
%\beginpgfgraphicnamed{\pathtotrunk diagrams/tikz/#1-external}%
\input{\pathtotrunk diagrams/tikz/InnerProducts/CC.tex}%
%\endpgfgraphicnamed
}\right).$$
The proof is very similar to that of Lemma \ref{lem:AA}.
Indeed,
these are both special cases of a general recursive formula.
The idea is to apply Lemma~\ref{lem:yetanother} to the copy of $f^{(2n+4)}$,
and then evaluate each of the resulting diagrams.

Consider the first term.
Here,
$f^{(2n+4)}$ is replaced by a copy of $f^{(2n+3)}$
together with a single vertical strand on the right.
We can use sphericality to drag this strand over to the left.
This results in a partial trace of $f^{(2n+3)}$,
which is
$$ \frac{[2n+4]}{[2n+3]} f^{(2n+2)}.$$
Using Lemma~\ref{lem:AA},
we obtain
$$\frac{[2n+4]}{[2n+3]} \langle A,A \rangle.$$

Now consider the terms in the sum over $a$ and $b$.
Here,
$f^{(2n+4)}$ is replaced by a copy of $f^{(2n+2)}$
together with a ``cup'' and a ``cap''
in positions given by $a$ and $b$.
In most cases,
the resulting diagram is zero because $S$ is uncappable.
A cup at the leftmost position also gives zero
since it connects two strands coming from the top right of $f^{(2n+2)}$.
Similarly,
a cup in the rightmost position gives zero,
as does a cap in the leftmost or rightmost position.
The only cases that give a non-zero diagram are when $a,b \in \{2,2n+2\}$.
We obtain
$$\frac{1}{[2n+3][2n+4]}
( -[2]^2 - [2n+2]^2 - [2][2n+2]\omega - [2][2n+2]\omega^{-1})
\langle A,A \rangle.$$

The result now follows by adding the above two expressions
and expanding the quantum integers in terms of $q$.
\end{proof}

\begin{lem}\label{lem:DD}
If Assumption~\ref{ass:a} holds and $S^2 = f^{(n)}$
then
$$\langle D,D \rangle
= \frac{[n+1]^2[2n+2]}{[2][n]^2[2n+3]}([n+3]-[2][n]).$$
\end{lem}

\begin{proof}
We must evaluate the closed diagram
$$\pairing{D}{D}=Z\left(\scalebox{0.8}{%
%\beginpgfgraphicnamed{\pathtotrunk diagrams/tikz/#1-external}%
\input{\pathtotrunk diagrams/tikz/InnerProducts/DD.tex}%
%\endpgfgraphicnamed
}\right).$$
Consider the expansion of $f^{(2n+4)}$ into
a linear combination of Temperley-Lieb diagrams $\beta \in TL_{2n+4}$.
For each such $\beta$,
we compute the coefficient
and the value of the corresponding diagram.
There are twelve values of $\beta$ that give a non-zero diagram.
Many of these are reflections or rotations of each other.
They are shown in Table~\ref{table:DD}.

\begin{table}[ht]
\center
\begin{tabular}{ c | c | c }
$\beta$ & $\mathrm{Coeff}_{f^{(2n+4)}}(\beta)$ & value of diagram \\
\hline
%
%\beginpgfgraphicnamed{\pathtotrunk diagrams/tikz/#1-external}%
\input{\pathtotrunk diagrams/tikz/TLn/IiiiiiiI.tex}%
%\endpgfgraphicnamed

     & $1$ \tablestrut
     & $\left(\frac{[n+1]}{[n]}\right)^2 \tr{S^2}$ \\
%
%\beginpgfgraphicnamed{\pathtotrunk diagrams/tikz/#1-external}%
\input{\pathtotrunk diagrams/tikz/TLn/IiuniiiI.tex}%
%\endpgfgraphicnamed
, %
%\beginpgfgraphicnamed{\pathtotrunk diagrams/tikz/#1-external}%
\input{\pathtotrunk diagrams/tikz/TLn/IiiiuniI.tex}%
%\endpgfgraphicnamed
 
     & $- \frac{[n+1] [n+3]} {[2n+4]}$ \tablestrut
     & $\frac{[n+1]}{[n]} \tr{S^4}$ \\
%
%\beginpgfgraphicnamed{\pathtotrunk diagrams/tikz/#1-external}%
\input{\pathtotrunk diagrams/tikz/TLn/IUNiiI.tex}%
%\endpgfgraphicnamed
, %
%\beginpgfgraphicnamed{\pathtotrunk diagrams/tikz/#1-external}%
\input{\pathtotrunk diagrams/tikz/TLn/IiiUNI.tex}%
%\endpgfgraphicnamed
  
     & $\frac{[n] [n+1] [n+2] [n+3]} {[2] [2n+3] [2n+4]}$ \tablestrut
     & $\left(\frac{[n+1]}{[n]}\right)^2 \tr{S^2}$ \\
%
%\beginpgfgraphicnamed{\pathtotrunk diagrams/tikz/#1-external}%
\input{\pathtotrunk diagrams/tikz/TLn/IUiiNI.tex}%
%\endpgfgraphicnamed
, %
%\beginpgfgraphicnamed{\pathtotrunk diagrams/tikz/#1-external}%
\input{\pathtotrunk diagrams/tikz/TLn/INiiUI.tex}%
%\endpgfgraphicnamed

     & $\frac{[n]^2 [n+1]^2} {[2] [2n+3] [2n+4]}$ \tablestrut
     & $\left(\frac{[n+1]}{[n]}\right)^2 \tr{S^2}$ \\
%
%\beginpgfgraphicnamed{\pathtotrunk diagrams/tikz/#1-external}%
\input{\pathtotrunk diagrams/tikz/TLn/IUinniI.tex}%
%\endpgfgraphicnamed
, %
%\beginpgfgraphicnamed{\pathtotrunk diagrams/tikz/#1-external}%
\input{\pathtotrunk diagrams/tikz/TLn/IinniUI.tex}%
%\endpgfgraphicnamed
,
%
%\beginpgfgraphicnamed{\pathtotrunk diagrams/tikz/#1-external}%
\input{\pathtotrunk diagrams/tikz/TLn/INiuuiI.tex}%
%\endpgfgraphicnamed
, %
%\beginpgfgraphicnamed{\pathtotrunk diagrams/tikz/#1-external}%
\input{\pathtotrunk diagrams/tikz/TLn/IiuuiNI.tex}%
%\endpgfgraphicnamed

     & $- \frac{[n] [n+1]^2 [n+2]} {[2n+3] [2n+4]} $ \tablestrut
     & $ \frac{[n+1]}{[n]} \tr{S^4}$\\
%
%\beginpgfgraphicnamed{\pathtotrunk diagrams/tikz/#1-external}%
\input{\pathtotrunk diagrams/tikz/TLn/IiununiI.tex}%
%\endpgfgraphicnamed

     & $ \frac{[2] [n+1]^2 [n+2]^2} {[2n+3] [2n+4]} $ \tablestrut
     & $ \tr{S^6}$
\end{tabular}
\vspace{6pt}
\caption{The terms of $f^{(2n+4)}$ that contribute to $\langle D,D \rangle$. }
\label{table:DD}
\end{table}

Now take the sum over all $\beta$ in the table
of the coefficient times the value of the diagram.
Since $S^2 = f^{(n)}$,
we have
$$\tr{S^2} = \tr{S^4} = \tr{S^6}.$$
\end{proof}

\begin{lem}\label{lem:CD}
If Assumption~\ref{ass:a} holds and $S^2 = f^{(n)}$ and $\omega = -1$
then
\begin{equation*}
%\label{eq:CD}
\langle C,D \rangle =
    Z(\Gamma)
    + \frac {2}{[n]}
    + \frac{1}{[2n+3]} \langle A,A \rangle.
\end{equation*}
\end{lem}

\begin{proof}
First we prove a formula for $D$.
\begin{equation}
\label{eq:Dexpandone}
D = \scalebox{0.8}{%
%\beginpgfgraphicnamed{\pathtotrunk diagrams/tikz/#1-external}%
\input{\pathtotrunk diagrams/tikz/InnerProducts/iD.tex}%
%\endpgfgraphicnamed
}.
\end{equation}
Apply a left to right reflection of Lemma~\ref{lem:basicrecursion}
to the copy of $f^{(2n+4)}$ in $D$.
The first term gives the desired diagram.
Now consider a term in the sum over $a$.
This contains a cup in a position given by $a$.
For all but two values of $a$,
this cup connects two strands from the same copy of $S$,
so the resulting diagram is zero.
The remaining two values of $a$ are $n+1$ and $n+3$.
For each of these,
the cup connects two different copies of $S$,
giving rise to a copy of $S^2$.
We can replace this with $f^{(n)}$,
which is a linear combination of Temperley-Lieb diagrams.
But any such Temperley-Lieb diagram 
gives rise to a cup connected to the top edge of $f^{(2n+3)}$,
and thus gives zero.
This completes the derivation of Equation~\eqref{eq:Dexpandone}.

Next we prove the following.
\begin{equation}
\label{eq:Dexpandtwo}
D = \scalebox{0.8}{%
%\beginpgfgraphicnamed{\pathtotrunk diagrams/tikz/#1-external}%
\input{\pathtotrunk diagrams/tikz/InnerProducts/iDi.tex}%
%\endpgfgraphicnamed
}
  + \frac{1}{[2n+3]}
    \scalebox{0.8}{%
%\beginpgfgraphicnamed{\pathtotrunk diagrams/tikz/#1-external}%
\input{\pathtotrunk diagrams/tikz/InnerProducts/Sflower.tex}%
%\endpgfgraphicnamed
}.
\end{equation}
To prove this,
apply Lemma~\ref{lem:basicrecursion}
to the copy of $f^{(2n+3)}$ in Equation~\eqref{eq:Dexpandone}.
The first term in the expansion
gives the first term in the desired expression for $D$.

It remains to show that the sum over $a$
is equal to the second term in the desired expression.
Consider a term for $a \not \in \{n,n+2\}$.
The cup connects two strands from the same copy of $S$,
giving zero.

Consider the term corresponding to $a = n+2$.
The position of the cup is such that
the right two copies of $S$ are connected by $n$ strands.
This is a copy of $S^2$,
which is equal to $f^{(n)}$,
which in turn is a linear combination of Temperley-Lieb diagrams.
Any such Temperley-Lieb diagram
results in a cap connected to the top edge of $f^{(2n+2)}$,
giving zero.

Now consider the term corresponding to $a = n$.
The coefficient of this term is
$$(-1)^{n+1} \frac{[n]}{[2n+3]}.$$
The left two copies of $S$ form a copy of $S^2$,
which is equal to a sideways copy of $f^{(n)}$,
which in turn we express as
a linear combination of Temperley-Lieb diagrams $\beta$.
Every such $\beta$ gives zero except
$$\beta = 
\begin{tikzpicture}[TL12]
        \draw (1,-2) arc (180:0:0.5);
        \draw (3,-2)--(1,2);
        \node at (4,0) {...};
        \draw (7,-2)--(5,2);
        \draw (6,2) arc (-180:0:0.5);
\end{tikzpicture} \, ,$$
which has coefficient
$$(-1)^{n+1} \frac{1}{[n]}$$
and gives the second diagram in the desired expression for $D$.
The total coefficient of this diagram
is the product of the above coefficients
for the term $a$ and the diagram $\beta$.
This completes the derivation of Equation~\eqref{eq:Dexpandtwo}.

Now we return to our computation of $\langle C,D \rangle$.
We must evaluate the expression
$$Z\left(\scalebox{0.8}{%
%\beginpgfgraphicnamed{\pathtotrunk diagrams/tikz/#1-external}%
\input{\pathtotrunk diagrams/tikz/InnerProducts/CD.tex}%
%\endpgfgraphicnamed
}\right).$$
Apply Equation~\eqref{eq:Dexpandtwo},
upside down,
to the bottom half of this diagram.
For the last term of this equation,
apply sphericality
and use Lemma~\ref{lem:AA} to reverse the shading.
We obtain the term
$$\frac{1}{[2n+3]} \langle A,A \rangle.$$
The first term from the equation gives
\begin{equation}
Z\left(\scalebox{0.8}{%
%\beginpgfgraphicnamed{\pathtotrunk diagrams/tikz/#1-external}%
\input{\pathtotrunk diagrams/tikz/InnerProducts/iCDi.tex}%
%\endpgfgraphicnamed
}\right).
\label{eq:iCDi}
\end{equation}

We expand $f^{(2n+2)}$ into a linear combination of
Temperley-Lieb diagrams $\beta$.
There are five values of $\beta$ we need to consider.
These are shown in Table~\ref{table:CD}.

\begin{table}[ht]
\center
\begin{tabular}{ c | c | c}
$\beta$ & $\mathrm{Coeff}_{f^{(2n+2)}}(\beta)$ & value of diagram \\
\hline
%
%\beginpgfgraphicnamed{\pathtotrunk diagrams/tikz/#1-external}%
\input{\pathtotrunk diagrams/tikz/TLn/iIiiIi.tex}%
%\endpgfgraphicnamed

     & $1$ \tablestrut
     & $Z(\Gamma)$ \\
%
%\beginpgfgraphicnamed{\pathtotrunk diagrams/tikz/#1-external}%
\input{\pathtotrunk diagrams/tikz/TLn/uIniIi.tex}%
%\endpgfgraphicnamed

     & $(-1)^n \frac{[n+2]}{[2n+2]}$ \tablestrut
     & $(-1)^{n+1} \frac{1}{[n]} \omega \tr{S^2}$ \\
%
%\beginpgfgraphicnamed{\pathtotrunk diagrams/tikz/#1-external}%
\input{\pathtotrunk diagrams/tikz/TLn/iIinIu.tex}%
%\endpgfgraphicnamed

     & $(-1)^n \frac{[n+2]}{[2n+2]}$ \tablestrut
     & $(-1)^{n+1} \frac{1}{[n]} \omega^{-1} \tr{S^2}$ \\
%
%\beginpgfgraphicnamed{\pathtotrunk diagrams/tikz/#1-external}%
\input{\pathtotrunk diagrams/tikz/TLn/IniiIu.tex}%
%\endpgfgraphicnamed
, %
%\beginpgfgraphicnamed{\pathtotrunk diagrams/tikz/#1-external}%
\input{\pathtotrunk diagrams/tikz/TLn/uIiinI.tex}%
%\endpgfgraphicnamed

     & $(-1)^n \frac{[n]}{[2n+2]}$ \tablestrut
     & $(-1)^{n+1} \frac{1}{[n]}  \tr{S^2}$
\end{tabular}
\vspace{6pt}
\caption{The terms of $f^{(2n+2)}$ that contribute to \eqref{eq:iCDi}.}
\label{table:CD}
\end{table}

Now take the sum over all $\beta$ in the table
of the coefficient times the value of the diagram.
\end{proof}

The following inner products involving $B'$
will help us to evaluate $Z(\Gamma)$.

\begin{lem}
\label{lem:ABprime}
If Assumption~\ref{ass:a} holds and $S^2 = f^{(n)}$ and $\omega = -1$
then
$$\langle A,B' \rangle
   = \frac{[n-1]}{[n+1]} \tr{\rho^{1/2}(S)^3}$$
\end{lem}

\begin{proof}
We must evaluate the closed diagram
$$Z\left( \scalebox{0.8}{%
%\beginpgfgraphicnamed{\pathtotrunk diagrams/tikz/#1-external}%
\input{\pathtotrunk diagrams/tikz/InnerProducts/ABprime.tex}%
%\endpgfgraphicnamed
} \right).$$
The proof is very similar to that of Lemma~\ref{lem:AB},
so we will omit the details.
The relevant table is shown in Table~\ref{table:ABprime}.

\begin{table}[ht]
\center
\begin{tabular}{ c | c | c }
$\beta$ & $\mathrm{Coeff}_{f^{(2n+2)}}(\beta)$ & value of diagram \\
\hline
%
%\beginpgfgraphicnamed{\pathtotrunk diagrams/tikz/#1-external}%
\input{\pathtotrunk diagrams/tikz/TLn/uinbIibI.tex}%
%\endpgfgraphicnamed

     & $\frac{[2n]}{[2n+2]} $ \tablestrut
     & $\tr{\rho^{1/2}(S)^3}$ \\
%
%\beginpgfgraphicnamed{\pathtotrunk diagrams/tikz/#1-external}%
\input{\pathtotrunk diagrams/tikz/TLn/inbIibIu.tex}%
%\endpgfgraphicnamed

     & $\frac{[2]}{[2n+2]}$ \tablestrut
     & $\omega^{-1} \tr{\rho^{1/2}(S)^3}$ \\
\end{tabular}
\vspace{6pt}
\caption{The terms of $f^{(2n+2)}$ that contribute to $\langle A,B' \rangle$. }
\label{table:ABprime}
\end{table}

\end{proof}

\begin{lem}
\label{lem:BBprime}
If Assumption~\ref{ass:a} holds and $S^2 = f^{(n)}$ and $\omega=-1$
then
$$\langle B,B' \rangle = Z(\Gamma) + \frac{[2][n+1]}{[n][n+2]}.$$
\end{lem}

\begin{proof}
We must evaluate the closed diagram
$$ Z\left( \scalebox{0.8}{%
%\beginpgfgraphicnamed{\pathtotrunk diagrams/tikz/#1-external}%
\input{\pathtotrunk diagrams/tikz/InnerProducts/BBprime.tex}%
%\endpgfgraphicnamed
} \right) .$$
Inspired by the proof of Lemma \ref{lem:BB},
we observe that
$$\langle B,B' \rangle
  = \langle B,B' \rangle
  - \frac{[n+1]}{[n+2]} \langle f^{(n+1)},B' \rangle.$$
We expand the copies of $f^{(2n+2)}$ on the right hand side.
By the same argument as for Lemma \ref{lem:BB},
all terms will cancel except for those coming from the identity diagram.

If we replace $f^{(2n+2)}$ by the identity in $\langle B,B' \rangle$
then we obtain $\Gamma$.
If we replace $f^{(2n+2)}$ by the identity in $\langle f^{(n+1)},B' \rangle$
then we obtain a diagram containing
two copies of $S$ and one copy of $f^{(n+1)}$.
We must now expand $f^{(n+1)}$ as a linear combination
of Temperley-Lieb diagrams $\beta$.
For all but one such diagram $\beta$,
the resulting diagram is zero because $S$ is uncappable.
The only diagram we need to consider is
$$\beta =
\begin{tikzpicture}[TL12]
        \draw (0,-3) arc (180:0:1.5);  
        \draw (1,-3) arc (180:0:0.5);
        \draw (4,-3)--(1,3);
        \node at (4.5,0) {...};
        \draw (8,-3)--(5,3);
        \draw (6,3) arc (-180:0:1.5);
        \draw (7,3) arc (-180:0:0.5);
\end{tikzpicture} \, ,$$
which has coefficient $\frac{[2]}{[n][n+1]}$
and gives the diagram $\omega^{-1} \tr{S^2}$.
\end{proof}

%-----------------------------------
\subsection{Proving relations}\label{sec:provingrelations}

We now use our inner products,
together with Assumptions \ref{ass:a} and \ref{ass:b},
to prove that the required relations hold.
Note that the assumption $\omega=-1$ implies
$$W_{2m} = \left( \frac{[2m]}{[m]} \right)^2.$$

\begin{lem}
\label{lem:multiplication}
If Assumptions \ref{ass:a} and \ref{ass:b} hold
then $S^2 = f^{(n)}$.
\end{lem}

\begin{proof}
The relevant inner products are as follows.
\begin{itemize}
\item $\langle S^2,S^2 \rangle = \tr{S^4} = [n+1]$,
\item $\langle S^2,f^{(n)} \rangle = \tr{S^2} = [n+1]$,
\item $\langle f^{(n)},f^{(n)} \rangle = \tr{f^{(n)}} = [n+1]$.
\end{itemize}
The inner product of $S^2-f^{(n)}$ with itself is zero,
and the result follows from the assumption that $\cP$ is positive definite.
\end{proof}

\begin{lem}
\label{lem:ABrel}
If Assumptions \ref{ass:a} and \ref{ass:b} hold
then $A = i \frac{\sqrt{[n][n+2]}}{[n+1]} B$.
\end{lem}

\begin{proof}
By Lemmas \ref{lem:AA}, \ref{lem:AB}, \ref{lem:BB},
and our values for the moments,
we have the following.
\begin{itemize}
\item $\langle A,A \rangle = \frac{[2n+2]}{[n+1]}$,
\item $\langle A,B \rangle = i\frac{[2n+2]}{\sqrt{[n][n+2]}}$,
\item $\langle B,B \rangle = \frac{[n+1][2n+2]}{[n][n+2]}$.
\end{itemize}
Thus
$$\langle A,A \rangle \langle B,B \rangle = |\langle A,B \rangle|^2.$$
Thus $A$ and $B$ are linearly dependent.
The precise relation is then
$$A = \frac{\langle A,B \rangle}{\langle B,B \rangle}B.$$
\end{proof}

\begin{lem}
If Assumptions \ref{ass:a} and \ref{ass:b} hold
then
$$Z(\Gamma) = \frac{[n-1][2n+2] - [2][n+1]}{[n][n+2]}.$$
\end{lem}

\begin{proof}
By Lemmas \ref{lem:ABprime} and \ref{lem:BBprime},
\begin{itemize}
\item $\langle A,B' \rangle
       = \frac{[n-1]}{[n+1]} i\frac{[2n+2]}{\sqrt{[n][n+2]}} $.
\item $\langle B,B' \rangle = Z(\Gamma) + \frac{[2][n+1]}{[n][n+2]}.$
\end{itemize}

By Lemma \ref{lem:ABrel},
$A = i \frac{\sqrt{[n][n+2]}}{[n+1]} B$.
Thus
$$\langle A,B' \rangle
   = i \frac{\sqrt{[n][n+2]}}{[n+1]} \langle B,B' \rangle.$$
The result follows by solving for $Z(\Gamma)$.
\end{proof}

\begin{lem}\label{lem:CDrel}
If Assumptions \ref{ass:a} and \ref{ass:b} hold then
$$C = \frac{[2][2n+4]}{[n+1][n+2]}D.$$
\end{lem}

\begin{proof}
By Lemmas \ref{lem:CC}, \ref{lem:CD}, \ref{lem:DD}
and our values for $Z(\Gamma)$ and the moments,
we have the following.
\begin{itemize}
\item $\langle C,C \rangle =
       \frac{[2][2n+2]^2[2n+4]}{[n+1][n+2]^2[2n+3]}$.
\item $\langle C,D \rangle = \frac{[2n+2]^2}{[n+2][2n+3]}$.
\item $\langle D,D \rangle = \frac{[n+1]^2[2n+2]}{[2][n]^2[2n+3]}
       ([n+3]-[2][n])$.
\end{itemize}
Here,
we have used quantum integer identities to simplify the expression for $\langle C,D \rangle$.

For arbitrary $n$, $m$, and $q$,
$$[n+m] = \frac{1}{[4]}([4-m][n] + [m][n+4]),$$
and
$$[2m] = [m]([m+1]-[m-1]).$$
By Lemma \ref{lem:quantum-identity} and the assumption $[2] = d_k$,
$$[n+4] = [3][n].$$
(This is the only time we use the assumption $[2]=d_k$.)
We can now express each of our inner products
in terms of $[n]$, $[2]$, $[3]$, and $[4]$.
After some computation we find that
$$\langle C,C \rangle \langle D,D \rangle = |\langle C,D \rangle|^2.$$
Thus $C$ and $D$ are linearly dependent.
The precise relation is then
$$C = \frac{\langle C,C \rangle}{\langle C,D \rangle} D.$$
\end{proof}

%% file: diagrams/tikz/TLn/iIiiIi.tex
\begin{tikzpicture}[TL12]
    \foreach \x in {1,3,4,6} \draw (\x cm , -1.5cm)--(\x cm, 1.5cm);
    \foreach \x in {2,5} \draw[ultra thick] (\x cm , -1.5cm)--(\x cm, 1.5cm);
\end{tikzpicture}

%% file: diagrams/tikz/InnerProducts/Tet.tex
\begin{tikzpicture}[STrain]
	\fill[shaded] (-1,0) -- (1,0) -- (0,1) -- cycle;
	\draw (-1,0) -- (1,0);
	\draw (0,0) -- (0,1);
	\foreach \x in {-1,1} {
		\draw (\x,0) -- (0,1);
	}
	\draw (-1,0) .. controls (-1,2) and (1,2) .. (1,0);
	\drawS{0}{1}{0}
	\foreach \x in {-1,0,1} {
		\drawS{\x}{0}{-90}
	}
	\node at (0.1,0.5) {\footnotesize{$2$}};
	\node at (0.9,1) {\footnotesize{$2$}};
	\foreach \x in {-1,1} {
		\node at (0.5*\x,-0.1) {\footnotesize{$n-1$}};
		\node at (0.6*\x,0.7) {\footnotesize{$n-1$}};
	}
\end{tikzpicture}

%% file: diagrams/tikz/InnerProducts/Bprime.tex
\begin{tikzpicture}[STrain]
	\STrainOver \STrainStrings{$n+1$}{$n-1$} \STrainOne
        \draw (0,0)--(0,-0.5);
        \node[anchor=west] at (0,-0.35) {\footnotesize$2n+2$};
	\JWPlusTwo
\end{tikzpicture}

%% file: diagrams/tikz/InnerProducts/AA.tex
\begin{tikzpicture}[STrain]
	\RainbowOne
	\upsidedown{\RainbowOne}
	\JWPlusTwo
\end{tikzpicture}

%% file: diagrams/tikz/InnerProducts/AB.tex
\begin{tikzpicture}[STrain]
	\RainbowOne
	\upsidedown{\STrainStrings{$n+1$}{$n+1$} \STrainOne}
	\JWPlusTwo
\end{tikzpicture}

%% file: diagrams/tikz/TLn/uiInIi.tex
\begin{tikzpicture}[TL12]
        \draw (1,1.5) arc (-180:0:0.5);
        \draw (1,-1.5) -- (3,1.5);
        \draw [ultra thick] (2,-1.5) -- (4,1.5);
        \draw (3,-1.5) arc (180:0:0.5);
        \draw [ultra thick] (5,-1.5) -- (5,1.5);
        \draw (6,-1.5) -- (6,1.5);
\end{tikzpicture}

%% file: diagrams/tikz/TLn/iInIiu.tex
\begin{tikzpicture}[TL12]
        \draw (1,-1.5) -- (1,1.5);
        \draw [ultra thick] (2,-1.5) -- (2,1.5);
        \draw (3,-1.5) arc (180:0:0.5);
        \draw [ultra thick] (5,-1.5) -- (3,1.5);
        \draw (6,-1.5) -- (4,1.5);
        \draw (5,1.5) arc (-180:0:0.5);
\end{tikzpicture}

%% file: diagrams/tikz/InnerProducts/BB.tex
\begin{tikzpicture}[STrain]
	\STrainStrings{$n+1$}{$n+1$} \STrainOne
	\upsidedown{\STrainStrings{$n+1$}{$n+1$} \STrainOne}
	\JWPlusTwo
\end{tikzpicture}

%% file: diagrams/tikz/InnerProducts/JWB.tex
\begin{tikzpicture}[STrain]
	\filldraw[shaded] (-0.7,0.2) arc (180:90:0.5 and 0.8) --
                          (0.2,1) arc (90:0:0.5 and 0.8) -- cycle;
	\upsidedown{\STrainStrings{$n+1$}{$n+1$} \STrainOne}
	\JWPlusTwo;
        \filldraw[fill=white,thick] (-0.2,0.5) rectangle (0.2,1.5);
        \node[rotate=-90] at (0,1) {$f^{(n+1)}$};
        \node[anchor=south west] at (-0.7,0.2) {\footnotesize$n+1$};
        \node[anchor=south west] at (0.7,0.2) {\footnotesize$n+1$};
\end{tikzpicture}

%% file: diagrams/tikz/InnerProducts/CC.tex
\begin{tikzpicture}[STrain]
	\RainbowTwo
	\upsidedown{\RainbowTwo}
	\JWPlusFour
\end{tikzpicture}

%% file: diagrams/tikz/InnerProducts/DD.tex
\begin{tikzpicture}[STrain]
	\STrainThreeStrings{$n+1$}{$2$}{$n+1$} \STrainOneOne
	\upsidedown{\STrainThreeStrings{$n+1$}{$2$}{$n+1$} \STrainOneOne}
	\JWPlusFour
\end{tikzpicture}

%% file: diagrams/tikz/TLn/IiiiiiiI.tex
\begin{tikzpicture}[TL12]
     \foreach \x in {2,...,7} \draw (\x cm , -1.5cm)--(\x cm, 1.5cm);
     \foreach \x in {1,8} \draw[ultra thick] (\x cm , -1.5cm)--(\x cm, 1.5cm);
\end{tikzpicture}

%% file: diagrams/tikz/TLn/IiuniiiI.tex
\begin{tikzpicture}[TL12]
        \draw[ultra thick] (1,-1.5)--(1,1.5);
        \draw (2,-1.5)--(2,1.5);
        \draw (3,1.5) arc (-180:0:.5);
        \draw (3,-1.5) arc (180:0:.5);
        \draw (5,-1.5)--(5,1.5);
        \draw (6,-1.5)--(6,1.5);
        \draw (7,-1.5)--(7,1.5);
        \draw[ultra thick] (8,-1.5)--(8,1.5);
\end{tikzpicture}

%% file: diagrams/tikz/TLn/IiiiuniI.tex
\begin{tikzpicture}[TL12]
        \draw[ultra thick] (1,-1.5)--(1,1.5);
        \draw (2,-1.5)--(2,1.5);
        \draw (3,-1.5)--(3,1.5);
        \draw (4,-1.5)--(4,1.5);
        \draw (5,1.5) arc (-180:0:.5);
        \draw (5,-1.5) arc (180:0:.5);
        \draw (7,-1.5)--(7,1.5);
        \draw[ultra thick] (8,-1.5)--(8,1.5);
\end{tikzpicture}

%% file: diagrams/tikz/TLn/IUNiiI.tex
\begin{tikzpicture}[TL12]
        \draw[ultra thick] (1,-1.5)--(1,1.5);
        \draw (2,1.5) arc (-180:0: 1.5 and 1.1);
        \draw (2,-1.5) arc (180:0: 1.5 and 1.1);
        \draw (3,1.5) arc (-180:0:.5);
        \draw (3,-1.5) arc (180:0:.5);
        \draw (6,-1.5)--(6,1.5);
        \draw (7,-1.5)--(7,1.5);
        \draw[ultra thick] (8,-1.5)--(8,1.5);
\end{tikzpicture}

%% file: diagrams/tikz/TLn/IiiUNI.tex
\begin{tikzpicture}[TL12]
        \draw[ultra thick] (1,-1.5)--(1,1.5);
        \draw (2,-1.5)--(2,1.5);
        \draw (3,-1.5)--(3,1.5);
        \draw (4,1.5) arc (-180:0: 1.5 and 1.1);
        \draw (4,-1.5) arc (180:0: 1.5 and 1.1);
        \draw (5,1.5) arc (-180:0:.5);
        \draw (5,-1.5) arc (180:0:.5);
        \draw[ultra thick] (8,-1.5)--(8,1.5);
\end{tikzpicture}

%% file: diagrams/tikz/TLn/IUiiNI.tex
\begin{tikzpicture}[TL12]
        \draw[ultra thick] (1,-1.5)--(1,1.5);
        \draw (2,1.5) arc (-180:0: 1.5 and 1.1);
        \draw (3,1.5) arc (-180:0:.5);
        \draw (4,-1.5)--(6,1.5);
        \draw (5,-1.5)--(7,1.5);
        \draw (6,-1.5) arc (180:0: 1.5 and 1.1);
        \draw (7,-1.5) arc (180:0:.5);
        \draw[ultra thick] (10,-1.5)--(10,1.5);
\end{tikzpicture}

%% file: diagrams/tikz/TLn/INiiUI.tex
\begin{tikzpicture}[TL12]
        \draw[ultra thick] (1,-1.5)--(1,1.5);
        \draw (2,-1.5) arc (180:0: 1.5 and 1.1);
        \draw (3,-1.5) arc (180:0:.5);
        \draw (6,-1.5)--(4,1.5);
        \draw (7,-1.5)--(5,1.5);
        \draw (6,1.5) arc (-180:0: 1.5 and 1.1);
        \draw (7,1.5) arc (-180:0:.5);
        \draw[ultra thick] (10,-1.5)--(10,1.5);
\end{tikzpicture}

%% file: diagrams/tikz/TLn/IUinniI.tex
\begin{tikzpicture}[TL12]
        \draw[ultra thick] (1,-1.5)--(1,1.5);
        \draw (2,1.5) arc (-180:0: 1.5 and 1.1);
        \draw (3,1.5) arc (-180:0:.5);
        \draw (2,-1.5) .. controls (2,-1) and (6,0.5) ..(6,1.5);
        \draw (3,-1.5) arc (180:0:.5);
        \draw (5,-1.5) arc (180:0:.5);
        \draw (7,-1.5)--(7,1.5);
        \draw[ultra thick] (8,-1.5)--(8,1.5);
\end{tikzpicture}

%% file: diagrams/tikz/TLn/IinniUI.tex
\begin{tikzpicture}[TL12]
        \draw[ultra thick] (1,-1.5)--(1,1.5);
        \draw (2,-1.5)--(2,1.5);
        \draw (3,-1.5) arc (180:0:.5);
        \draw (5,-1.5) arc (180:0:.5);
        \draw (7,-1.5) .. controls (7,-1) and (3,0.5) ..(3,1.5);
        \draw (4,1.5) arc (-180:0: 1.5 and 1.1);
        \draw (5,1.5) arc (-180:0:.5);
        \draw[ultra thick] (8,-1.5)--(8,1.5);
\end{tikzpicture}

%% file: diagrams/tikz/TLn/INiuuiI.tex
\begin{tikzpicture}[TL12]
        \draw[ultra thick] (1,-1.5)--(1,1.5);
        \draw (2,-1.5) arc (180:0: 1.5 and 1.1);
        \draw (3,-1.5) arc (180:0:.5);
        \draw (2,1.5) .. controls (2,1) and (6,-0.5) ..(6,-1.5);
        \draw (3,1.5) arc (-180:0:.5);
        \draw (5,1.5) arc (-180:0:.5);
        \draw (7,-1.5)--(7,1.5);
        \draw[ultra thick] (8,-1.5)--(8,1.5);
\end{tikzpicture}

%% file: diagrams/tikz/TLn/IiuuiNI.tex
\begin{tikzpicture}[TL12]
        \draw[ultra thick] (1,-1.5)--(1,1.5);
        \draw (2,-1.5)--(2,1.5);
        \draw (3,1.5) arc (-180:0:.5);
        \draw (5,1.5) arc (-180:0:.5);
        \draw (7,1.5) .. controls (7,1) and (3,-0.5) ..(3,-1.5);
        \draw (4,-1.5) arc (180:0: 1.5 and 1.1);
        \draw (5,-1.5) arc (180:0:.5);
        \draw[ultra thick] (8,-1.5)--(8,1.5);
\end{tikzpicture}

%% file: diagrams/tikz/TLn/IiununiI.tex
\begin{tikzpicture}[TL12]
        \draw[ultra thick] (1,-1.5)--(1,1.5);
        \draw (2,-1.5)--(2,1.5);
        \draw (3,1.5) arc (-180:0:.5);
        \draw (3,-1.5) arc (180:0:.5);
        \draw (5,1.5) arc (-180:0:.5);
        \draw (5,-1.5) arc (180:0:.5);
        \draw (7,-1.5)--(7,1.5);
        \draw[ultra thick] (8,-1.5)--(8,1.5);
\end{tikzpicture}

%% file: diagrams/tikz/InnerProducts/iD.tex
\begin{tikzpicture}[STrain]
        \filldraw[shaded] (-1.5,-0.5) .. controls (-1.5,0.5) .. (-1,1)
                          -- (0,1) -- (0,-0.5);
        \node[anchor=west] at (0,-0.35) {\footnotesize$2n+3$};
	\STrainThreeStrings{$n$}{$2$}{$n+1$} \STrainOneOne
	\JWwide{2n+3}
\end{tikzpicture}

%% file: diagrams/tikz/InnerProducts/iDi.tex
\begin{tikzpicture}[STrain]
        \filldraw[shaded] (-1.5,-0.5) .. controls (-1.5,0.5) .. (-1,1)
                          -- (1,1) .. controls (1.5,0.5) .. (1.5,-0.5);
	\STrainThreeStrings{$n$}{$2$}{$n$} \STrainOneOne
        \draw (0,0)--(0,-0.5);
        \node[anchor=west] at (0,-0.35) {\footnotesize$2n+2$};
	\JWwide{2n+2}
\end{tikzpicture}

%% file: diagrams/tikz/InnerProducts/Sflower.tex
\begin{tikzpicture}[STrain]
	\filldraw[shaded] (-1.2,-0.5) -- (-1.2,0.2) arc (180:0:0.2) -- 
                          (-0.8,0) -- (0.8,0) --
                          (0.8,0.2) arc (180:0:0.2) -- (1.2,-0.5);
	\draw (0,-0.5) -- (0,1);
	\node[anchor=west] at (0,0.5) {\footnotesize{$2n$}};
        \node[anchor=west] at (0,-0.35) {\footnotesize$2n+4$};
	\drawS{0}{1}{90}
	\JWPlusTwo
\end{tikzpicture}

%% file: diagrams/tikz/InnerProducts/CD.tex
\begin{tikzpicture}[STrain]
	\RainbowTwo
	\upsidedown{\STrainThreeStrings{$n+1$}{$2$}{$n+1$} \STrainOneOne}
	\JWPlusFour
\end{tikzpicture}

%% file: diagrams/tikz/InnerProducts/iCDi.tex
\begin{tikzpicture}[STrain]
   \filldraw[shaded] (-1,-1) .. controls (-1.5,-1) ..  (-1.5,0.2)
         arc (180:0:1.5) .. controls (1.5,-1) ..(1,-1)
         -- (1,0.5) arc (0:180:1)  -- cycle;
   \draw (0,0) -- (0,1);
   \drawS{0}{1}{180}
   \node at (0.15,0.5) {\footnotesize$2n$};
   \upsidedown{\STrainThreeStrings{$n$}{$2$}{$n$} \STrainOneOne}
   \JWwide{2n+2}
\end{tikzpicture}

%% file: diagrams/tikz/TLn/uIniIi.tex
\begin{tikzpicture}[TL12]
        \draw (1,1.5) arc (-180:0:0.5);
        \draw[ultra thick] (1,-1.5)--(3,1.5);
        \draw (2,-1.5) arc (180:0:0.5);
        \draw (4,-1.5)--(4,1.5);
        \draw[ultra thick] (5,-1.5)--(5,1.5);
        \draw (6,-1.5)--(6,1.5);
\end{tikzpicture}

%% file: diagrams/tikz/TLn/iIinIu.tex
\begin{tikzpicture}[TL12]
        \draw (1,-1.5)--(1,1.5);
        \draw[ultra thick] (2,-1.5)--(2,1.5);
        \draw (3,-1.5)--(3,1.5);
        \draw (4,-1.5) arc (180:0:0.5);
        \draw[ultra thick] (6,-1.5)--(4,1.5);
        \draw (5,1.5) arc (-180:0:0.5);
\end{tikzpicture}

%% file: diagrams/tikz/TLn/IniiIu.tex
\begin{tikzpicture}[TL12]
        \draw[ultra thick] (1,-1.5)--(1,1.5);
        \draw (2,-1.5) arc (180:0:0.5);
        \draw (4,-1.5)--(2,1.5);
        \draw (5,-1.5)--(3,1.5);
        \draw[ultra thick] (6,-1.5)--(4,1.5);
        \draw (5,1.5) arc (-180:0:0.5);
\end{tikzpicture}

%% file: diagrams/tikz/TLn/uIiinI.tex
\begin{tikzpicture}[TL12]
        \draw (1,1.5) arc (-180:0:0.5);
        \draw[ultra thick] (1,-1.5)--(3,1.5);
        \draw (2,-1.5)--(4,1.5);
        \draw (3,-1.5)--(5,1.5);
        \draw (4,-1.5) arc (180:0:0.5);
        \draw[ultra thick] (6,-1.5)--(6,1.5);
\end{tikzpicture}

%% file: diagrams/tikz/InnerProducts/ABprime.tex
\begin{tikzpicture}[STrain]
	\RainbowOne
	\upsidedown{ \STrainOver \STrainStrings{$n+1$}{$n-1$} \STrainOne}
	\JWPlusTwo
\end{tikzpicture}

%% file: diagrams/tikz/TLn/uinbIibI.tex
\begin{tikzpicture}[TL12]
        \draw (1,1.5) arc (-180:0:0.5);
        \draw (1,-1.5)--(3,1.5);
        \draw (2,-1.5) arc (180:0:0.5);
        \draw[ultra thick] (4,-1.5)--(4,1.5);
        \draw (5,-1.5)--(5,1.5);
        \draw[ultra thick] (6,-1.5)--(6,1.5);
\end{tikzpicture}

%% file: diagrams/tikz/TLn/inbIibIu.tex
\begin{tikzpicture}[TL12]
        \draw (1,-1.5)--(1,1.5);
        \draw (2,-1.5) arc (180:0:0.5);
        \draw[ultra thick] (4,-1.5)--(2,1.5);
        \draw (5,-1.5)--(3,1.5);
        \draw[ultra thick] (6,-1.5)--(4,1.5);
        \draw (5,1.5) arc (-180:0:0.5);
\end{tikzpicture}

%% file: diagrams/tikz/InnerProducts/BBprime.tex
\begin{tikzpicture}[STrain]
	\STrainOneStrings \STrainOne
	\upsidedown{ \STrainOver \STrainStrings{$n+1$}{$n-1$} \STrainOne}
	\JWPlusTwo
\end{tikzpicture}

%% file: text/propertiesofS.tex
In this section we construct an $8$-box $S \in GPA(\Ha{1}^p)$ that satisfies the hypotheses of Proposition \ref{prop:SgivesSPA}. The planar algebra generated by $S$ is the desired extended Haagerup planar algebra, thus completing the proof of Theorem \ref{thm:existence}.

We start with a brief description of how we found $S$,
since the definition of $S$ is not very enlightening on its own.
The goal was to find $S \in GPA(\Ha{1}^p)_8$ satisfying
the first three relations of Definition \ref{defn:relations},
which say that $\rho(S)=-S$, $S$ is uncappable, and $S^2=f^{(8)}$.
The dimension
of $GPA(\Ha{1}^p)_8$ is the number of loops of length $16$ 
based at even vertices of $\Ha{1}^p$, which is equal to $148375$.
We found the $19$-dimensional space of solutions to the first two relations,
then tried to solve the equation $S^2 = f^{(8)}$.
This one equation in the $8$-box space
is actually $148375$ separate equations over $\Complex$. We expect of course that there are many redundancies amongst these equations.

At this point the problem sounds quite tractable,
but we were still unable to solve it by general techniques.
We then used various ad hoc methods. First, we searched for quadratics that are perfect squares and solved those. This reduced the problem from $19$ variables to $9$.
We then chose a small collection of quadratics, corresponding to certain `extremal' loops in the basis, and found numerical approximations to a solution of these, using Newton's method. Approximating such a solution by algebraic numbers, we could then go back and check that all the quadratic equations are satisfied exactly.

\begin{rem}
Our solution $S$ need not be unique.
Although the subfactor planar algebra with principal graphs $\Ha{1}$ is unique,
there may be more than one way to embed it in its graph planar algebra.
Indeed, $-\overline{S}$ is also a solution,
which corresponds to applying the graph automorphism,
by Lemma \ref{lem:graphsymmetry}.
Due to the approximate nature of our search for $S$,
we cannot say whether there are any other solutions.
\end{rem}

The above description may sound daunting, but we manage to give a definition of $S$ that involves specifying only $21$ arbitrary looking numbers, and we reduce all the conditions we need to check on this element to computing certain powers of two matrices. The definition of $S$ is still somewhat overwhelming, and the verification of its properties is done by computer. Those intrepid readers who continue reading this section will have to read a short Mathematica program in order to fully verify some of the steps.

We will use the notation for vertices of $\Ha{1}^p$ given in Section \ref{sec:names}.
Let
$$\lambda = \sqrt{2-d^2},$$
$$C=-21516075 \lambda ^4+8115925 \lambda ^2+45255025.$$
For each $w \in \{0,1,2\}^8$ we will define an element $p_w \in \Integer[\lambda]$ below. 

%We now define the element $S \in GPA(\Ha{1})_{8,+}$, in Definition \ref{defn:S}, and state Lemma \ref{lem:Srelations} and Lemma \ref{lem:moments} giving the desired properties. The proofs of these statements will occupy the rest of the section.

\begin{defn}
\label{defn:S}
Suppose $\gamma$ is a loop of length $16$ in $\Ha{1}^p$.
Let the {\em collapsed loop} $\widehat{\gamma}$
be the sequence in  $\{0,1,2\}^8$
such that $\gamma_{2i-1}$ is in arm number $\widehat{\gamma}_i$ (in the notation of Section \ref{sec:names})
of $\Ha{1}^p$.
Further, define $\sigma(\gamma)$ to be $-1$ raised to the number of times
the vertices $v_0, w_0, z_i$ and $a_i$ appear in $\gamma$.
\begin{equation}
\label{eq:S}
S(\gamma) = C \sigma(\gamma) p_{\widehat{\gamma}} \frac{1}{\sqrt{d_{\gamma_1} d_{\gamma_9}}} \prod_{i=1}^{16} \frac{1}{\sqrt{d_{\gamma_i}}},
\end{equation} where the $p_{\widehat{\gamma}}$ are defined below.
\end{defn}

\begin{rem}
The main reason we write $S$ in the above form is that the $p_{\widehat{\gamma}}$ have much better number theoretic properties than the $S(\gamma)$.  In particular, the $p_{\widehat{\gamma}}$ all live in a degree $6$ extension of $\Rational$ while the $S(\gamma)$ live in a much larger number field.  As a consequence it is easier to do exact arithmetic using the $p_{\widehat{\gamma}}$.  There is a general reason for this phenomenon: the convention that subfactor planar algebras be spherical is not the best convention from the point of view of number theory.   A much more convenient convention is the ``lopsided" one where shaded circles count for $1$ while unshaded circles count for the index $d^2$.  Furthermore, this convention is also well-motivated from the perspective of subfactor theory where ${}_N M_M$ has as its left von Neumann dimension the index while its right von Neumann dimension is $1$.  These issues warrant further investigation (see \cite{1002.0168,gpa}).
\end{rem}

We make an apparently ad-hoc definition of $21$ elements of $\Integer[\lambda]$.
\begin{align*}
\input{text/polynomials}
\end{align*}
These elements are also defined in a \MMA notebook, available along with the sources for this article on the  \code{arXiv} (as the file \code{extra/code/Generator.nb}), or at \url{http://tqft.net/EH/notebook}. A PDF printout of the notebook is available by following this URL, then replacing \code{.nb} with \code{.pdf}. Everything that follows in this section is paralleled in the notebook, and in particular each of the statements below that requires checking some arithmetic has a corresponding test defined in the notebook.

%Notice that all these elements lie in the $\Integer$-lattice generated by $p_{00000001}$, $p_{00000011}$, $p_{00000111}$, $p_{00001011}$, $p_{00010111}$ and $p_{00101011}$. (You can verify this using the function \code{VerifyZ$\lambda$Lattice} in the \MMA notebook.) We only point this out for interest's sake; we don't use this fact later.

\begin{defnlem} We can consistently extend these definitions to every $p_w$ for $w \in \{0,1\}^8$ by the rules
\begin{align}
p_{abcdefgh} & = - p_{bcdefgha}, \label{eq:rot} \\
p_{abcdefgh} & = \overline{p}_{ahgfedcb}, \label{eq:reverse} \\
\intertext{and}
p_{00000000} & = 0, \label{eq:0s} \\
p_{abcd1111} & = 0. \label{eq:1s}
\end{align}
\end{defnlem}
\begin{proof}
For example, one can get from $p_{00110011}$ to $p_{01100110}$ either by rotating, or by reversing; fortunately $p_{00110011}$ is purely imaginary. Under the operations implicit in Equations \eqref{eq:rot} and \eqref{eq:reverse} each orbit in $\{0,1\}^8$ contains exactly one of the elements on which $p$ is defined above or in Equations \eqref{eq:0s} and \eqref{eq:1s}. The \MMA notebook provides functions \code{VerifyRotation} and \code{VerifyConjugation} to check that these rules hold uniformly.
\end{proof}

We further extend these definitions to every $p_w$ for $w \in \{0,1,2\}^8$ by the rules
\begin{align}
 p_{x 0 y} + p_{x 1 y} + p_{x 2 y} = 0. \label{eq:2s}
\end{align}

\begin{lem}
\label{lem:p2s} For every $abcd \in \{0,1,2\}^4$
\begin{align}
p_{abcd2222} = 0.
\end{align}
\end{lem}

\begin{proof}
This is a direct computation of $16$ cases for $abcd \in \{0,1\}^4$ using Equation \eqref{eq:2s}, after which the general case of $abcd \in \{0,1,2\}^4$ follows, again from \eqref{eq:2s}. The \MMA notebook provides a function \code{Verify2sVanish} that checks this Lemma.
\end{proof}

A final interesting note on the $p_w$: %(which one could use to reduce the explicit definitions above from $21$ elements to $10$):

\begin{lem}
\label{lem:graphsymmetry}
If $w'$ is obtained from $w$ by exchanging all $1$s and $2$s, then $p_w = - \overline{p_{w'}}$.
\end{lem}
\begin{proof}
A direct computation which you can verify using the \MMA function \code{VerifyGraphSymmetry}.
\end{proof}

This ends the definition of $S$. We now prepare to prove that it has the properties required to generate a subfactor planar algebra with principal graph $\Ha{1}$.

\begin{lem}\label{lem:Srelations}
The generator $S$ is self-adjoint,  has rotational eigenvalue $-1$ and is uncappable:
\begin{align*}
 S^* & = S, & \rho(S) & =-S, &\epsilon_i(S)&=0 \text{ for } i=1,\ldots, 2k.
\end{align*}
\end{lem}

\begin{lem}\label{lem:moments}
The generator $S$ and its ``one-click'' rotation $\rho^{1/2}(S)$ have the following moments:
\begin{align*}
\dtr{S^2} & =[9], \\
\dtr{S^3} & =0, \\
\dtr{S^4} & =[9], \\
\dtr{\rho^{1/2}(S)^3} & = i\frac{[18]}{\sqrt{[8][10]}}.
\end{align*}
Note that the scalars on the right sides of these equations actually refer to scalar multiples of the empty diagram in $GPA(\Ha{1}^p)_{0,\pm}$. In particular, each of them is a constant function on the even (or odd in the last case) vertices of the graph $\Ha{1}^p$.
\end{lem}

\begin{proof}[Proof of Lemma \ref{lem:Srelations}.]
Self-adjointness follows immediately from Equation \eqref{eq:reverse}. 

To show $\rho(S(\gamma))=-S(\gamma)$, first note that $\sigma(\gamma)$ and $\prod_{i=1}^{16} \frac{1}{\sqrt{d_{\gamma_i}}}$ are independent of the order of vertices appearing in 
$\gamma= \gamma_1 \cdots \gamma_{16} \gamma_1$. 
Thus, recalling Example \ref{eg:rotation} from \S \ref{sec:gpa},
\begin{align*}
\rho(S)(\gamma) 
%	& = \inputtikz{rho8Slabelled} \\ 
	& = \sqrt \frac{d_{\gamma_3}d_{\gamma_{11}}}{d_{\gamma_9}d_{\gamma_{1}}}
	     S(\gamma_3 \ldots \gamma_{16}\gamma_{1} \gamma_{2} \gamma_{3} ) \\
	& = \sqrt \frac{d_{\gamma_3}d_{\gamma_{11}}}{d_{\gamma_9}d_{\gamma_{1}}}
	     C \sigma(\gamma) 
	     p_{\widehat{\gamma_3 \ldots \gamma_{16}\gamma_{1} \gamma_{2} \gamma_{3}}} 
	     \frac{1}{\sqrt{ d_{\gamma_{11}} d_{\gamma_3} }}  
	     \prod_{i=1}^{16} \frac{1}{\sqrt{d_{\gamma_i}}} \\
	& =  C \sigma(\gamma) 
	     (- p_{\widehat{\gamma}}) 
	     \frac{1}{\sqrt{ d_{\gamma_{9}} d_{\gamma_1} }}  
	     \prod_{i=1}^{16} \frac{1}{\sqrt{d_{\gamma_i}}} \\
	& = -S(\gamma)	
\end{align*}
where we used Equation \eqref{eq:rot} in the second to last step.

Next consider $\epsilon_i(S)$, the result of attaching a cap on strands $i$ and $i+1$. Since we know $S$ is a rotational eigenvector, we only need to check $\epsilon_1(S) = \epsilon_2(S) = 0$. In particular, we don't need to explicitly treat the more complicated cases of $\epsilon_8(S)$ and $\epsilon_{16}(S)$, in which the cap is attached ``around the side'' of $S$, and the coefficients coming from critical points in the graph planar algebra are more complicated. 

Let $\Gamma_k$ be the set of length-$k$ loops on $\Ha{1}^p$, and $\gamma_i$ denote the $i^{th}$ vertex of 
$\gamma \in \Gamma_k$.  
The graph planar algebra formalism tells us that for $\varphi \in \Gamma_{14}$,
\begin{equation*}
\epsilon_i(S)(\varphi) = \sum_{\substack{\text{$\gamma \in \Gamma_{16}$ with} \\ \gamma_{j \leq i} = \varphi_j, \gamma_{i+2} = \varphi_i \\ \text{and $\gamma_{j \geq i+3} = \varphi_{j-2}$}}} \sqrt{\frac{d_{\gamma_{i+1}}}{d_{\gamma_i}}} S(\gamma)
\end{equation*}

We consider three cases, depending on whether the valence of $\varphi_i$ is $1$, $2$ or $3$.

If $\varphi_i $ has valence $1$, that is, it is an endpoint, then there is just one term in the sum: if $\gamma_i = v_0$, then $\gamma_{i+1} = w_0$, and if $\gamma_i = z_j$, then $\gamma_{i+1}=a_j$. In the  first case, the collapsed loop $\widehat{\gamma}$ must be $00000000$, so $S(\gamma)=0$ by Equation \eqref{eq:0s}. In the second case, if $\gamma_i = z_1$ then $\widehat{\gamma}$ must contain at least $4$ consecutive $1$s, so $S(\gamma)=0$ by Equation \eqref{eq:1s}. If $\gamma_i = z_2$, then $\widehat{\gamma}$ must contain at least $4$ consecutive $2$s, so $S(\gamma)=0$ by Lemma \ref{lem:p2s}.

If $\varphi_i$ has valence $2$, that is, it lies on one of the arms, then there are two terms in the sum, say $\gamma^+$ and $\gamma^-$. Moreover, the collapsed loops for the two terms are the same, and $\sigma(\gamma^+) = - \sigma(\gamma^-)$. Thus
\begin{align*}
\epsilon_i(S)(\varphi) & = 
C 
p_{\widehat{\gamma^\pm}}  
\frac{1}{\sqrt{d_{\gamma^\pm_1}d_{\gamma^\pm_9}}}
\sqrt{\prod_{j=1}^{14} \frac{1}{d_{\varphi_j}}} \times \\
& \qquad
\left( 
	\sigma(\gamma^+) 
		\frac{1}{\sqrt{d_{\gamma^+_{i+1}}d_{\gamma^+_{i+2}}}} 
		\sqrt{\frac{d_{\gamma^+_{i+1}}}{d_{\gamma^+_i}}} + 
	\sigma(\gamma^-)  
		\frac{1}{\sqrt{d_{\gamma^-_{i+1}}d_{\gamma^-_{i+2}}}} 
		\sqrt{\frac{d_{\gamma^-_{i+1}}}{d_{\gamma^-_i}}}
\right) \\
& = 
C 
p_{\widehat{\gamma^\pm}} 
 \frac{1}{\sqrt{d_{\gamma^\pm_1}d_{\gamma^\pm_9}}}
 \sqrt{\prod_{j=1}^{14} \frac{1}{d_{\varphi_j}}}
\left( 
	\frac{\sigma(\gamma^+) }{\sqrt{d_{\gamma^+_{i}}d_{\gamma^+_{i+2}}}} 
	+
	\frac{\sigma(\gamma^-) }{\sqrt{d_{\gamma^-_{i}}d_{\gamma^-_{i+2}}}} 
\right) =0.
\end{align*}
Since $\gamma^+_i = \gamma^-_i = \gamma^+_{i+2} = \gamma^-_{i+2} = \varphi_i$, the two terms in the parentheses cancel exactly. 

Finally, if $\varphi_i$ has valence $3$, then it must be the triple point, $c$. There are then three terms, say $\gamma^0$, $\gamma^1$ and $\gamma^2$, with $\gamma^j_{i+1} = b_j$. Now the collapsed paths differ; $\widehat{\gamma^j} = w_1 j w_2$ for some fixed words $w_1$ and $w_2$. On the other hand, the signs $\sigma(\gamma^j)$ are all equal. Thus we obtain
\begin{align*}
\epsilon_i(S)(\varphi) & = 
C 
\sigma(\gamma^j)  
 \frac{1}{\sqrt{d_{\gamma^j_1}d_{\gamma^j_9}}}
\sqrt{\prod_{j=1}^{14} \frac{1}{d_{\varphi_j}}}
\left(
	\frac{	p_{\widehat{\gamma^0}} }{\sqrt{d_{\gamma^0_{i}}d_{\gamma^0_{i+2}}}} 
+	\frac{	p_{\widehat{\gamma^1}} }{\sqrt{d_{\gamma^1_{i}}d_{\gamma^1_{i+2}}}} 
+	\frac{	p_{\widehat{\gamma^2}} }{\sqrt{d_{\gamma^2_{i}}d_{\gamma^2_{i+2}}}} 
\right).
\end{align*}
Since $\gamma^j_i=\gamma^j_{i+2}= \varphi_i$,  the three terms in the parentheses cancel exactly by Equation \eqref{eq:2s}.
\end{proof}

We now verify the moments of $S$ are the Haagerup moments.

\begin{proof}[Computer-assisted proof of Lemma \ref{lem:moments}.]
We first treat the moments of $S$, and later describe the changes required to calculate the moments of $\rho^{1/2}(S)$.

With multiplication given by the multiplication tangle from Figure \ref{fig:timestracetensor}, the vector space $GPA(\Ha{1})_{8,+}$ becomes a finite-dimensional semisimple associative algebra, which of course must just be a multimatrix algebra. It is easy to see that the simple summands are indexed by pairs of even vertices, and that the minimal idempotents in the summand indexed by $(s,t)$ are given by symmetric loops of length $16$, which go from $s$ to $t$ in $8$ steps, then return the same way. Since there are $8$ even vertices ($v_0, x_0,z_0,b_0,b_1,b_2,z_1$ and $z_2$), there are $64$ simple summands $\mathcal{A}_{s,t}$, although four of these ($\mathcal{A}_{v_0,z_1}, \mathcal{A}_{z_1,v_0}, \mathcal{A}_{v_0,z_2}$, and $\mathcal{A}_{z_2,v_0}$) are trivial because $s$ and $t$ are more than $8$ edges apart. Moreover, the trace tangle from Figure \ref{fig:timestracetensor} composed with the partition function puts a trace on each of these matrix algebras. We write $\mathcal{A}_{s,t} = \left(M_{k\times k}, \frac{d_t}{d_s} \right)$ to indicate there are $k$ paths of length $8$ from $s$ to $t$, and that the trace of the identity in $M_{k \times k}$ is $\frac{d_t}{d_s} k$. We find that
\begin{equation*}
(GPA(\Ha{1})_{8,+}, \text{multiplication tangle}) \iso \DirectSum_{\substack{s,t \\ \text{even vertices}}} \mathcal{A}_{s,t}.
\end{equation*}

Now to compute the required moments, we just need to identify the image of $S$ in this multimatrix algebra, compute the appropriate powers via matrix multiplication, and take weighted traces.
It turns out that the necessary calculation, namely taking $k^{th}$ powers of the matrices for $S$, for $k=2, 3$ and $4$, is actually computationally difficult! First notice that some of the matrices are quite large, up to $118 \times 118$. Worse than this, the entries are quite complicated numbers, involving square roots of dimensions, and so the arithmetic step of simplifying matrix entries after multiplication turns out to be extremely slow. One can presumably do these calculations directly with the help of a computer, using exact arithmetic, but our implementation in Math\-e\-mat\-ic\-a took more than a day attempting to simplify the matrix entries in $S^4$ before we stopped it. Instead, we choose a matrix (really, a multimatrix) $A$ so that all the entries of $A S A^{-1}$ lie in the number field $\Rational(\lambda)$; this matrix certainly has the same moments as $S$, but once the computer can do its arithmetic inside a fixed number field, everything happens much faster. In particular, the moments required here take less than an hour to compute, using Mathematica 7 on a 2.4Ghz Intel Core 2 Duo. See the remark following Definition \ref{defn:S} for an explanation of why this trick works:  we cooked up the matrix $A$ with the desired property by comparing the usual definition of the graph planar algebra with an alternative definition that produces the corresponding ``lopsided'' planar algebra.

The matrix $A$ is defined by
\begin{equation}
\label{eq:A}%
(A_{s,t})_{\pi,\epsilon} = \delta_{\pi=\epsilon} \prod_{i=1}^{8} \sqrt{d_{\pi_i}}
\end{equation}
recalling that the matrix entries in $A_{s,t}$ are indexed by pairs of paths $\pi,\epsilon$ from $s$ to $t$, so $\pi=\pi_1 \cdots \pi_9$ and $\epsilon=\epsilon_1 \cdots \epsilon_9$ with $\pi_1= \epsilon_1 = s$ and $\pi_9=\epsilon_9 =t$. Notice that the index in the product ranges from $1$ to $8$, leaving out the endpoint $t$.

\begin{lem}\label{ASAi}
The entries of  $A S A^{-1}$ lie in $\Rational(\lambda)$.
\end{lem}
(The proof appears below.)

The second half of the Mathematica notebook referred to above produces the matrices for $A S A^{-1}$ (these, and the corresponding matrices for $\rho^{1/2}(S)$ described below, are also available at \url{http://tqft.net/EH/matrices} in machine readable form and as a PDF typeset for an enormous sheet of paper) and actually does the moment calculation. Any reader wanting to check the details should look there. Here, we'll just indicate the schematic calculation:
\begin{align*}
\tr{S^2}(s) & = \sum_{t} \frac{d_t}{d_s} \tr{(S_{s,t})^2} \\
            & = \sum_{t} \frac{d_t}{d_s} \tr{(A_{s,t} S_{s,t} A_{s,t}^{-1})^2} \\
            & \text{approximately 8 minutes later...} \\
            & = [9] \displaybreak[1] \\ 
\tr{S^3}(s) & = \sum_{t} \frac{d_t}{d_s} \tr{(S_{s,t})^3} \\
            & = \sum_{t} \frac{d_t}{d_s} \tr{(A_{s,t} S_{s,t} A_{s,t}^{-1})^3} \\
            & \text{approximately 16 minutes later....} \\
            & = 0  \displaybreak[1] \\ 
\tr{S^4}(s) & = \sum_{t} \frac{d_t}{d_s} \tr{(S_{s,t})^4} \\
            & = \sum_{t} \frac{d_t}{d_s} \tr{(A_{s,t} S_{s,t} A_{s,t}^{-1})^4} \\
            & \text{approximately 24 minutes later.....} \\
            & = [9]. 
\end{align*}
Note that in each case above we're actually computing $8$ potentially different numbers, as $s$ ranges over the even vertices of the graph.

The moments of $\rho^{1/2}(S)$ can be calculated by a very similar approach. The other $8$-box space $GPA(\Ha{1})_{8,-}$ becomes a multimatrix algebra with summands indexed by pairs of odd vertices on the graph $\Ha{1}$.

\begin{lem}\label{ArhohalfSAi}
The entries of $A_{s,t} \rho^{1/2}(S)_{s,t} A_{s,t}^{-1}$ lie in $d \cdot \Rational(\lambda)$.
\end{lem}
(Again, the proof appears below.)

We thus compute 
\begin{align*}
\tr{\rho^{1/2}(S)^3} & = d^3 \tr{(d^{-1} A_{s,t} \rho^{1/2}(S)_{s,t} A_{s,t}^{-1})^3}.
\end{align*}
As before, this is implemented in Mathematica. The calculation takes slightly longer than in the first case. The details can be found in the notebook. 
\end{proof}

\begin{proof}[Proof of Lemma \ref{ASAi}]
Let $\sqcup$ denote concatenation of paths,
and $\bar\epsilon$ be the reverse of the path $\epsilon$.
We readily calculate
\begin{align*}
(A_{s,t} S_{s,t} A_{s,t}^{-1})_{\pi,\epsilon}
	& =   C \sigma_{\pi \sqcup \bar\epsilon} p_{\widehat{\pi \sqcup\bar\epsilon}} \frac{1}{\sqrt{d_s d_t}} \prod_{i=1}^8 \frac{1}{\sqrt{d_{\pi_i}}} \prod_{i=2}^9 \frac{1}{\sqrt{d_{\epsilon_i}}} \frac{\prod_{i=1}^8 \sqrt{d_{\pi_i}}}{\prod_{i=1}^8 \sqrt{d_{\epsilon_i}}} \\
& = \frac{C \sigma_{\pi \sqcup \bar\epsilon} p_{\widehat{\pi \sqcup\bar\epsilon}} }{\prod_{i=1}^{9} d_{\epsilon_i}}
\end{align*}

Most of the factors in this product are already in $\Rational(\lambda)$; 
the one in question is $\prod_{i=1}^{9} d_{\epsilon_i}$.
All even dimensions are in $\Rational[d^2] = \Rational[\lambda^2]$,
and all odd dimensions are in $ d \cdot \Rational[\lambda^2]$. So the product $\prod_{i=1}^{9} d_{\epsilon_i}$,
 a product of five even and four odd dimensions,
lies in $d^4 \Rational[d^2] \subset \mathbb{Q}(\lambda)$
and
$$(A S A^{-1})_{\gamma,\epsilon} \in \mathbb{Q}(\lambda).$$ 
\end{proof}

\begin{proof}[Proof of Lemma \ref{ArhohalfSAi}]
First, we have 
\begin{align*}
\rho^{1/2}(S)(\gamma_1\gamma_2 \cdots \gamma_{16}\gamma_1) 
 & = \sqrt{\frac{d_{\gamma_2} d_{\gamma_{10}}}{d_{\gamma_9}d_{\gamma_1}} } 
 S(\gamma_2 \cdots \gamma_{16}\gamma_{1} \gamma_2) \\
 & =   \sqrt{\frac{d_{\gamma_2} d_{\gamma_{10}}}{d_{\gamma_9}d_{\gamma_1}} } 
  C \sigma(\gamma) 
 p_{\widehat{\gamma_2 \cdots \gamma_{16}\gamma_{1} \gamma_2}}
 \frac{1}{\sqrt{d_{\gamma_2}d_{\gamma_10}}} 
 \prod_{i=1}^{16} \frac{1}{\sqrt{d_{\gamma_i}}} \\
 & =  C \sigma(\gamma) 
 p_{\widehat{\gamma_2 \cdots \gamma_{16}\gamma_{1} \gamma_2 }} 
 \frac{1}{\sqrt{d_{\gamma_1}d_{\gamma_9}}} \prod_{i=1}^{16} \frac{1}{\sqrt{d_{\gamma_i}}}.
\end{align*}
Be careful here: although this looks very similar to the formula in Equation \eqref{eq:S} for $S$, the path $\gamma$ here starts at an odd vertex.

We now conjugate by a multimatrix $A$ that has exactly the same formula for its definition as appears in Equation \eqref{eq:A}, except again the paths $\gamma$ and $\epsilon$ start and finish at odd vertices. We obtain 
\begin{align*}
(A_{s,t} \rho^{1/2}(S)_{s,t} A_{s,t}^{-1})_{\pi,\epsilon}
& = \frac{r \sigma_{\pi \sqcup \bar\epsilon} p_{\widehat{\epsilon_2 \pi_1 \cdots \pi_8 \epsilon_9 \cdots \epsilon_2}} }{\prod_{i=1}^{9} d_{\epsilon_i}}.
\end{align*}
One readily checks that these matrix entries are in $d \cdot \Rational(\lambda)$.
\end{proof}

We've finally shown the existence and uniqueness of the extended Haagerup subfactor. Uniqueness is Theorem \ref{thm:uniqueness}.  By Lemma \ref{lem:Srelations} and Lemma \ref{lem:moments}, $S$ satisfies the hypotheses of Proposition \ref{prop:SgivesSPA}.  Therefore $\PA(S)$ is a subfactor planar algebra with principal graphs $\Ha{1}$.

%% file: text/polynomials.tex
%% This file is automatically generated by Mathematica, using the WordToDualLoop notebook.
p_{00000001} & = -2 \lambda ^4-\lambda ^2+9 &
p_{00000011} & = -\lambda ^5-\lambda ^3+3 \lambda  \displaybreak[1] \\
p_{00000101} & = 2 \lambda ^4+\lambda ^2-9 &
p_{00000111} & = 1 \displaybreak[1] \\
p_{00001001} & = \lambda ^5-\lambda ^3-3 \lambda  &
p_{00001011} & = \lambda ^3-1 \displaybreak[1] \\
p_{00010001} & = 2 \lambda ^4+\lambda ^2-9 & 
p_{00010011} & = \lambda ^5-\lambda ^4+\lambda ^2-3 \lambda +4 \displaybreak[1] \\
p_{00010101} & = \lambda ^4-2 \lambda ^2+1 &
p_{00010111} & = -\lambda ^4+1 \displaybreak[1] \\
p_{00011011} & = \lambda ^4-\lambda ^2-3 &
p_{00100101} & = -2 \lambda ^4+5 \displaybreak[1] \\
p_{00100111} & = \lambda ^2+1 &
p_{00101011} & = -\lambda ^5-\lambda ^3+\lambda +1 \displaybreak[1] \\
p_{00101101} & = \lambda ^5-\lambda  &
p_{00110011} & = 2 \lambda ^5+5 \lambda ^3+4 \lambda  \displaybreak[1] \\
p_{00110111} & = -\lambda ^5-2 \lambda ^3-4 \lambda ^2-\lambda -5 &
p_{01010101} & = -4 \lambda ^4+3 \lambda ^2+7 \displaybreak[1] \\
p_{01010111} & = \lambda ^4+\lambda ^2 &
p_{01011011} & = \lambda ^4-2 \lambda ^2-4 \displaybreak[1] \\
p_{01110111} & = \lambda ^4+6 \lambda ^2+6 & &

%% file: text/fusion.tex
\section{Fusion categories coming from the extended Haagerup subfactor}
\label{sec:fusion}
The even parts of a subfactor are the unitary tensor categories of $N-N$ and $M-M$ bimodules respectively.  Hence every finite depth subfactor yields two unitary fusion categories.  In terms of the planar algebra, the simple objects in these categories are the irreducible projections in the box spaces $P_{2m,\pm}$ for some $m$.

In the case of extended Haagerup, the global dimension of each of these fusion categories is the largest real root of $x^3 -585 x^2 +8450 x -21125$ (approximately 570.247).  The fusion tables are given in Figures \ref{fig:NNfusion} and \ref{fig:MMfusion}.  

\newcommand{\p}{\!+\!}

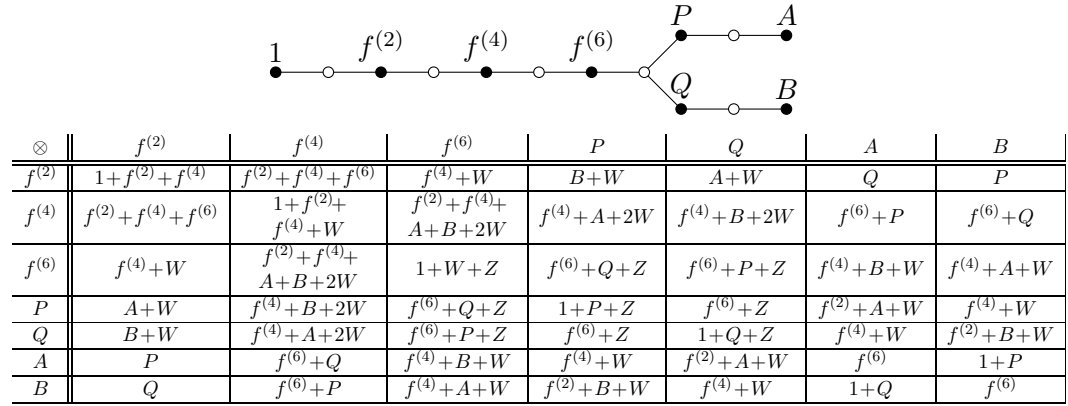
\begin{figure}[!ht]
$$\begin{tikzpicture}[baseline=-1mm, scale=.7]
	\draw (0,0) -- (7,0);
	\draw (7,0)--(7.7,.7)--(9.7,.7);
	\draw (7,0)--(7.7,-.7)--(9.7,-.7);

	\filldraw              (0,0) node [above] {$1$} circle (1mm);
	\filldraw [fill=white] (1,0) node {} circle (1mm);
	\filldraw              (2,0) node [above] {$\JW{2}$}  circle (1mm);
	\filldraw [fill=white] (3,0) node {}  circle (1mm);
	\filldraw              (4,0) node [above] {$\JW{4}$}  circle (1mm);
	\filldraw [fill=white] (5,0) node {}  circle (1mm);
	\filldraw              (6,0) node [above] {$\JW{6}$}  circle (1mm);
	\filldraw [fill=white] (7,0) node {}  circle (1mm);

	\filldraw              (7.7,.7) node [above] {$P$}  circle (1mm);
	\filldraw [fill=white] (8.7,.7) node [above] {}  circle (1mm);
	\filldraw              (9.7,.7) node [above] {$A$}  circle (1mm);
	\filldraw              (7.7,-.7) node [above] {$Q$}  circle (1mm);
	\filldraw [fill=white] (8.7,-.7) node [above] {}  circle (1mm);
	\filldraw              (9.7,-.7) node [above] {$B$}  circle (1mm);
\end{tikzpicture}$$

%\begin{equation*}
%\scriptsize
%\begin{array}{cccccccc}
% 1 & \JW{2} & \JW{4} & \JW{6} & P & Q & A & B \\
% \JW{2} & \JW{2}+\JW{4}+1 & \JW{2}+\JW{4}+\JW{6} & \JW{4}+W & A+W & B+W & P & Q \\
% \JW{4} & \JW{2}+\JW{4}+\JW{6} & \JW{2}+\JW{4}+W+1 & A+B+\JW{2}+\JW{4}+2 W & B+\JW{4}+2 W & A+\JW{4}+2 W & \JW{6}+Q & \JW{6}+P \\
% \JW{6} & \JW{4}+W & A+B+\JW{2}+\JW{4}+2 W & \JW{6}+P+Q+Z+1 & \JW{6}+Q+Z & \JW{6}+P+Z & B+\JW{4}+W & A+\JW{4}+W \\
% P & A+W & B+\JW{4}+2 W & \JW{6}+Q+Z & P+Z+1 & \JW{6}+Z & A+\JW{2}+W & \JW{4}+W \\
% Q & B+W & A+\JW{4}+2 W & \JW{6}+P+Z & \JW{6}+Z & Q+Z+1 & \JW{4}+W & B+\JW{2}+W \\
% A & P & \JW{6}+Q & B+\JW{4}+W & A+\JW{2}+W & \JW{4}+W & P+1 & \JW{6} \\
% B & Q & \JW{6}+P & A+\JW{4}+W & \JW{4}+W & B+\JW{2}+W & \JW{6} & Q+1
%\end{array}
%\end{equation*}

\scalebox{0.7}{
\begin{tabular}{c||c|c|c|c|c|c|c|}
$\tensor$  & $\JW{2}$ & $\JW{4}$ & $\JW{6}$ & $P$ & $Q$ & $A$ & $B$ \\
\hline\hline
 $\JW{2}$ & $1\p \JW{2}\p\JW{4}$ & $\JW{2}\p\JW{4}\p\JW{6}$ & $\JW{4} \p W$ & $B \p W$ & $A\p W$ & $Q$ & $P$ \\
 \hline
 $\JW{4}$ & $\JW{2} \p \JW{4} \p \JW{6}$ & $\begin{array}{c}1 \p \JW{2} \p \\ \JW{4} \p W \end{array}$ & $\begin{array}{c} \JW{2} \p \JW{4} \p \\ A \p B \p 2 W \end{array}$ & $\JW{4} \p A \p 2 W$ & $\JW{4} \p B \p 2 W$ & $\JW{6} \p P$ & $\JW{6} \p Q$ \\
 \hline
 $\JW{6}$ & $\JW{4} \p W$ & $\begin{array}{c}\JW{2} \p \JW{4} \p \\ A \p B \p 2 W \end{array}$ & $1 \p W \p Z$ & $\JW{6} \p Q \p Z$ & $\JW{6} \p P \p Z$ & $\JW{4} \p B \p W$ & $\JW{4} \p A \p W$ \\
 \hline
 $P$ & $A \p W$ & $\JW{4} \p B \p 2 W$ & $\JW{6} \p Q \p Z$ & $1 \p P \p Z$ & $\JW{6} \p Z$ & $\JW{2} \p A \p W$ & $\JW{4} \p W$ \\
 \hline
 $Q$ & $B \p W$ & $ \JW{4} \p A \p 2 W$ & $\JW{6} \p P \p Z$ & $\JW{6} \p Z$ & $1 \p Q \p Z$ & $\JW{4} \p W$ & $\JW{2} \p B \p W$ \\
 \hline
 $A$ & $P$ & $\JW{6} \p Q$ & $ \JW{4} \p B \p W$ & 
  $\JW{4} \p W$ & $\JW{2} \p A \p W$ &
  $\JW{6}$ & $1 \p P$ \\
 \hline
 $B$ & $Q$ & $\JW{6} \p P$ & $\JW{4} \p A \p W$ & 
  $\JW{2} \p B \p W$ & $\JW{4} \p W$ &
 $1 \p Q$ & $\JW{6}$ \\ \hline
\end{tabular}
}
\caption{The simple objects and fusion rules for the $N-N$ fusion category coming from the extended Haagerup subfactor. We use the abbreviations $W  = \JW{6}+P+Q$ and 
$Z  =A + B + \JW{2} + 2 \JW{4} + 3 \JW{6} + 3 P + 3 Q$. }
\label{fig:NNfusion}
\end{figure}

\begin{figure}
$$\begin{tikzpicture}[baseline=-1mm, scale=.7]
	\draw (0,0) -- (7,0);
	\draw (7,0)--(7.7,.7);
	\draw (7,0)--(7.7,-.7);
	\draw (7.7,.7)--(8.4,1.4);
	\draw (7.7,.7)--(8.4,0);

	\filldraw              (0,0) node [above] {$1$} circle (1mm);
	\filldraw [fill=white] (1,0) node {} circle (1mm);
	\filldraw              (2,0) node [above] {$\JW{2}$}  circle (1mm);
	\filldraw [fill=white] (3,0) node {}  circle (1mm);
	\filldraw              (4,0) node [above] {$\JW{4}$}  circle (1mm);
	\filldraw [fill=white] (5,0) node {}  circle (1mm);
	\filldraw              (6,0) node [above] {$\JW{6}$}  circle (1mm);
	\filldraw [fill=white] (7,0) node {}  circle (1mm);

	\filldraw              (7.7,.7) node [above] {$P'$}  circle (1mm);
	\filldraw              (7.7,-.7) node [below] {$Q'$}  circle (1mm);
	\filldraw [fill=white] (8.4,1.4) node [above] {}  circle (1mm);
	\filldraw [fill=white] (8.4,0) node [above] {}  circle (1mm);
\end{tikzpicture}$$

\scalebox{.73}{
\begin{tabular}{c||c|c|c|c|c|}
 $\otimes$ & $\JW{2}$ & $\JW{4}$ & $\JW{6}$ & $P'$ & $Q'$ \\ \hline \hline
 $\JW{2}  $ & $1 \p \JW{2} \p \JW{4}  $ & $ \JW{2} \p \JW{4} \p \JW{6}  $ & $ \JW{4} \p \JW{6} \p P' \p Q'  $ & $ \JW{6} \p 2 P' \p Q'  $ & $ \JW{6} \p P' $\\ \hline
$ \JW{4}  $ & $ \JW{2} \p \JW{4} \p \JW{6}  $ & $\begin{array}{c} 1 \p \JW{2} \p \JW{4} \\  \p \JW{6} \p P' \p Q' \end{array}$ & $\begin{array}{c} \JW{2} \p \JW{4} \p 2 \JW{6} \\ \p 3 P' \p Q'  \end{array}$ & $ \begin{array}{c} \JW{4} \p 3 \JW{6} \\ \p 3 P' \p 2 Q'  \end{array}$ & $\begin{array}{c} \JW{4} \p \JW{6} \\ \p 2 P' \p Q' \end{array}$ \\ \hline
$\JW{6}$ & $\JW{4} \p \JW{6} \p P' \p Q'  $ & $\begin{array}{c} \JW{2} \p \JW{4} 2 \JW{6} \\ \p 3 P' \p Q'  \end{array}$ & $\begin{array}{c} 1 \p\JW{2} \p 2 \JW{4}  \\ \p 4 \JW{6} \p 5 P' \p 3 Q'  \end{array}$ & $\begin{array}{c} \JW{2} \p 3 \JW{4} \p 5 \JW{6}  \\ \p 6 P' \p 3 Q' \end{array} $ & $\begin{array}{c} \JW{2} \p \JW{4} \p 3 \JW{6} \\ \p 3 P' \p 2 Q' \end{array}$ \\ \hline
$P'$ & $ \JW{6} \p 2 P' \p Q' $ & $\begin{array}{c} \JW{4} \p 3 \JW{6} \\ \p 3 P' \p 2 Q'\end{array}$ & $\begin{array}{c} \JW{2} \p 3 \JW{4} \p 5 \JW{6} \\ \p 6 P' \p 3 Q'  \end{array}$ & $\begin{array}{c} 1 \p 2 \JW{2} \p 3 \JW{4} \\ \p 6 \JW{6} \p 7 P' \p 4 Q'   \end{array}$ & $ \begin{array}{c} \JW{2} \p 2 \JW{4} \p 3 \JW{6} \\ \p 4 P' \p 2 Q' \end{array}$\\ \hline
 $Q'  $ & $ \JW{6} \p P'  $ & $\begin{array}{c} \JW{4} \p \JW{6} \\ \p 2 P' \p Q'  \end{array}$ & $\begin{array}{c} \JW{2} \p \JW{4} \p 3 \JW{6} \\ \p 3 P' \p 2 Q' \end{array}$ & $\begin{array}{c} \JW{2} \p 2 \JW{4} \p 3 \JW{6} \\ \p 4 P' \p 2 Q' \end{array}$ & $\begin{array}{c}1 \p \JW{4} \p 2 \JW{6}  \\ \p 2 P' \p Q' \end{array}$ \\ \hline
\end{tabular}}
\caption{The simple objects and fusion rules for the $M-M$ fusion category coming from the extended Haagerup subfactor.}
\label{fig:MMfusion}
\end{figure}

%\begin{equation*}
%\scriptsize
%\begin{array}{cccccc}
% 1 & \JW{2} & \JW{4} & \JW{6} & P' & Q' \\
% \JW{2} & \JW{2}+\JW{4}+1 & \JW{2}+\JW{4}+\JW{6} & \JW{4}+\JW{6}+P'+Q' & \JW{6}+2 P'+Q' & \JW{6}+P' \\
% \JW{4} & \JW{2}+\JW{4}+\JW{6} & \JW{2}+\JW{4}+\JW{6}+P'+Q'+1 & \JW{2}+\JW{4}+2 \JW{6}+3 P'+Q' & \JW{4}+3 \JW{6}+3 P'+2 Q' & \JW{4}+\JW{6}+2 P'+Q' \\
% \JW{6} & \JW{4}+\JW{6}+P'+Q' & \JW{2}+\JW{4}+2 \JW{6}+3 P'+Q' & \JW{2}+2 \JW{4}+4 \JW{6}+5 P'+3 Q'+1 & \JW{2}+3 \JW{4}+5 \JW{6}+6 P'+3 Q' & \JW{2}+\JW{4}+3 \JW{6}+3 P'+2 Q' \\
% P' & \JW{6}+2 P'+Q' & \JW{4}+3 \JW{6}+3 P'+2 Q' & \JW{2}+3 \JW{4}+5 \JW{6}+6 P'+3 Q' & 2 \JW{2}+3 \JW{4}+6 \JW{6}+7 P'+4 Q'+1 & \JW{2}+2 \JW{4}+3 \JW{6}+4 P'+2 Q' \\
 %Q' & \JW{6}+P' & \JW{4}+\JW{6}+2 P'+Q' & \JW{2}+\JW{4}+3 \JW{6}+3 P'+2 Q' & \JW{2}+2 \JW{4}+3 \JW{6}+4 P'+2 Q' & \JW{4}+2 \JW{6}+2 P'+Q'+1
%\end{array}
%\end{equation*}